\documentclass[11pt,a4paper]{article}

\usepackage[a4paper,margin=30mm]{geometry}
\usepackage[T1]{fontenc}
\usepackage[utf8]{inputenc}
\usepackage{mathpazo}
\usepackage{microtype}
\usepackage{amsmath,amssymb,amsthm,mathtools}
\usepackage{enumitem}
\usepackage{xcolor}
\usepackage{hyperref}
\hypersetup{
  colorlinks=true,
  linkcolor=blue!55!black,
  citecolor=blue!55!black,
  urlcolor=blue!55!black,
  pdftitle={The singular Paneitz equation on four-manifolds: refined blow-up analysis and degree jumps},
  pdfauthor={Mohameden Ahmedou, Lina Wu, Lei Zhang},
  pdfsubject={Singular Paneitz equation, refined blow-up analysis and nonresonant degree}
}
\numberwithin{equation}{section}
\newtheorem{theorem}{Theorem}[section]
\newtheorem{proposition}[theorem]{Proposition}
\newtheorem{lemma}[theorem]{Lemma}
\newtheorem{corollary}[theorem]{Corollary}
\theoremstyle{definition}
\newtheorem{definition}[theorem]{Definition}

\newcommand{\R}{\mathbb R}
\newcommand{\Sph}{\mathbb S}
\newcommand{\cD}{\mathcal D}
\newcommand{\cF}{\mathcal F}
\newcommand{\cK}{\mathcal K}
\newcommand{\Crit}{\operatorname{Crit}}

\newcommand{\cT}{\mathcal T}
\newcommand{\Conf}{\operatorname{Conf}}
\newcommand{\morse}{\operatorname{morse}}
\newcommand{\dd}{\,dV_g}
\newcommand{\grad}{\nabla}
\newcommand{\supp}{\operatorname{supp}}
\newcommand{\Ent}{\operatorname{Ent}}

\newcommand{\shorttitle}{The singular Paneitz equation on four-manifolds}
\title{\bfseries \shorttitle:\\[0.4em]
\large refined blow-up analysis and degree jumps}
\author{
Mohameden Ahmedou\thanks{Department of Mathematics, Justus Liebig University Giessen, Germany.}
\and
Lina Wu\thanks{Department of Mathematics, Beijing Jiaotong University, China.}
\and
Lei Zhang\thanks{Department of Mathematics, University of Florida, USA.}}
\date{20 November 2026}

\begin{document}

\maketitle
\pagestyle{myheadings}
\markboth{\shorttitle}{\shorttitle}

\begin{abstract}
We compute the Leray--Schauder degree of the singular Paneitz equation on a
closed four-manifold with conical orders strictly between minus one and
zero. We assume that the kernel of the Paneitz operator consists of the
constants. In each positive nonresonant chamber the degree is a partial sum
of an explicit generating series determined by the Euler characteristic and
the conical orders. We match the heights of exact blow-up solutions and
expand the gradient near the corrected projected bubbles. A localized
dilation identity gives its common-scale component on the zero set of the
relative-scale equations. For generic data a fixed-mass homotopy then
gives the signed degree of each concentrating configuration. Assign ordinary points weight one and conical points weight one plus
their order. If the largest concentrating weight exceeds one half, a finite
second moment gives the coefficient. At one half the expansion acquires a
logarithm, and below one half the coefficient is a global finite part. On
the round four-sphere with two antipodal points of equal weight below one
half, perturbing the prescribed curvature away from both points selects the
side of bubbling.
\end{abstract}



{\setcounter{tocdepth}{1}\tableofcontents}

\section{The geometric equation and the degree problem}
For the prescribed $Q$-curvature equation on a four-manifold with conical
points, quantization identifies the masses at which compactness can fail.
Several configurations of ordinary and conical points can carry the same
mass. For negative conical orders and generic curvature data, we determine
the side of concentration and the signed contribution of each configuration
to the degree jump. Compactness then gives the degree in every chamber for
all smooth positive data. Perturbing the curvature away from the blow-up
points can change the side of concentration.

Let $(M,g)$ be smooth, closed and connected, of volume $V$. Fix a
smooth positive function $h$ and distinct points
$\Sigma=\{q_1,\ldots,q_N\}$ with orders $-1<\gamma_j<0$. Write
$D=\sum_j\gamma_j\delta_{q_j}$ and assume $\ker P_g=\R$.
We consider the Paneitz equation with compensated point sources
\begin{equation}\label{eq:original}
 P_gv_P+2Q_g
 =2he^{4v_P}-8\pi^2\sum_j\gamma_j(\delta_{q_j}-V^{-1}).
\end{equation}
Integration fixes the mass $\int_M he^{4v_P}\dd$ at
$\kappa_P=\int_M Q_g\dd$. Removing the logarithmic singularities
gives a positive weight $K_D$ whose order at $q_j$ is
$d_g(x,q_j)^{4\gamma_j}$. We let the mass vary and study the
mean-zero equation
\begin{equation}\label{eq:normalized}
 P_gu=2\kappa(\mu_u-V^{-1}dV_g),\qquad
 \mu_u=\frac{K_De^{4u}\,dV_g}{\int_M K_De^{4u}\dd},\qquad
 \int_Mu\dd=0 .
\end{equation}
In distributional equations we identify functions with their volume
densities. Section~\ref{sec:green} gives the weight $K_D$ and the
bijection between solutions of \eqref{eq:normalized} at
$\kappa=\kappa_P>0$ and solutions of \eqref{eq:original}.

Put $\rho=\kappa/(8\pi^2)$ and $\mu_j=1+\gamma_j$. Quantization
\cite[Theorem~5.1 and Corollaries~5.2--5.3]{AWZ} confines loss of
compactness to the positive parameters in
\begin{equation}
 \mathcal W=\left\{p+\sum_{j\in I}\mu_j:
 p\in\mathbb N_0,\ I\subset\{1,\ldots,N\}\right\},\qquad
 \cK(M,D)=8\pi^2(\mathcal W\setminus\{0\}).
\end{equation}
We call a parameter in $\cK(M,D)$ a wall and a component of its
complement in $(0,\infty)$ a chamber. For a wall $\kappa_*$, a
concentration type $(p,I)$ satisfies
\begin{equation}
 \kappa_*=8\pi^2\left(p+\sum_{j\in I}\mu_j\right).
\end{equation}
In the source measure $2\kappa\mu_u$, the quantization theorem allows
an ordinary atom of mass in $16\pi^2\mathbb N$. A localized dilation
identity (Lemma~\ref{lem:virial}) gives mass exactly $16\pi^2$ at
each ordinary point, so $p$ counts distinct ordinary points. At an
active mark $q_j$, $j\in I$, the atom has mass $16\pi^2\mu_j$ in the
limiting source measure and mass $8\pi^2\mu_j/\kappa_*$ in the limiting
probability measure.

\begin{theorem}[Degree in the nonresonant chambers]\label{thm:main-degree}
For $\kappa>0$ outside $\cK(M,D)$, let $d_{\rm fp}(\kappa)$ denote the
Leray--Schauder degree of the fixed-point equation obtained by applying
$P_g^{-1}$ to \eqref{eq:normalized} on the volume-mean-zero
subspace of $H^2(M)$, with the identity orientation. Then
\begin{equation}\label{eq:generating}
 d_{\rm fp}(\kappa)=\sum_{\substack{\sigma\in\mathcal W\\ \sigma<\rho}} b_\sigma,
 \qquad
 \sum_{\sigma\in\mathcal W} b_\sigma x^\sigma
 =(1-x)^{\chi(M)-N}\prod_{j=1}^N(1-x^{\mu_j}).
\end{equation}
If $d_{\rm fp}(\kappa)\ne0$, equation \eqref{eq:normalized} has a solution.
\end{theorem}

The factor $(1-x)^{\chi(M)-N}$ records the Euler characteristics of
configurations of ordinary points in $M\setminus\Sigma$. Each marked
point can occur at most once, which accounts for its factor
$1-x^{\mu_j}$. The degree starts at one and $b_\sigma$ is its jump
across $8\pi^2\sigma$ for each positive $\sigma\in\mathcal W$.
The geometric parameter for prescribing $Q$-curvature off the marks is
$\kappa_D=\kappa_P+4\pi^2\sum_j\gamma_j$ by the singular
Chern--Gauss--Bonnet formula \cite[Proposition~3.2]{FangMa}.
The following consequence gives explicit existence criteria for both
curvature masses.

\begin{corollary}[Nonzero degree and existence]\label{cor:nonzero-chambers}
Take $\kappa=\kappa_D$ or $\kappa=\kappa_P$ and put $\rho=\kappa/(8\pi^2)$.
Order the marked weights so that
$\mu_{(1)}\le\cdots\le\mu_{(N)}$, with empty sums taken as zero.
\begin{enumerate}[label=\textup{(\roman*)},leftmargin=*]
\item If for some integer $k\ge0$,
\[
 k<\rho< k+1-\sum_{j=1}^{\min(k+1,N)}(1-\mu_{(j)}),
\]
then $d_{\rm fp}(\kappa)=(-1)^k\binom{\chi(M)-1}{k}\ne0$.
\item Suppose $0<\rho<1$, $\rho<\mu_i+\mu_j$ for every $i<j$,
and $\rho\ne\mu_j$ for every $j$. If $r=\#\{j:\mu_j<\rho\}$,
then $d_{\rm fp}(\kappa)=1-r$.
\end{enumerate}
In case \textup{(i)}, and in case \textup{(ii)} when $r\ne1$,
\eqref{eq:normalized} has a solution for every smooth positive $h$.
At $\kappa=\kappa_D$ it yields a conformal metric with divisor $D$
and $Q$-curvature $h$ on $M\setminus\Sigma$. At $\kappa=\kappa_P$
it yields a solution of \eqref{eq:original}.
\end{corollary}

For example, if $N>0$ and exactly $m$ marks have the smallest weight, the degree
in the chamber immediately above the first wall is $1-m$.

\smallskip
\noindent\textit{The mass defect.}
Give each concentrating point $z_s$ weight $\mu_s$, with $\mu_s=1$
at an ordinary point, and put $m=\max_s\mu_s$. Theorem~\ref{thm:strongentry}
matches the heights of an exact blow-up sequence: if $\lambda_s$ is
the inverse length scale of its $s$-th peak, then $\lambda_s^{\mu_s}$
is comparable to a common scale $\Lambda\to\infty$. With the
normalization $\Lambda=e^{T_0}$ of Section~\ref{sec:rate-low} and limiting
configuration $A$, Theorem~\ref{thm:allrates} gives
\[
 \frac{\kappa-\kappa_*}{\kappa_*}
 = C(A)\,s_m(\Lambda)+o\bigl(s_m(\Lambda)\bigr),\qquad
 s_m(\Lambda)=
 \begin{cases}
 \Lambda^{-2/m}, &m>1/2,\\
 \Lambda^{-4}\log\Lambda, &m=1/2,\\
 \Lambda^{-4}, &m<1/2 .
 \end{cases}
\]
The transition comes from the radial second moment, whose tail has
density proportional to $r^{1-4\mu_s}\,dr$. In the first two cases
the coefficient uses the regular part of the outer field and the metric
at the points of largest weight. In the third case all concentrating
points are marked, and $C(A)=D_I^{(4)}/\kappa_*$ is the finite part of a
global integral. Multiplying $h$ by $e^v$, with $v\ge0$ smooth,
nonzero and supported away from the marks, strictly increases this
finite part, since every counterterm in \eqref{eq:D4} stays fixed.
This determines the side of concentration when the coefficient is nonzero.
On the round sphere with two antipodal marks of equal weight below one
half, constant curvature data admit bubbling at fixed mass. A perturbation
away from the marks moves bubbling to either prescribed side while
preserving the curvature data near both poles, by \eqref{eq:antipodal-variation}.

\smallskip
\noindent\textit{From a configuration to its degree.}
The singular bubble has a single bounded Jacobi field, generated by
dilation (Theorem~\ref{thm:kernel}). Projecting the first angular
correction absorbs the linear jet of the outer field and gives a normal
estimate uniform in the component scales
(Propositions~\ref{prop:corrected-defect} and~\ref{prop:balancedgap}).
Theorem~\ref{thm:strongentry} places each exact blowing-up solution in this
chart with $(\alpha_s-1)\sqrt{\log\lambda_s}\to0$.
On the shrinking neighborhood of Proposition~\ref{prop:sharp-fixed},
put $\Lambda=e^T$. If $|\kappa-\kappa_*|=O(s_m(\Lambda))$ and the
relative-scale equations vanish, the common-scale component equals
$C(A)s_m(\Lambda)-(\kappa-\kappa_*)/\kappa_*+o(s_m(\Lambda))$. The homotopy
\eqref{eq:fixed-wall-homotopy} retains the relative-scale equations
until the common-scale component has reached its model value.
For generic data an unordered critical configuration $A$ of the
interaction function $\cF_{p,I}$ contributes
$(-1)^{p+|I|+\morse(\cF_{p,I},A)}$ to the fixed-point degree jump
(Theorem~\ref{thm:jump}).

In dimension two, Chen--Lin \cite{ChenLinSharp,ChenLinDegree} obtained
sharp blow-up estimates and a generating-series formula for the singular
mean field equation. At fixed resonance, Ahmedou--Ben Ayed--Lucia
\cite{ABL} studied the Dirichlet problem by critical points at infinity,
and Ahmedou--Hu--Wang--Zhang \cite{AHWZ} obtained a degree formula for
the Neumann problem. For the fourth-order Navier problem on bounded domains in
$\R^4$, Lin--Wei \cite{LinWei} derived the leading mass defect at
configurations of ordinary blow-up points.
For smooth closed four-manifolds, Malchiodi \cite[Theorem~1.3]{MalchiodiDegree}
computed the degree at nonresonant total $Q$-curvature using the
variational orientation. Section~\ref{sec:degree} gives the factor
relating that orientation to the fixed-point orientation of
Theorem~\ref{thm:main-degree}: each negative Paneitz eigenvalue,
counted with multiplicity, contributes a factor $-1$. Ndiaye
\cite{NdiayeResonance} studied the degree at nondegenerate resonance
in even dimensions. For positive conical orders, Jevnikar--Sire--Yang
\cite[Theorem~1.1]{JSY} proved supercritical existence under a
topological retraction hypothesis. Here the sharp mass
defect for mixed configurations and the local degree at each end connect
our conical quantization theorem \cite{AWZ} to the chamber formula.

An independent variational argument gives an existence criterion in
terms of the topology of the concentration measures. Give an ordinary point cost $\omega(x)=1$ and
a marked point cost $\omega(q_j)=\mu_j$. Let $\omega_0=\min_{x\in M}\omega(x)$.
For $\rho>0$ the weighted barycenter space is
\[
 B_\rho=\left\{\sum_{i=1}^k t_i\delta_{x_i}:
 k\ge1,\ t_i>0,\ \sum_i t_i=1,\ x_i\text{ distinct},\
 \sum_i\omega(x_i)<\rho\right\},
\]
with the weak topology of probability measures.

\begin{theorem}[Nonresonant existence]\label{th:exist}
Suppose $\rho>0$ and $\rho\notin\mathcal W$. Equation \eqref{eq:normalized}
has a solution if $\rho<\omega_0$ or if
$\widetilde H_j(B_\rho;\mathbb F_2)\ne0$ for some $j\ge0$.
If $P_g\ge0$ and $\rho<\omega_0$, a solution minimizes the mean-zero
energy $J_\kappa$ in \eqref{eq:energy}. These conclusions apply to \eqref{eq:original} when
$\kappa=\kappa_P>0$.
\end{theorem}

\smallskip
\noindent\textit{Acknowledgements.} Lina Wu is partially supported by
the China Scholarship Council (No.~201806210165). Lei Zhang is partially
supported by a Simons Foundation Collaboration Grant.

\smallskip
\noindent\textit{Use of artificial intelligence.}
The authors used OpenAI's ChatGPT (Codex) to check calculations and
references and to edit the exposition and \LaTeX{} source.

\section{Normalization and Green profiles}\label{sec:green}
Our sign convention is $\Delta_g=\operatorname{div}_g\grad$, and $\delta$
is the formal adjoint of $d$ on functions. We use
\[
 P=P_g=\Delta_g^2+\delta\big((\tfrac23 R_g g-2\operatorname{Ric}_g)d\,\cdot\,\big),
 \qquad Q_g=-\tfrac1{12}\big(\Delta_gR_g-R_g^2+3|\operatorname{Ric}_g|^2\big).
\]
Thus $P_gv+2Q_g=2Q_{e^{2v}g}e^{4v}$ for a smooth conformal factor.
Let $G_0$ be the mean-zero Green function of $P_g$ and $f$ the mean-zero
solution of $P_gf=2\kappa_P/V-2Q_g$. Define
\begin{equation}\label{eq:weight}
 K_D=h\exp\Big(4f-32\pi^2\sum_j\gamma_jG_0(\cdot,q_j)\Big),\quad
 dm_D=\frac{K_D\,dV_g}{\int_M K_D\dd},\quad Z(u)=\int_M e^{4u}\,dm_D .
\end{equation}
Then $\mu_u=e^{4u}dm_D/Z(u)$. At $\kappa=\kappa_P>0$ the lift is
\begin{equation}\label{eq:lift}
 v_P=u+f-8\pi^2\sum_j\gamma_jG_0(\cdot,q_j)
     +\tfrac14\log\frac{\kappa_P}{\int_M K_De^{4u}\dd}.
\end{equation}
The constant imposes $\int_Mhe^{4v_P}\dd=\kappa_P$, and applying
$P_g$ gives \eqref{eq:original}. Conversely, remove the fixed singular
and smooth terms from a solution of \eqref{eq:original} and subtract the
mean to recover $u$.

In our normalization the smooth Chern--Gauss--Bonnet formula reads
$\kappa_P=4\pi^2\chi(M)-\tfrac18\int_M|W_g|^2\dd$, where $W_g$
is the Weyl tensor. To obtain the lift at $\kappa_D$, take a solution of
\eqref{eq:normalized} at any $\kappa>0$, replace the numerator
$\kappa_P$ in the logarithm of \eqref{eq:lift} by $\kappa$ and
denote the resulting lift by $v_\kappa$. The equations for $f$ and $G_0$
give
\[
 P_gv_\kappa+2Q_g
 =2he^{4v_\kappa}-8\pi^2\sum_j\gamma_j\delta_{q_j}
       +2(\kappa_D-\kappa)/V .
\]
At $\kappa=\kappa_D>0$ the lift has the prescribed curvature off the
marks. Its singular term at $q_j$ is $\gamma_j\log d_g(x,q_j)$.
The regularity of $u$ and the Green remainder below give an $H^2$
remainder, smooth away from the marks. In the total-curvature formula
of \cite[Proposition~3.2]{FangMa}, $Q$ is twice our $Q_g$.

For any $\kappa>0$, use the mass normalization
\begin{equation}\label{eq:mass-gauge}
 v=u-\tfrac14\log\int_M K_De^{4u}\dd,\qquad \int_M K_De^{4v}\dd=1 .
\end{equation}
For a solution it gives
\begin{equation}
 P_gv=2\kappa\big(K_De^{4v}-V^{-1}\big).
\end{equation}
Since $u$ has mean zero,
\begin{equation}
 \tau=e^{4\bar v}=\left(\int_M K_De^{4u}\dd\right)^{-1},
 \qquad \bar v=V^{-1}\int_M v\dd .
\end{equation}
For the blow-up analysis we use $G=16\pi^2G_0$, so
\begin{equation}\label{eq:green}
 P_gG(a,\cdot)=16\pi^2(\delta_a-V^{-1}),\qquad \int_MG(a,x)\,dV_g(x)=0 .
\end{equation}
In geodesic normal coordinates and near the diagonal,
\begin{equation}\label{eq:green-split}
                   G(a,x)=2\log\frac1{d_g(a,x)}+H(a,x).
\end{equation}

In particular $K_D(x)=d_g(x,q_j)^{4\gamma_j}K_j(x)$ near $q_j$,
with $K_j>0$ of class $C^{2,\eta}$ for every $0<\eta<1$. We prove the
claimed regularity of the remainder $H(a,\cdot)$ in \eqref{eq:green-split}. Here the particular lower-order term of $P_g$ matters.
Write $r=d_g(a,x)$, $\theta=\exp_a^{-1}(x)/r$,
$B(\theta)=\operatorname{Ric}_a(\theta,\theta)$ and $R=R_g(a)$.
The normal volume density is $1-B(\theta)r^2/6+O(r^3)$.
Consequently $\Delta_g r=3/r-B(\theta)r/3+O(r^2)$ and
$\Delta_g\log r=2/r^2-B(\theta)/3+O(r)$. The remainders have the
corresponding differentiated bounds by the smooth polar expansion.
Since $\Delta_{\mathbb S^3}B=-8(B-R/4)$, for $r>0$ one obtains
\[
 \Delta_g^2\log r=r^{-2}(4B-2R/3)+O(r^{-1}).
\]
Indeed $\Delta_g(2r^{-2})=4B/(3r^2)+O(r^{-1})$, and the angular
term contributes $8(B-R/4)/(3r^2)$. For
$A=\tfrac23R_g g-2\operatorname{Ric}_g$, the leading term of
$\delta(A\,d\log r)$ is
$-r^{-2}(\operatorname{tr}A-2A(\theta,\theta))
=r^{-2}(2R/3-4B)$. Thus $P_g\log r=O(r^{-1})$ away from $a$.
Distributionally its delta coefficient is the flat one: integrating
over a shrinking sphere, the metric corrections to the principal flux
and the flux of the lower-order term are $O(r^2)$.
Subtracting the cutoff logarithm in \eqref{eq:green} therefore leaves
an $L^p$ source for every $1<p<4$. Elliptic regularity and the Sobolev
embedding give the stated regularity. In smoothly recentered coordinates
the cancellation holds at every pole; differentiating the smooth polar
coefficients with respect to that parameter preserves the $O(r^{-1})$
bound. After pullback to a fixed ball, derivatives of the operator
acting on the preceding $W^{4,p}$ remainder belong to $L^p$ as well.
Interior estimates (or the mean-zero inverse after extending the
coordinate change by the identity) give the same local regularity for
the recentered parameter derivatives. The same argument works for a
smooth metric parameter while $\ker P_g=\R$: in moving normal charts
the $r^{-2}$ coefficient cancels identically for every parameter, so the
differentiated error is still $O(r^{-1})$. Equivalently one can invert
$P_g+\Pi_{\R,g}$, where $\Pi_{\R,g}$ is projection onto the constants,
on fixed Sobolev spaces after pullback; its inverse varies smoothly
from $L^p$ to $W^{4,p}$. This supplies the derivative of the regular
factor used in Theorem~\ref{thm:D4}. In particular the regular factors
in \eqref{eq:Theta} admit a quadratic Taylor expansion at the pole.
We shall also differentiate the kernel at observation points distinct
from the source. Uniformly for $0<d_g(a,x)=r$ below a fixed small radius,
\[
 |\nabla_x^kG(a,x)|\le C_k r^{-k},\qquad 1\le k\le3 .
\]
Indeed on $B_{r/4}(x)$ the function $G(a,\cdot)-2\log(1/r)$ is bounded
by \eqref{eq:green-split}. On rescaling this ball to unit size, its
equation has uniformly smooth elliptic coefficients and a constant
right-hand side of order $r^4$. Interior estimates and rescaling give
the bounds. In particular these three derivatives are locally
integrable in the pole variable, uniformly in the observation point.

Put $\mathcal H=\{u\in H^2(M):\int_Mu\dd=0\}$ and
$s_*=\min_j(-1/\gamma_j)$, with $s_*=\infty$ for $N=0$.
Then $\inf_M K_D>0$ and $K_D\in L^s$ for $1<s<s_*$. Fontana's
inequality \cite[Theorem~1.7]{fontana} in dimension four, with two
derivatives, has sharp constant $32\pi^2$. Completing the square gives
\begin{equation}\label{eq:adams}
 \log\int_M e^{4a(v-\bar v)}\dd
 \le C+\frac{a^2}{8\pi^2}\|\Delta_gv\|_2^2,\qquad v\in H^2(M).
\end{equation}
A bounded sequence in $\mathcal H$ has an almost-everywhere convergent
subsequence; \eqref{eq:adams} at a larger exponent then gives strong
convergence of its exponentials in every finite $L^a$, $a\ge1$.
For $1<s<b<s_*$, H\"older with $K_D\in L^b$ yields strong convergence
of the weighted densities in $L^s$. Jensen gives
$\int K_De^{4u}\dd\ge V\inf K_D$ on $\mathcal H$. Consequently
$u\mapsto d\mu_u/dV_g$ maps bounded sets into relatively compact sets
in $L^s$. Its derivatives on $\mathcal H$ can be taken under the
integral: Taylor's remainder for $e^{4z}$ after order $j$ is bounded
by $C_j|z|^{j+1}e^{4|z|}$. H\"older, \eqref{eq:adams} and the
finite-exponent Sobolev embeddings bound the weighted $L^s$ remainder
by $C\|z\|_{H^2}^{j+1}$ on bounded sets with $\|z\|_{H^2}\le1$.
On $\mathcal H$ the variational functional and its first variation are
\begin{equation}\label{eq:energy}
 J_\kappa(u)=\langle P_gu,u\rangle-\kappa\log Z(u),\qquad
 \tfrac12D J_\kappa(u)[f]=\langle\mathcal R_\kappa(u),f\rangle ,
\end{equation}
where $\mathcal R_\kappa(u)=P_gu-2\kappa(\mu_u-dV_g/V)$.
Every weak solution of \eqref{eq:normalized} belongs
to $W^{4,s}$ for $1<s<s_*$ by elliptic inversion, and is smooth off
the marks by local bootstrap.

At a wall $\kappa_*$, a source atom of mass $16\pi^2m$ contributes
$8\pi^2m/\kappa_*$ to the limiting probability and $mG$ to the outer
field. Lemma~\ref{lem:virial} identifies $m$ with the weight of the
blow-up point.

For a type $(p,I)$ let $A=(a_1,\ldots,a_p)$ be the ordered ordinary
centers. Index components by the disjoint union of these $p$ labels and
$I$. Write $(z_s,\gamma_s,\mu_s)=(a_i,0,1)$ or $(q_j,\gamma_j,1+\gamma_j)$
accordingly, and $M_s=16\pi^2\mu_s$. The outer profile is
\begin{equation}\label{eq:outer}
 \Phi_{A,I}(x)=\sum_s\mu_sG(z_s,x).
\end{equation}
Near an ordinary center set $\widehat K_s=K_D$ and near an active
mark set $\widehat K_s=K_j$. With the centers fixed define
\begin{equation}\label{eq:Theta}
 \Theta_s^{A,I}(x)=\widehat K_s(x)
 \exp\!\left(4\mu_sH(z_s,x)+4\sum_{t\ne s}\mu_tG(z_t,x)\right).
\end{equation}
In a ball about a component the outer density is exactly
$r_s^{-4-4\mu_s}\Theta_s^{A,I}(x)$, where $r_s=d_g(x,z_s)$.
At an inactive mark $q_j$ it is
$d_g(x,q_j)^{-4+4\mu_j}K_j(x)e^{4\Phi_{A,I}(x)}$,
which is locally integrable since $\mu_j>0$.
The mixed interaction function is
\begin{equation}\label{eq:Hamiltonian}
 \cF_{p,I}(A)=\sum_s\mu_s\left[
 \log\widehat K_s(z_s)+2\mu_sH(z_s,z_s)
 +2\sum_{t\ne s}\mu_tG(z_s,z_t)\right].
\end{equation}
Only the regular points are variables, and direct differentiation gives
\begin{equation}\label{eq:gradF}
             \nabla_{a_i}\cF_{p,I}(A)
             =\left.\nabla_x\log\Theta_i^{A,I}(x)\right|_{x=a_i},\qquad 1\le i\le p.
\end{equation}

\section{The entire singular bubble}\label{sec:entire}

Fix $\gamma\in(-1,0)$ and put $\mu=1+\gamma$. The endpoint case of
Hyder--Mancini--Martinazzi \cite[Theorem~1.2]{HMM}, with $n=4$ and
$\alpha=\gamma$, gives a radial normal solution of
$\Delta^2 w=|x|^{4\gamma}e^{4w}$ with mass $16\pi^2\mu$.
Subtracting $\frac14\log6$ and then applying the dilation
\eqref{eq:scale} fixes its value at zero. Thus there exists a normal
entire solution $U_\gamma$, normalized by $U_\gamma(0)=0$, of
\begin{equation}\label{eq:entire}
 \Delta^2U_\gamma=6|x|^{4\gamma}e^{4U_\gamma}\quad\hbox{in }\R^4,
 \qquad \int_{\R^4}6|x|^{4\gamma}e^{4U_\gamma}\,dx=16\pi^2\mu.
\end{equation}
It is radial and unique up to
\begin{equation}\label{eq:scale}
 U_{\gamma,\lambda}(x)=U_\gamma(\lambda x)+\mu\log\lambda.
\end{equation}
Theorem~1.3 and Corollary~1.4 of \cite{AWZ} give this uniqueness and
\begin{equation}\label{eq:far-field}
 U_\gamma(x)=-2\mu\log|x|+c_\gamma+
                  O\bigl(|x|^{-\min(4\mu,1)}\bigr)
 \qquad(|x|\to\infty).
\end{equation}
For a regular core we use $U_0(y)=-\log(1+|y|^2/4)$, with
$U_0(0)=0$, $c_0=\log4$ and source mass $16\pi^2$.

Set
\begin{equation}
 W_\gamma(x)=|x|^{4\gamma}e^{4U_\gamma(x)},\qquad
 Z_\gamma(x)=x\cdot\nabla U_\gamma(x)+\mu,
\end{equation}
and
\begin{equation}
 L_\gamma\phi=\Delta^2\phi-24W_\gamma\phi.
\end{equation}
A component subscript $s$ on $U,W,Z,c$ or the model probability
denotes the profile of order $\gamma_s$. Differentiation in the
dilation parameter gives $L_\gamma Z_\gamma=0$.  Put $r=e^t$ and
\begin{equation}
                  V(t)=U_\gamma(e^t)+\mu t.
\end{equation}
Then
\begin{equation}\label{eq:Vode}
 V^{(4)}-4V''=6e^{4V},\qquad
 V'(-\infty)=\mu,\quad V'(+\infty)=-\mu.
\end{equation}

\begin{lemma}[Strict concavity]\label{lem:concavity}
With $p=-V''$ one has
\begin{equation}\label{eq:p-convolution}
 p(t)=\frac32\int_\R e^{-2|t-s|}e^{4V(s)}\,ds>0.
\end{equation}
Moreover, for $\sigma=\min(4\mu,2)$,
\begin{equation}\label{eq:p-endpoint}
 p(t)+|p'(t)|\le C(1+\mathbf1_{\{4\mu=2\}}|t|)e^{-\sigma|t|}.
\end{equation}
\end{lemma}

\begin{proof}
Equation \eqref{eq:Vode} gives
$(-D_t^2+4)p=6e^{4V}$.  Since $V(t)=-\mu|t|+O(1)$, $p$ is tempered and
the whole-line inverse gives \eqref{eq:p-convolution}. Convolving
$O(e^{-4\mu|t|})$ with $e^{-2|t|}$ gives \eqref{eq:p-endpoint}
(the factor $1+|t|$ occurs at equal exponents); $|p'|\le2p$.
Integrability of $p$ and the linear asymptotes of $V$ give the
two slopes in \eqref{eq:Vode}; the distance of $V'$ from its
limiting slope is bounded by the corresponding tail integral of $p$.
\end{proof}

\section{Full bounded Jacobi nondegeneracy}\label{sec:kernel}

\begin{definition}
A bounded Jacobi field is an $L^\infty(\R^4)$ distributional
solution of $L_\gamma\phi=0$.
\end{definition}
Since $W_\gamma\phi\in L^p_{\rm loc}$ for $1<p<-1/\gamma$, every
such field is continuous through the origin by elliptic regularity,
and smooth away from it.

\subsection{The radial mode}

\begin{theorem}[Radial orientation of the singular model]\label{thm:radpos}
Put $m=\int_{\R}e^{4V(t)}dt$, let $P(dt)=e^{4V(t)}dt/m$, and take a
radial variation $\phi(e^t)=f(t)$ with $u=f'\in H^1(\R)$. The form
\[
 B(u)=\int_\R \{(u')^2+4u^2\}\,dt-24m\,\mathrm{Var}_P(f)
\]
is nonnegative. Its nullspace consists of multiples of $V''$. On
their orthogonal complement in the kinetic inner product there is a
positive relative spectral gap.
\end{theorem}
\begin{proof}
Set $h=-V''>0$ and $A=-D_t^2+4$, so that $Ah=6e^{4V}$.
The probability $P$ has potential $-4V+\log m$ and potential
Hessian $4h>0$. For smooth compactly supported $u=f'$ the
Brascamp--Lieb variance inequality \cite[(1.3)]{CCL} gives
\[
 24m\,\mathrm{Var}_P(f)\le\int_\R (Ah/h)u^2\,dt .
\]
Writing $u=hz$ and integrating by parts yields
$B(u)\ge\int h^2|z'|^2$. This inequality extends on every compact
interval to $u\in H^1$ by density. Indeed, fixing $f(0)=0$,
$|f(t)|\le |t|^{1/2}\|u\|_2$; the exponential tails of $P$ and local
compactness make $u\mapsto f$ compact from $H^1$ to $L^2(P)$.
Consequently $B(u)=0$ implies $(u/h)'=0$. Conversely the primitive
$-V'$ of $h$ has zero $P$-mean, and integration of
$(V'e^{4V})'$ gives $4\int (V')^2 e^{4V}=\int h e^{4V}$.
Thus $B(h)=0$ by $Ah=6e^{4V}$ and \eqref{eq:p-endpoint}.
The compact self-adjoint spectral theorem gives the relative gap
on the kinetic orthogonal complement of $h$.
\end{proof}

\subsection{Angular modes}

For a spherical harmonic of degree $\ell$ put
$H_\ell=(-D_t^2+\ell^2)(-D_t^2+(\ell+2)^2)$ and $q=24e^{4V}$.
Its free quadratic form is
\[
 E_\ell(f)=\int_\R \big\{|f''|^2+(\ell^2+(\ell+2)^2)|f'|^2+
                             \ell^2(\ell+2)^2|f|^2\big\}\,dt .
\]
Indeed $r^2\Delta$ on this mode acts as
$D_t^2+2D_t-\ell(\ell+2)$. A second application, with the factor
$r^{-2}$, gives $r^4\Delta^2=H_\ell$. Thus the cylinder
linearization is $H_\ell-q$.

\begin{theorem}[Angular spectral gap]\label{thm:angular}
For $0<\mu<1$, $\ell\ge1$ and $f\in H^2(\R)$,
\begin{equation}\label{eq:Q1}
 E_\ell(f)-\int_\R q f^2\,dt
        \ge \frac{1-\mu}{1+\mu}\,E_\ell(f).
\end{equation}
Every bounded solution of $(H_\ell-q)f=0$ is zero.
\end{theorem}

\begin{proof}
The radial coefficient of $-\partial_iU_\gamma$ in the first
harmonic is $v=e^{-t}(\mu-V')>0$. Differentiating
\eqref{eq:entire} off the pole and using
$0<\mu-V'<2\mu$ gives
\[
 H_1 v=q(v-\gamma e^{-t})\ge
          \frac{1+\mu}{2\mu}\,qv .
\]
The inverse of $H_1$ has positive convolution kernel
$K(t)=\tfrac18(e^{-|t|}/2-e^{-3|t|}/6)$.
Both $v$ and $H_1v$ are $O(e^{a|t|})$ for some $a<1$ by
\eqref{eq:p-endpoint} and $q=O(e^{-4\mu|t|})$.
Consequently $v=K*(H_1v)$: the difference solves the homogeneous
equation, whose exponents are $\pm1,\pm3$, and has the same growth
bound. For the symmetric positive kernel of $S=q^{1/2}H_1^{-1}q^{1/2}$,
the weight $g=q^{1/2}v>0$ thus satisfies
$Sg\le[2\mu/(1+\mu)]g$. The weight is finite at every finite $t$. The weighted Schur test
\cite[Section~2]{GTschur} gives $\|S\|\le2\mu/(1+\mu)$.
Since $\|q^{1/2}H_1^{-1/2}\|^2=\|S\|$ and $E_\ell\ge E_1$,
this proves \eqref{eq:Q1} on every angular mode.
If $f$ is bounded and $H_\ell f=qf$, then $qf\in L^2$.
The Fourier multiplier $H_\ell^{-1}$ maps $L^2$ to $H^4$ and
the homogeneous equation has no nonzero tempered solution.
Hence $f\in H^4$; testing by $f$ in \eqref{eq:Q1} gives $f=0$.
\end{proof}

\begin{theorem}[Bounded singular Jacobi kernel]\label{thm:kernel}
Every bounded Jacobi field has the form
\begin{equation}
                              \phi=cZ_\gamma.
\end{equation}
\end{theorem}

\begin{proof}
Expand $\phi$ in spherical harmonics on compact annuli. All coefficients
are bounded on the cylinder. The radial coefficient
$f(t)$ satisfies $f^{(4)}-4f''=qf$, with $q=24e^{4V}$. Since
$f$ is bounded, $f''$ is tempered. Inverting $-D^2+4$ on tempered
distributions gives $f''=-\tfrac14 e^{-2|\cdot|}*(qf)$.
Thus $D^2f,D^3f$ decay exponentially at both ends. The limits of
$f'$ are zero by boundedness of $f$, so $f'\in H^1(\R)$.
Integration of the equation gives $\int q f=0$; testing with $f$
gives $B(f')=0$ for the form of Theorem~\ref{thm:radpos}.
Consequently $f=cV'+c_0$ and $c_0=0$, since $V'$ also has zero
$q$-mean. Theorem~\ref{thm:angular} eliminates the other modes;
$V'(\log r)=Z_\gamma(r)$ is itself a Jacobi field by dilation.
\end{proof}

\subsection{The first singular corrector}

The linear jet of the coefficient determines an odd correction to
the singular bubble. Its projection is included when differentiating
the bubble with respect to its parameters. For $r\ge0$ set
\begin{equation}
 \omega_\mu(r)=
 \begin{cases}
 r,&0\le r\le1,\\
 1+r^{(1-4\mu)_+},&r>1.
 \end{cases}
\end{equation}

\begin{proposition}[First-harmonic corrector]\label{prop:first-corrector}
For $b\in\R^4$ the normal logarithmic potential in the first angular
mode gives the unique solution in that class satisfying
\begin{equation}\label{eq:Psi-equation}
 L_\gamma\Psi_b=6W_\gamma(b\cdot x)\quad\hbox{in }\R^4,
 \qquad |\Psi_b(x)|+|x||\nabla\Psi_b(x)|
       \le C|b|\omega_\mu(|x|).
\end{equation}
By oddness, for every $R>0$,
\begin{equation}\label{eq:odd-cancellation}
 \int_{B_R}W_\gamma(x)\bigl(4\Psi_b(x)+b\cdot x\bigr)\,dx=0.
\end{equation}
The normal representation uses the logarithmic potential of the
full biharmonic source $6W_\gamma(b\cdot x+4\Psi_b)$, with the
biharmonic polynomial set to zero.
\end{proposition}

\begin{proof}
Put $D=d/dt$, $t=\log |x|$, and
$A_\mu(r)=\mu+Z_\gamma(r)=2\mu+rU_\gamma'(r)$. We claim that
\begin{equation}\label{eq:explicit-odd}
             \Psi_b(x)=-\frac{b\cdot x}{4\gamma}\,A_\mu(|x|).
\end{equation}
In the first spherical harmonic the free cylinder operator is
$H_1=(D^2-1)(D^2-9)$. Conjugating by $e^t$ gives
$e^{-t}H_1 e^t=(D+4)D(D^2-4)$. Applying this to $V'+\mu$ and
using $V^{(4)}-4V''=6e^{4V}$ yields
\[
 (H_1-24e^{4V})(e^t(V'+\mu))=-24\gamma e^t e^{4V}.
\]
This proves the equation off zero. There is no distribution supported
at the pole: after the linear term in \eqref{eq:explicit-odd} the
remainder and its first three derivatives have orders
$O(r^{1+\sigma-j}(1+|\log r|))$, $j=0,1,2,3$,
$\sigma=\min(4\mu,2)>0$, by \eqref{eq:p-convolution} and its equation.
All boundary terms tend to zero at the pole.

At infinity $A_\mu(e^t)=\int_t^\infty(-V''(s))\,ds$. Equation \eqref{eq:p-endpoint} gives
$|A_\mu(r)|+r|A_\mu'(r)|\le C r^{-\sigma}
(1+\mathbf1_{\{4\mu=2\}}\log r)$ for $r\ge2$. At zero $A_\mu$
is bounded and $rA_\mu'$ tends to zero. These estimates imply
\eqref{eq:Psi-equation} and sublinear growth at infinity.

In a first spherical harmonic the logarithmic kernel has
angular coefficient bounded by $C\min(r/s,s/r)$ at radii $r,s$.
For $s<r/2$ or $r<s/2$ subtract the angle-independent logarithm
and use the mean value bound. For comparable radii use the uniform
spherical $L^1$ bound on the logarithmic singularity.
The flat potential of $6W_\gamma(b\cdot x+4\Psi_b)$ consequently
converges in this harmonic and is sublinear: its large-radius bound
is $C(1+r^{(1-4\mu)_+})(1+\log r)$. The difference from
\eqref{eq:explicit-odd} is an entire sublinear biharmonic function,
hence constant by the interior derivative estimate on expanding balls;
oddness removes the constant. This proves the stated normal representation.
For uniqueness the difference of two such solutions is the potential
of $24W_\gamma$ times itself. If it is $O(1+r^a)$ at infinity,
$0<a<1$, the same split with this new source gives a bound by
$C(1+r^{\max(0,a-2\mu)})$: when $a=4\mu$ the kernel bound above
gives $O(1)$ directly. Iteration makes the difference bounded.
Theorem~\ref{thm:angular} then makes it zero. Radial symmetry of the
weight gives \eqref{eq:odd-cancellation}.
\end{proof}

\section{Sharp distributional tail expansion}\label{sec:tails}

Write the entire probability and its dilation as
\begin{equation}\label{eq:nu}
 p_\gamma(y)=\frac{6W_\gamma(y)}{16\pi^2\mu},\qquad
 d\nu_\gamma(y)=p_\gamma(y)\,dy,\qquad
 d\nu_{\gamma,\lambda}(x)=\lambda^4p_\gamma(\lambda x)\,dx .
\end{equation}
Equation \eqref{eq:far-field} gives
\begin{equation}\label{eq:nu-tail}
 \frac{d\nu_\gamma}{dy}(y)
  =A_\gamma|y|^{-4-4\mu}\bigl(1+o(1)\bigr),\qquad
 A_\gamma=\frac{3e^{4c_\gamma}}{8\pi^2\mu}.
\end{equation}

\begin{theorem}[Sharp moment trichotomy]\label{thm:moment}
Let $\varphi\in C_c^2(\R^4)$. As $\lambda\to\infty$:
\begin{enumerate}[label=\textup{(\roman*)}]
\item if $\mu>1/2$, then
\begin{equation}\label{eq:moment1}
 \int\varphi\,d\nu_{\gamma,\lambda}
 =\varphi(0)+\frac{M_{2,\gamma}}{8\lambda^2}\Delta\varphi(0)
       +o(\lambda^{-2}),
 \quad M_{2,\gamma}=\int|y|^2\,d\nu_\gamma<\infty;
\end{equation}
\item if $\mu=1/2$, then
\begin{equation}\label{eq:moment2}
 \int\varphi\,d\nu_{\gamma,\lambda}
 =\varphi(0)+\frac{|\Sph^3|A_\gamma}{8}
       \frac{\log\lambda}{\lambda^2}\Delta\varphi(0)
       +o(\lambda^{-2}\log\lambda);
\end{equation}
\item if $0<\mu<1/2$, then, for a radial $\chi\in C_c^\infty(\R^4)$ equal to one near $0$,
\begin{equation}\label{eq:moment3}
 \begin{split}
 \int\varphi\,d\nu_{\gamma,\lambda}
 &=\varphi(0)+A_\gamma\lambda^{-4\mu}\cD_\mu\varphi(0)
       +o(\lambda^{-4\mu}),\\
 \cD_\mu\varphi(0)
 &=\int_{\R^4}\frac{\varphi(x)-\varphi(0)
   -\chi(x)\nabla\varphi(0)\cdot x}{|x|^{4+4\mu}}\,dx.
 \end{split}
\end{equation}
The last integral is absolutely convergent and independent of the radial
cutoff.
\end{enumerate}
\end{theorem}

\begin{proof}
Let $m(r)$ be the spherical average of $\varphi(r\theta)$, put
$c=\Delta\varphi(0)/8$ and $v(r)=m(r)-m(0)$. Spherical symmetry in
dimension four gives $v(r)=cr^2+o(r^2)$ and
$|v(r)|\le C\min(r^2,1)$. The defect is $\int v(|y|/\lambda)d\nu_\gamma$.
For $4\mu>2$, \eqref{eq:nu-tail} gives a finite second moment, so
dominated convergence gives \eqref{eq:moment1}. At $4\mu=2$ the
second moment truncated at $R$ is $|\Sph^3|A_\gamma\log R+o(\log R)$.
On $|y|<d\lambda$ the quadratic Taylor error costs
$o_d(1)O(\lambda^{-2}\log\lambda)$, while the outer tail costs
$O_d(\lambda^{-2})$. Letting first $\lambda\to\infty$ and then
$d\downarrow0$ gives \eqref{eq:moment2}.

For $4\mu<2$, $\lambda^{4+4\mu}p_\gamma(\lambda x)
\le C|x|^{-4-4\mu}$ on all of $\R^4\setminus\{0\}$ by the
origin bound and \eqref{eq:nu-tail}. After subtracting
$\varphi(0)+\chi(x)\nabla\varphi(0)\cdot x$, the numerator is
$O(|x|^2)$ at zero and bounded at infinity. Dominated convergence
in the physical variable gives \eqref{eq:moment3}. The subtracted
linear term averages to zero for every radial cutoff.
\end{proof}

For $L^{-1}\Lambda\le\lambda_s^{\mu_s}\le L\Lambda$ with fixed $L$,
substitution into Theorem~\ref{thm:moment} gives the respective
orders $\Lambda^{-2/\mu_s}$, $\Lambda^{-4}\log\Lambda$ and
$\Lambda^{-4}$. At a regular core $\mu_s=1$ and the second moment
is finite. Hence the largest weight determines the common order
in \eqref{eq:Eorder}.

\section{The fourth-order global finite part}\label{sec:D}

For a pure singular type $(0,I)$ put $\Phi_I=\Phi_{\varnothing,I}$ in \eqref{eq:outer}.

\begin{definition}[Fourth-order negative-tail invariant]\label{def:D4}
Assume $\mu_j<1/2$ for all $j\in I$.  Define
\begin{equation}\label{eq:D4}
 D_I^{(4)}=\lim_{r\downarrow0}\left[
 \int_{M\setminus\bigcup_{j\in I}B_r(q_j)}
 K_De^{4\Phi_I}\,dV_g
 -\sum_{j\in I}\frac{|\Sph^3|}{4\mu_j}C_j(I)r^{-4\mu_j}
 \right],
\end{equation}
where $C_j(I)=\Theta_j^{\varnothing,I}(q_j)$.
\end{definition}

\begin{theorem}[Existence and canonicity of $D_I^{(4)}$]\label{thm:D4}
The limit in \eqref{eq:D4} exists and is finite. One may delete
$B_{r_j}(q_j)$ with $r_j\downarrow0$ independently, using $r_j^{-4\mu_j}$
in the corresponding counterterm. The value is independent of auxiliary
cutoffs. It is $C^1$ under smooth variation of the data with separated
marks, fixed active set $I$, weights in the stated ranges and
$\ker P_g=\R$.
\end{theorem}

\begin{proof}
In normal coordinates $x$ at an active mark put $r=|x|$. The outer density is
\begin{equation}\label{eq:local-D-density}
 K_De^{4\Phi_I}\,dV_g
 =C_j(I)r^{-4-4\mu_j}
 \left(1+\nabla\log\Theta_j(q_j)\cdot x+O(r^2)\right)dx.
\end{equation}
The Jacobian has zero linear jet and the regular factor is $C^{2,\eta}$
by Section~\ref{sec:green}. The constant term in \eqref{eq:local-D-density}
integrates on $r<|x|<r_0$ to
$|\Sph^3|C_j(I)(r^{-4\mu_j}-r_0^{-4\mu_j})/(4\mu_j)$.
The linear term has zero spherical mean. The remaining radial integral
over $(0,r_j)$ is $O(r_j^{2-4\mu_j})$, proving convergence with
independent radii and independence of the cutoff decomposition.

For varying data choose smoothly moving normal charts and fixed radial
cutoffs at all marks. After subtracting the same jets, the active
remainder and its parameter derivative are bounded by
$Cr^{-2-4\mu_j}(1+|\log r|)$ relative to $dx$.
This follows from the recentered Green estimates in Section~\ref{sec:green};
differentiation of a varying exponent produces the logarithm.
At an inactive mark the corresponding bound is
$Cr^{-4+4\mu_j}(1+|\log r|)$. Both admit integrable majorants
locally in the parameters. The constant counterterms are explicit and
the exterior is smooth, so dominated differentiation gives $C^1$ dependence.
\end{proof}

\section{Projected bubbles and mixed configurations}\label{sec:chart}

On a compact separated set of centers fix a common ball radius below
the injectivity radius so the normal balls are disjoint and each avoids
every other marked point. Transport a fixed radial cutoff, equal to
one on the half-ball, in these charts. Define the transported probability
measures by
\begin{equation}\label{eq:cutoff-measure}
 \int_M f\,d\nu_{j,\lambda}
 =\frac{\displaystyle
   \int_{\R^4}\chi_j(y/\lambda)f(\exp_{q_j}(y/\lambda))
          \,d\nu_{\gamma_j}(y)}
 {\displaystyle\int_{\R^4}\chi_j(y/\lambda)
          \,d\nu_{\gamma_j}(y)}.
\end{equation}
The regular measure $d\nu_{a,\lambda}$ is defined by the same formula with
$\gamma=0$, a cutoff centered at $a$, and the standard regular bubble.
The normalizing denominator is $1+o(1)$. We also write $\nu_s^M$ for a
transported probability when the scale and center are understood. Define
mean-zero projected bubbles by
\begin{align}
 P_g\varphi_{a,\lambda}
 &=16\pi^2\left(\nu_{a,\lambda}-\frac{dV_g}{V}\right),
 &\int_M\varphi_{a,\lambda}\dd=0,\\
 P_g\varphi_{j,\lambda}
 &=16\pi^2\mu_j\left(\nu_{j,\lambda}-\frac{dV_g}{V}\right),
 &\int_M\varphi_{j,\lambda}\dd=0.
\end{align}
Write $\varphi_s$ for either projection. Proposition~\ref{prop:green-sizes}
gives the rate of convergence to $\mu_sG(z_s,\cdot)$ in $C^3$ away from $z_s$.

Suppose a possibly signed density $g$ on $\R^4$ satisfies
\[
 |g(z)|\le C|z|^{-4+\delta}\quad (|z|\le1),\qquad
 |g(z)|\le C|z|^{-4-\delta}\quad (|z|>1)
\]
for some $0<\delta<4$. Then, uniformly for these bounds,
\begin{equation}\label{eq:core-log-potential}
\sup_y\int\left|\log\frac{1+|y|}{|y-z|}\right|\,|g(z)|\,dz\le C' .
\end{equation}
For bounded $y$, local H\"older with $1<p<4/(4-\delta)$ and the
integrable logarithmic tail give the bound. If $R=|y|>2$, the logarithm
is bounded on $|z|<R/2$. On $R/2\le |z|\le2R$, rescaling gives
$O(R^{-\delta})$ by local integrability of the translated logarithm.
The remaining tail is bounded by $C\int_{2R}^\infty r^{-1-\delta}\log(2r/R)\,dr$.

Let $\bar k$ be the number of negative eigenvalues of $P_g$,
counted with multiplicity, and let $v_1,\ldots,v_{\bar k}$ be an $L^2$-orthonormal basis of the negative
eigenspace,
\begin{equation}\label{eq:negative}
                     P_gv_r=-\vartheta_rv_r,\qquad \vartheta_r>0.
\end{equation}
On the mean-zero space use $P_g^+=|P_g|$ with inner product
\begin{equation}\label{eq:Pplus}
 \langle w,z\rangle_{P^+}
 =\langle P_gw,z\rangle_{L^2}
 +2\sum_{r=1}^{\bar k}\vartheta_r
       \langle w,v_r\rangle_{L^2}\langle z,v_r\rangle_{L^2}.
\end{equation}
Ellipticity and $\ker P_g=\R$ make its norm equivalent to the
$H^2$ norm on $\mathcal H$.

Each component has an amplitude and a scale. Only ordinary points vary in
position, in compact subsets of the ordered configuration space of
$M\setminus\Sigma$.

\section{The extended singular projection}\label{sec:ext}
Use the normal chart at each fixed singular point, a radial cutoff as in
\eqref{eq:cutoff-measure}, and the entire first harmonic of
Proposition~\ref{prop:first-corrector}. Put $k_j(A)=\nabla_x\log
\Theta_j^{A,I}(q_j)$ in an orthonormal frame. In particular $k_j$ varies
with the regular centers, even though the singular pole is fixed.
Ahmedou--Ben Ayed \cite[Section~2, (2.1)--(2.3)]{AhmedouBenAyed}
project the first angular correction in the planar Dirichlet problem
and include its center derivatives in the tangent space. Here the linear
jet of $\Theta_j$ forces the fourth-order profile of
Proposition~\ref{prop:first-corrector}. With $d\nu_\gamma$ as in
\eqref{eq:nu}, put $y=\lambda\exp_q^{-1}x$,
$c_\lambda=\int\chi(y/\lambda)d\nu_\gamma$, and
$B_b(y)=b\cdot y+4\Psi_b(y)$. Define the signed measure $S_{q,\lambda,b}$ by
\begin{equation}\label{eq:extended-source}
 \langle S_{q,\lambda,b},f\rangle =
 \frac{6}{\lambda c_\lambda}\int_{\R^4}
 \chi(y/\lambda)W_\gamma(y)B_b(y)f(\exp_q(y/\lambda))\,dy .
\end{equation}
Its integral vanishes by oddness on every concentric ball. The integration in
$dy$ specifies the volume Jacobian in the transported source. Let
$P_g\psi^{P}_{q,\lambda,b}=S_{q,\lambda,b}$, $\int_M\psi^P\,dV_g=0$, and set
\begin{equation}
 \widetilde\varphi_{j,\lambda,A}
 =\varphi_{j,\lambda}+\psi^P_{q_j,\lambda,k_j(A)},\qquad
 \widetilde U=\sum_{i=1}^p\alpha_i\varphi_{a_i,\lambda_i}
  +\sum_{j\in I}\alpha_j\widetilde\varphi_{j,\lambda_j,A}
  +\sum_r\beta_r v_r .
\end{equation}

Define the model error order
\begin{equation}\label{eq:e-mu}
 \mathfrak e_\mu(\lambda)=
 \begin{cases}
 \lambda^{-2},&\mu>1/2,\\
 \lambda^{-2}\log\lambda,&\mu=1/2,\\
 \lambda^{-4\mu},&0<\mu<1/2.
 \end{cases}
\end{equation}

\begin{proposition}[Size and derivatives of the extension]\label{prop:ext-size}
Fix a marked pole $q$ of order $\gamma\in(-1,0)$ and a compact set
of separated regular centers. With $\mu=1+\gamma$, write
\[
 d_\mu(\lambda)=
 \begin{cases}
 \lambda^{-1},&4\mu>1,\\
 \lambda^{-1}(1+\log\lambda)^{3/4},&4\mu=1,\\
 \lambda^{-4\mu},&4\mu<1 .
 \end{cases}
\]
Let $\ell=\log\lambda$ and let $b=b(A)$ be $C^2$ in local coordinates
of the regular centers. For $\lambda$ sufficiently large and
$i+|\nu|\le2$ one has
\[
 \|\partial_\ell^i\partial_A^\nu\psi^P_{q,\lambda,b(A)}\|_{H^2(M)}
 \le C|\partial_A^\nu b(A)|d_\mu(\lambda).
\]
On each compact $K\subset M\setminus\{q\}$ the same estimate holds in
$C^3(K)$ with $\mathfrak e_\mu(\lambda)$ from \eqref{eq:e-mu} in
place of $d_\mu(\lambda)$.
\end{proposition}

\begin{proof}
Set $D=y\cdot\nabla_y$. By \eqref{eq:explicit-odd},
$B_b(y)=(b\cdot y)(1-A_\mu(|y|)/\gamma)$ with
$A_\mu(e^t)=\mu+V'(t)$. Its radial factor has bounded Euler derivatives
through order two by \eqref{eq:p-endpoint}, which bounds $V''$ and $V^{(3)}$.
Since $D\log W_\gamma=4(Z_\gamma-1)$ and $DZ_\gamma=V''(\log|y|)$,
\[
 \big|D^j(W_\gamma B_b)(y)\big|
       \le C_\gamma |b|\,|y| W_\gamma(y),\qquad 0\le j\le2.
\]
All three functions are odd. In the fixed coordinate $v=\exp_q^{-1}x$,
\eqref{eq:extended-source} has density
$6\chi(v)c_\lambda^{-1}\lambda^3(W_\gamma B_b)(\lambda v)\,dv$.
Differentiating the scaled factor in $\ell$ acts by $3+D$. For $j=1,2$,
$\partial_\ell^j c_\lambda=O(\lambda^{-4\mu})$, since the differentiated
cutoff in its defining integral is supported on $|y|\asymp\lambda$.
The same holds for $c_\lambda^{-1}$ since $c_\lambda\to1$.
Each center derivative acts only on $b$. By linearity it suffices below
to estimate the pure scale derivatives.

On $A_k=\{2^k<|y|<2^{k+1}\}$ subtract the annular mean of
$f(\exp_q(y/\lambda))$, using oddness. Rescaled Poincar\'e and H\"older
bound the mean absolute oscillation by
$C\|\nabla f\|_{L^4(\exp_q(A_k/\lambda))}$. Take the chart slightly
larger than the support so that the last annulus is included.
For each scale derivative through order two the factor multiplying
this norm is bounded by
\[
 C|b|\lambda^{-1}
 \begin{cases}
  2^{(1+4\mu)k},&k<0,\\
  2^{(1-4\mu)k},&0\le k\le\log_2(C\lambda),
 \end{cases}
\]
by the endpoint bounds for $W_\gamma$. The $\ell^{4/3}$ norm of
these coefficients is $O(|b|d_\mu(\lambda))$. Disjointness of the
shells, $H^2(M)\hookrightarrow W^{1,4}(M)$ and mean-zero Paneitz
inversion give the $H^2$ estimate.

For $x\in K$, localize the source to a smaller fixed ball about $q$
disjoint from $K$ by a radial cutoff. Subtract the Green kernel's value
at the pole. Its oscillation, also after up to three $x$-derivatives,
is $O(2^k/\lambda)$. The resulting bound is
$C|b|\lambda^{-2}(1+\sum_{0\le k\le\log_2(C\lambda)}2^{(2-4\mu)k})
=O(|b|\mathfrak e_\mu(\lambda))$. The remaining fixed annulus has
density $O(|b|\lambda^{-4\mu})$ for each scale derivative. The Green
kernel and its first three observation derivatives are uniformly
integrable there, proving the $C^3(K)$ estimate.
\end{proof}

\begin{proposition}[The extended projection on the rescaled core]\label{prop:ext-core}
For a fixed singularity, put $D_\lambda(x)=\psi^P_{q,\lambda,b}(x)$.
Then $|D_\lambda(q)|\le C|b|\mathfrak e_\mu(\lambda)$ and, for each
fixed $R$, uniformly for $|y|\le R$,
\[
 \lambda\{D_\lambda(\exp_q(y/\lambda))-D_\lambda(q)\}
       =\Psi_b(y)+o_R(1)|b|.
\]
At fixed $x$ and $b$, with $l=\log\lambda$,
\begin{equation}\label{eq:ext-sup-der}
 \|D_\lambda\|_\infty+\|\partial_l D_\lambda\|_\infty
 \le C|b|(1+l)^2\lambda^{-\min(1,4\mu)} .
\end{equation}
\end{proposition}
\begin{proof}
At $q$ the logarithmic singularity of $G(q,\exp_q(z/\lambda))$ is
radial. Its pairing with \eqref{eq:extended-source} vanishes. The
regular part, after subtraction of its value at $q$, gains
$C|z|/\lambda$. The shell sum of Proposition~\ref{prop:ext-size}
gives the first bound.

Subtract the two Green potentials and multiply by $\lambda$. For
$|z|>L>2R$ the difference kernel is bounded by $C R/|z|$ by the
first Green derivative estimate and the triangle inequality. Its tail
against $|W_\gamma B_b|$ is $\le C_R|b| L^{-4\mu}$. On $|z|\le L$
the kernel tends, with the inverse normalization, to
$(8\pi^2)^{-1}\log(|z|/|y-z|)$. Uniform convergence of its integrals
in $|y|\le R$ follows by removing small balls about zero and $y$:
$W_\gamma B_b$ has a uniform local $L^{1+\epsilon}$ bound for some
$\epsilon>0$ and the logarithmic kernels have uniformly small local
$L^{(1+\epsilon)/\epsilon}$ norms on those balls. The metric is
asymptotically Euclidean on the complement. By the normal representation in Proposition~\ref{prop:first-corrector},
the limiting potential of $6W_\gamma B_b$, subtracted at zero,
is $\Psi_b(y)$ since $\Psi_b(0)=0$. Let $L\to\infty$.

For the last claim, on a physical shell of radius
$r_\ell=2^\ell/\lambda$ the average absolute value of the logarithmic
Green kernel in the source variable is at most
$C(1+|\log r_\ell|)$ uniformly in the observation point, by dilation
and local integrability of the logarithm. For the source and its log-scale derivative the same shell supremum
bounds hold by Proposition~\ref{prop:ext-size}. Summing the
negative shells gives $O(|b|\lambda^{-1}(1+\log\lambda))$.
Summing the positive shells gives at most
$C|b|(1+\log\lambda)^2\lambda^{-\min(1,4\mu)}$. This proves \eqref{eq:ext-sup-der}.
\end{proof}
\medskip
Let $E$ be the $P^+$-orthogonal complement of the parameter derivatives
of $\widetilde U$, and $\Pi_E$ the orthogonal projection onto $E$.
The regular center and singular scale derivatives are
\begin{align}
 \partial_{a_i}\widetilde U&=\alpha_i\partial_{a_i}\varphi_{a_i,\lambda_i}
       +\sum_{j\in I}\alpha_j\psi^P_{q_j,\lambda_j,\partial_{a_i}k_j},\\
 \partial_{\log\lambda_j}\widetilde U
 &=\alpha_j\big(\partial_{\log\lambda_j}\varphi_{j,\lambda_j}
              +\partial_{\log\lambda_j}\psi^P_{q_j,\lambda_j,k_j}\big).
\end{align}
The amplitude derivative at $q_j$ is $\widetilde\varphi_j$.
For example after normalization by $\lambda_i^{-1}$ the extra center
terms have norm $O(\lambda_i^{-1}\sum_j d_{\mu_j}(\lambda_j))$ on
compact separated parameter sets.

\subsection{The normal component of the gradient}
Fix disjoint balls $B_s$ inside the regions where the projection
cutoffs are one. Put $n=p+|I|$, $l_s=\log\lambda_s$,
$a=\max_s|\alpha_s-1|$, $B=\sum_j\beta_jv_j$ and
$k_s=\nabla\log\Theta_s(z_s)$. The two height factors are
\begin{equation}\label{eq:individual-height}
 L_s=\Theta_s(z_s)e^{-4c_s}\lambda_s^{4\mu_s},\qquad
 L'_s=L_s\exp\{8\mu_s(\alpha_s-1)l_s+4B(z_s)\}.
\end{equation}
In rescaled variables $y=\lambda_s\exp_{z_s}^{-1}x$ write
$K_De^{4\widetilde U}dV_g=L'_s g_s(y)dy$ on $B_s$.

For $l=\log\lambda$ put
\[
 E^{(2)}_\mu=(1+l)^2\mathfrak e_\mu(\lambda),\qquad
 E^{(1)}_\mu=\begin{cases}
 \lambda^{-1},&4\mu>1,\\
 \lambda^{-1}(1+l)^2,&4\mu=1,\\
 \lambda^{-4\mu}(1+l),&4\mu<1 .
 \end{cases}
\]
Each order indexed by $s$ is evaluated at $\lambda_s$.

\begin{proposition}[Normal defect of the projected approximation]
\label{prop:corrected-defect}
Let the centers range in a compact separated set, $\beta$ remain
bounded, $a\to0$ and $\min_s\lambda_s\to\infty$.
For every $w\in E$ with $\|w\|=o(1)$ the unnormalized mass of
each core satisfies
\begin{equation}\label{eq:unnormalized-core}
 \int_{B_s}K_D e^{4(\widetilde U+w)}\,dV_g
          =L'_s\bigl(M_s/6+o(1)\bigr).
\end{equation}
If in addition the $L'_s$ have bounded ratios, set
$T=(4n)^{-1}\sum_s\log(L'_s/3)$ and let $\kappa$ range in a compact
positive interval.
For $U=\widetilde U$ and the normal space of its full derivatives,
\begin{equation}\label{eq:source-defect}
 \left\|\Pi_E(P_g^+)^{-1}\mathcal R_\kappa(U)\right\|
 \le C\left(a+e^{-4T}+\sum_s E^{(2)}_{\mu_s}
        +|\beta|\sum_s E^{(1)}_{\mu_s}
        +\sum_{i=1}^p\lambda_i^{-1}|\nabla_{a_i}\cF_{p,I}|\right).
\end{equation}
All norms of functions in this estimate are $P^+$ norms.
\end{proposition}

\begin{proof}
We first compare the projected potentials with their entire profiles
throughout a fixed smaller core. In normal coordinates $x$ about $z_s$,
$r=|x|$, $y=\lambda_s x$, and with $b=k_s$ at a mark, one has
\begin{align}
 \varphi_s(\exp_{z_s}x)&=2\mu_s l_s+U_s(y)-c_s+
     \mu_sH(z_s,\exp_{z_s}x)+O((1+l_s)\mathfrak e_{\mu_s}),\notag\\
 \psi^P_{z_s,\lambda_s,b}(\exp_{z_s}x)&=
      \lambda_s^{-1}\Psi_b(y)+O(|b|(1+l_s)\mathfrak e_{\mu_s}).
 \label{eq:physical-projection}
\end{align}
To prove the first line compare the Green integral with its flat
logarithmic part. For source coordinate $v$ in the projection ball,
\[
 \left|\log\frac{d_g(\exp x,\exp v)}{|x-v|}\right|\le C|v|^2 .
\]
For $|x|\le2|v|$ this follows by normal metric comparison on a ball of
radius $C|v|$. For $|x|>2|v|$ the squared-distance error and its
first $v$-derivative vanish at $v=0$ by first variation. Its second
$v$-derivative on the segment to $v$ is $O(|x|^2)$ by smoothness of
squared distance and its Euclidean Taylor polynomial through order
three. Divide by $|x-v|^2$. The integral of $|v|^2$ against the
truncated model is $O(\mathfrak e_{\mu_s})$ by summing shells.
The smooth kernel has the same error after subtracting its constant
and its linear Taylor term in the source; the latter integrates to
zero. In physical coordinates the removed flat tail contributes
$O(\lambda_s^{-4\mu_s})$, since the observation core is strictly
inside the cutoff and $\int_1^\infty (1+\log r)r^{-1-4\mu_s}dr<\infty$.
The cutoff normalization differs from one by $O(\lambda_s^{-4\mu_s})$
and the flat logarithmic potential is $O(1+l_s)$. Its entire
representation gives the first line. At fixed $x$, its
log-scale derivative satisfies
\begin{equation}\label{eq:core-dilation}
 \partial_{l_s}\varphi_s(\exp_{z_s}x)
   =\mu_s+Z_s(\lambda_s x)+O((1+l_s)\mathfrak e_{\mu_s}).
\end{equation}
Indeed the uncut radial source derivative is $4Z_s$ times that
source and has zero mass. The cutoff is fixed in physical
coordinates, its normalization has derivative $O(\lambda_s^{-4\mu_s})$,
and the differentiated densities obey the same absolute shell bounds
and linear-moment cancellations. Differentiating the Green integrals
therefore preserves the error estimate.

For the second line the signed source on a physical shell of radius
$r'$ is dominated by $C|b|r'$ times the uncorrected density, by
\eqref{eq:Psi-equation}. Its mass on every sphere is zero.
The constant regular kernel drops out, its first-order error costs
$C|b|\mathfrak e_{\mu_s}$, and the metric error costs at most the
same on the fixed ball. In the omitted tail the first harmonic of
the flat logarithm is bounded by $Cr/r'$ for $r'$ outside the cutoff
transition. Integrating $C|b|r\lambda_s^{-4\mu_s}
(r')^{-1-4\mu_s}dr'$ gives $O(|b|\lambda_s^{-4\mu_s})$.
The full potential is $\lambda_s^{-1}\Psi_b(y)$, uniformly bounded
on the core by \eqref{eq:Psi-equation}. The normalization error is
again absorbed. This proves \eqref{eq:physical-projection}.

Divide the nonlinear density on core $s$ by $L'_s$ and express it in
$dy$. For unit amplitudes and $B=0$, \eqref{eq:physical-projection}
and the exterior estimates in Propositions~\ref{prop:ext-size}
and~\ref{prop:green-sizes} give
\begin{equation}
 \frac{g_s(y)}{W_s(y)}
 =1+\lambda_s^{-1}\bigl(k_s\cdot y+4\Psi_{k_s}(y)\bigr)
  +O\left(r^2+\sum_j(1+l_j)\mathfrak e_{\mu_j}\right).
\end{equation}
Set $\Psi=0$ at an ordinary point. Indeed the logarithm has the
displayed linear terms, the volume Jacobian is $1+O(r^2)$ and
$|\lambda_s^{-1}\Psi_{k_s}(y)|\le Cr$. At a mark both linear terms
occur in the projected source. At an ordinary point $k_s\cdot x$
remains. Varying the amplitudes, after the height extraction,
adds $O(a(1+|U_s(y)|))$ to the logarithm, by the first line of
\eqref{eq:physical-projection} and separation. The negative
eigenfunctions add $4(B(\exp x)-B(z_s))=O(|\beta|r)$.

Thus $g_s\to W_s$ uniformly on compact rescaled annuli, and on the
entire ball, by \eqref{eq:physical-projection} and \eqref{eq:far-field},
\begin{align}
 \varphi_s(\exp_{z_s}(y/\lambda_s))
   &=2\mu_s\log\frac{\lambda_s}{1+|y|}+O(1),\label{eq:raw-core-potential}\\
 g_s(y)&\le C|y|^{4\gamma_s}(1+|y|)^{-8\mu_s+8\mu_s a}.
 \label{eq:individual-shape}
\end{align}
To prove \eqref{eq:unnormalized-core}, use \eqref{eq:individual-shape}
to obtain summable geometric tails in
both shell directions, also in the scaled $L^{1+\epsilon}$ norms of
\eqref{eq:weighted-shell} for small $\epsilon>0$. Hence
$\int g_s\to\int W_s=M_s/6$. Under this conditional law,
$\int|e^{4w}-1|=O(\|w\|)$ by \eqref{eq:weighted-exp} and
$|e^z-1|\le |z|e^{|z|}$. The testing bound follows from the annular estimates and amplitude
orthogonality, independently of height ratios.

Test the difference on shell $2^j<|y|<2^{j+1}$
against $|f|$ for $f\in E$, $\|f\|=1$, or against $1$.
Lemma~\ref{lem:tight} bounds the shell mean by $C(1+|j|)^{3/4}$.
The amplitude error costs $O(a)$ after summation, since
$|e^z-1|\le |z|e^{|z|}$ and the model shell mass decays as
$2^{-4\mu_s|j|}$, absorbing $e^{Ca|j|}$. In making this split,
first remove the amplitude factor, leaving the unperturbed shell
mass for the other errors. The $r^2$ error on positive shells costs
at most
\[
 C\lambda_s^{-2}\sum_{0\le j\le l_s/\log2+C}
           2^{(2-4\mu_s)j}(1+j)^{3/4}\le C E^{(2)}_{\mu_s}.
\]
The analogous sum with $r$ in place of $r^2$ is bounded by
$CE^{(1)}_{\mu_s}$, hence by $C/\lambda_s$ at ordinary points. Negative shells are summable
with either power. Replacing the smaller core by the cutoff support
in a transported source costs $O(E^{(2)}_{\mu_s})$, including its signed
part. The same comparisons with test $1$ bound the change of the
normalizing integral: the corrected model has total mass one after
division by $M_s/6$. Its total variation paired with $|f|$ is bounded
by the same shell estimates. Thus conditional normalization preserves
the stated error orders. On the fixed exterior $U$ is bounded and
$K_D\in L^{1+\epsilon}$, so the normalized exterior pairing costs
$O(e^{-4T})$. Finally set $\widetilde\nu_s=\nu_s^M+S_s/M_s$, with $S_s$ from
\eqref{eq:extended-source} and $S_s=0$ at ordinary points.
Then $\int f\,d\widetilde\nu_s=0$ for each corrected source by amplitude and negative-mode orthogonality.
The kinetic term also annihilates $f$. Weighting the conditional
estimates by the actual core probabilities gives
\eqref{eq:source-defect}.
\end{proof}

\section{The constrained energy and localization}\label{sec:coerc}
\begin{lemma}[Energy of separated rescaled limits]\label{lem:core-energy}
Let $w_k\rightharpoonup w$ in $\mathcal H$. For finitely many uniformly
separated centers $z_{s,k}$ and scales $\lambda_{s,k}\to\infty$,
suppose $w_k(\exp_{z_{s,k}}(\cdot/\lambda_{s,k}))$, after subtraction
of constants, converges weakly in $H^2$ on compact annuli to $f_s$. Then
\[
 \|w\|_{P^+}^2+\sum_s\int_{\R^4}|\Delta f_s|^2\,dy
              \le\liminf_{k\to\infty}\|w_k\|_{P^+}^2 .
\]
\end{lemma}
\begin{proof}
For $v_k=w_k-w$, Hilbert norm splitting and Rellich compactness
in $H^1$ give
$\|\Delta_g v_k\|_2^2=\|w_k\|_{P^+}^2-\|w\|_{P^+}^2+o(1)$,
as the lower-order and finite-rank terms vanish. On every fixed
rescaled annulus $R^{-1}<|y|<R$ the metrics converge in $C^1$.
The rescaled Laplacian of $w_k$ converges weakly to $\Delta f_s$;
that of $w$ tends strongly to zero by absolute continuity of its
energy integral. In dimension four these annular energies are scale
invariant. Weak lower semicontinuity on the disjoint annuli and
then $R\to\infty$ prove the inequality.
\end{proof}

\begin{lemma}[Annular estimates for rescaled tests]\label{lem:tight}
In a normal chart of a separated core, let
$A_\ell=\{2^\ell<|y|<2^{\ell+1}\}$,
$f_\lambda(y)=f(\exp_z(y/\lambda))$, and let $m_\ell$ denote its
Euclidean average on $A_\ell$. For $1\le q<\infty$,
\begin{equation}\label{eq:entry-shell-test}
 \left(\operatorname{avg}_{A_\ell}|f_\lambda-m_0|^q\right)^{1/q}
 \le C_q(1+|\ell|)^{3/4}\|f\|_{H^2(M)} .
\end{equation}
For a unit normal test, $|m_0|\le C$. More generally, on any
$H^2$-bounded sequence with bounded $m_0$, the model integral of
$|f_\lambda|^2$ on shells $|\ell|\ge J$ inside the cutoff is bounded by
$C\sum_{|\ell|\ge J}(1+|\ell|)^{3/2}2^{-4\mu|\ell|}$. Consequently
local $L^2$ convergence to zero on compact annuli gives convergence
to zero in the model $L^2$ norm.
\end{lemma}
\begin{proof}
Scaled Poincar\'e on $\{2^\ell<|y|<2^{\ell+2}\}$ bounds $|m_{\ell+1}-m_\ell|$ by a constant times the
$L^4$ norm of the gradient there. Summing and applying discrete
H\"older gives $|m_\ell-m_0|\le C|\ell|^{3/4}\|f\|_{H^2}$ by bounded
overlap and $H^2\hookrightarrow W^{1,4}$. Scaled Sobolev--Poincar\'e
on the current annulus controls the normalized $L^q$ oscillation about
$m_\ell$. Together these bounds prove \eqref{eq:entry-shell-test}.

The model density is bounded on each shell by
$C|A_\ell|^{-1}2^{-4\mu|\ell|}$, with $\mu=1$ at a regular core.
For unit normal $f$, amplitude and negative-mode orthogonality give
$|\int f\,d\nu_s^M|=o(1)$ by Proposition~\ref{prop:ext-size}.
Integrating \eqref{eq:entry-shell-test} against the model gives
$|m_0|\le C$. With bounded $m_0$, apply the same estimate with $q=2$
and sum the weighted shell bounds. Local convergence deals with
the remaining finitely many shells.
\end{proof}

\begin{lemma}[Uniform weighted testing]\label{lem:weighted-testing}
Fix a type and a compact separated set of centers. Let all scales tend to
infinity and let $E$ be the corrected normal space. On core $s$, let
$\eta_{s,\lambda}$ be a conditional probability whose rescaled density
$\rho_{s,\lambda}$ satisfies, for some $p>1$ and $c>0$,
\begin{equation}\label{eq:weighted-shell}
 |A_\ell|^{1-1/p}\|\rho_{s,\lambda}\|_{L^p(A_\ell)}
       \le C2^{-c|\ell|},\qquad
 A_\ell=\{2^\ell<|y|<2^{\ell+1}\},
\end{equation}
and converges in local $L^p$ on compact annuli to the entire probability
$d\nu_s/dy$. The bounds may instead have a fixed subexponential factor
$\exp(C(1+|\ell|)^{3/4})$, after decreasing $c$. Then
\begin{equation}\label{eq:weighted-operator}
 \sup_{\substack{w,z\in E\\\|w\|,\|z\|\le1}}
 \left|\int wz\,(d\eta_{s,\lambda}-d\nu_s^M)\right|\longrightarrow0.
\end{equation}
The analogous assertion holds for first moments. For every fixed $a>0$,
\begin{equation}\label{eq:weighted-exp}
 \sup_{w\in E,\ \|w\|\le1}\int e^{a|w|}\,d\eta_{s,\lambda}
       \le C_a.
\end{equation}
All constants are uniform for these bounds, fixed strengths and separated
centers. The norm is the $P^+$ norm.
\end{lemma}
\begin{proof}
Lemma~\ref{lem:tight} bounds the normalized $L^q$ norm of a unit
normal test on $A_\ell$ by $C_q(1+|\ell|)^{3/4}$ for every finite $q$.
Let $p'=p/(p-1)$. On each shell, H\"older and
\eqref{eq:weighted-shell} give
\[
 \begin{split}
 \int_{A_\ell}|wz|\rho_{s,\lambda}\,dy
 &\le C2^{-c|\ell|}
 \left(\operatorname{avg}_{A_\ell}|w|^{2p'}\right)^{1/(2p')}
 \left(\operatorname{avg}_{A_\ell}|z|^{2p'}\right)^{1/(2p')}\\
 &\le C(1+|\ell|)^{3/2}2^{-c|\ell|}\|w\|\|z\|.
 \end{split}
\]
The same bound holds for the model. A fixed subexponential factor is
absorbed by decreasing $c$. On finitely many shells, local $L^p$
convergence and the same $L^{2p'}$ bounds give convergence uniformly
in $w,z$. This proves
\eqref{eq:weighted-operator}, first sending the scales to infinity and
then $L\to\infty$. The first-moment assertion is identical.

On each annulus the rescaled $H^2$ norm of $w-m_\ell$ is bounded:
Poincar\'e controls its lower derivatives and the Hessian energy is
scale invariant in dimension four. Extension from the fixed annulus
and Adams' inequality consequently give
$\operatorname{avg}_{A_\ell}e^{b|w-m_\ell|}\le C_b$ for every fixed
$b$. H\"older and \eqref{eq:weighted-shell} now give the summable
bound $C_a2^{-c'|\ell|}\exp(C_a(1+|\ell|)^{3/4})$ for each shell.
Summing proves \eqref{eq:weighted-exp}. Changes of charts and the last
incomplete shell are handled in a slightly larger fixed chart.
\end{proof}

\begin{proposition}[Balanced coercivity on the current normal space]\label{prop:balancedgap}
Fix a type, fixed strengths and a positive lower bound on the separation
of the centers. In the extended charts suppose $\min_s\lambda_s\to\infty$,
$\alpha_s\to1$ and $|\beta|$ is bounded. Set
$\theta_s=8\pi^2\mu_s/\kappa_*$ and
$\eta^0=\sum_s\theta_s\nu_{s,\lambda_s}^{M}$.
On the corrected normal space $E$ of these same parameters one has
\begin{equation}\label{eq:reference-gap}
 \|w\|_{P^+}^2-8\kappa_*\operatorname{Var}_{\eta^0}(w)
       \ge c_0\|w\|_{P^+}^2,\qquad w\in E ,
\end{equation}
eventually along every such sequence, with $c_0>0$ depending only on
the fixed type, metric and separation.
\end{proposition}
\begin{proof}
Take unit vectors $w_k$ orthogonal to the corrected tangents. Pass to a
subsequence with weak local rescaled limits $f_s$ on annuli. The bounds on the annular means and local elliptic Sobolev
compactness give strong $L^2$ convergence on each compact annulus.
The same shell argument bounds local $L^2$ norms through the pole;
in particular the limiting profiles have finite second-derivative
energy and gradients in $L^4$. They extend across the pole in the
energy class: integration by parts with cutoffs on a dyadic shell incurs,
in the distributional derivative identification, factors bounded by
$Cr^2(1+|\log r|)$ for the terms with two cutoff derivatives and by
$Cr^2\|\nabla f_s\|_{L^4}$ for the gradient term. Both vanish.
The amplitude equations give zero model means. Log-scale orthogonality
gives $\int f_s Z_s\,d\nu_s=0$, since the log-scale derivative of the
uncut density is $4Z_s$ times the density. At a regular center its four
normalized center derivatives give orthogonality to the four translation
Jacobi fields by the same change of variables. In fixed coordinates their
leading density derivatives are the four translates of the rescaled
regular density; smooth chart, cutoff and Jacobian derivatives are $o(1)$
after the factor $\lambda_i^{-1}$, by first truncating to finitely many
shells. The density logarithmic derivatives of the regular radial factor
are bounded on these shells and in the tail; changing the moving frame
only adds a rotation of a radial density. The extensions contribute $o(1)$ by Proposition~\ref{prop:ext-size}.

After extraction $w_k$ also has a weak limit in $\mathcal H$.
Lemma~\ref{lem:core-energy} gives $\sum_s\int|\Delta f_s|^2\le1$.
Lemma~\ref{lem:tight} removes the model tails uniformly, so along
the extracted subsequence
\[
 8\kappa_*\lim_k\int w_k^2\,d\eta^0
   =\sum_s24\int_{\R^4}W_s f_s^2\,dy .
\]
For the radial average $f(t)$ of a limiting profile, finite flat
energy gives $g=f''+2f'\in L^2(\R)$, and the $L^4$ gradient bound
gives $f'\in L^4(\R)$. The solution of $(D+2)u=g$ obtained by
convolution with $e^{-2t}\mathbf1_{t>0}$ belongs to $H^1$.
Its difference from $f'$ is $ce^{-2t}$, excluded by the $L^4$
condition. Thus $f'\in H^1$ and
$\int|f''+2f'|^2=\int(|f''|^2+4|f'|^2)$ by integration by parts.
This is the kinetic norm in Theorem~\ref{thm:radpos}. With $h=-V''$, its kinetic pairing with
$f'$ equals $6\int f'e^{4V}=-24\int fV'e^{4V}=0$ by the dilation
orthogonality and the exponential tails. Each nonconstant angular
coefficient $\psi_\ell$ belongs to $L^4(dt)$ by spherical
Poincar\'e and the $L^4$ gradient bound. Its flat Laplacian gives
$(D-\ell)(D+\ell+2)\psi_\ell\in L^2(dt)$. This constant-coefficient
operator has a bounded inverse $L^2\to H^2$ by its Fourier symbol.
The two free exponentials are excluded by $L^4$. Thus its angular
form is defined on $H^2(dt)$. Summing the orthogonal modes and applying
the radial and angular gaps proves the constrained singular estimate.

For a regular core the sphere has Paneitz eigenvalues
$\ell(\ell+3)(\ell(\ell+3)+2)$. Its constant and five first harmonics
are exactly the removed mean, dilation and translations. To justify
stereographic transfer for the limiting profile $f_s$, take a cutoff
equal to one on $B_R$, zero outside $B_{2R}$, and set
$g_R=c_R+\chi_R(f_s-c_R)$, where $c_R$ is the mean on the intervening
annulus. Scaled Poincar\'e, $\nabla f_s\in L^4$ and
$\nabla^2 f_s\in L^2$ give $\|\Delta(g_R-f_s)\|_2\to0$ by the
product rule and vanishing tail energies. Also $g_R\to f_s$ in
$L^2(d\nu_s)$ by the shell bounds, since $|c_R|\le
C(1+\log R)^{3/4}$ for dyadic $R$. Each $g_R$ is constant near the
point at infinity. Local smooth approximation and conformal
invariance identify its spherical energy with
$\int|\Delta g_R|^2$. Projection away from the constant and first
harmonics costs $o(1)$ in the limit by the model orthogonalities.
On the unit sphere these are the modes with eigenvalues $0$ and $24$;
the next eigenvalue is $120$, giving the regular relative gap $1-24/120$.

Let $\delta_{s,\mathrm{rad}}>0$ be the radial relative gap of
Theorem~\ref{thm:radpos}. For each singular component set
$\delta_s=\min\{\delta_{s,\mathrm{rad}},(1-\mu_s)/(1+\mu_s)\}$,
and set $\delta_s=4/5$ for a regular component. The finite minimum
$\delta=\min_s\delta_s$ is positive for the fixed strengths.
Therefore
\[
 8\kappa_*\limsup_k\operatorname{Var}_{\eta^0}(w_k)
 \le\sum_s24\int W_sf_s^2
 \le(1-\delta)\sum_s\int|\Delta f_s|^2\le1-\delta ,
\]
after extraction along the limsup. This proves \eqref{eq:reference-gap} with $c_0=\delta/2$.
\end{proof}
\subsection{Modulation in the extended family}
At an ordinary point write $\widetilde\varphi_s=\varphi_s$.

\begin{lemma}[Uniform local modulation]\label{lem:modulate}
Suppose all scales tend to infinity, $1/2<\alpha_s<2$ and $|\beta|$
is bounded. In a fixed local ordering of the regular points, functions
within a sufficiently small fixed $P^+$ distance of the reference sum
admit a unique orthogonal representation in a fixed small parameter
ball about that sum. In this assertion increments are measured by
\[
 \sqrt{\log\lambda_s}\,\delta\alpha_s,\quad
 \delta\log\lambda_s,\quad
 \lambda_i\,\delta a_i\ (i\ {\rm regular}),\quad \delta\beta .
\]
The norm of the parameter increment and of the normal error is at most
a constant times the initial distance. The constants are uniform for
a fixed type, a lower separation bound, and reference amplitudes in
the indicated interval.
\end{lemma}
\begin{proof}
Freeze the normalizations at the reference point. Pair each normalized
tangent with its source. The projected bubble on its own core is
$2\mu_s l_s+O(1+\log(1+|y|))$, so the normalized amplitude Gram entry
tends to $2\mu_sM_s=32\pi^2\mu_s^2$. Distinct amplitude pairings
before normalization are bounded by separation.
By \eqref{eq:core-dilation} and \eqref{eq:ext-sup-der}, the log-scale
tangent has rescaled limit $\mu_s+Z_s$ with source $4M_sZ_s\,d\nu_s$.
At an ordinary pole the normalized center tangent has limit
$-\partial_jU_0$ with source $-4M_s\partial_jU_0\,d\nu_0$, by
differentiating the translated density. Chart and cutoff errors on
fixed annuli are $O(\lambda_i^{-1})$, and the tail bound below is summable.
Thus the scale and center block, up to the amplitude factors, tends to
$4M_s$ times the model Gram matrix of dilation and translations.
These mean-zero modes are independent. Cross-core pairings vanish in
the limit by separation and the exterior estimates. The bounded core
tangents pair to $o(1)$ with the normalized amplitude, as do the smooth
negative modes, whose Gram matrix is $\mathrm{diag}(\vartheta_r)$.
Their pairings with the scale and center tangents tend to zero:
the differentiated sources have zero mass and concentrate at their poles.
Proposition~\ref{prop:ext-size} controls the extension derivatives.
The normalized Gram matrices are therefore uniformly positive and bounded.

For the uniform Taylor estimate, each first or second log-scale or
normalized regular-center derivative of a transported source has total
mass zero. On a physical dyadic shell $|x-z_s|\asymp2^\ell/\lambda_s$
its density times the shell volume is bounded by
$C2^{-4\mu_s|\ell|}$, with a bounded constant on the fixed transition
annulus after division by that bound. Indeed in the flat factor the
dilation derivatives are generated by $4+y\cdot\nabla$; $Z_s$ and
$y\cdot\nabla Z_s$ are bounded. Center derivatives occur only for
the explicit regular density, whose first two scaled derivatives
satisfy the same bounds. Derivatives of the inverse exponential chart
in normalized center coordinates are bounded after the source rescaling,
and orthonormal frame rotations annihilate the radial density. Smooth
Jacobian and cutoff factors and their normalized derivatives are
bounded. For a unit $H^2$ test function subtract its first-shell mean
from the zero-mass pairing. The mean absolute deviation on shell $\ell$
is $O(1+|\ell|^{3/4})$ by the two-shell Poincar\'e estimate in
Lemma~\ref{lem:tight}. Summation proves a uniform dual bound
for these derivatives, with a summable bound for the tails used above.
Proposition~\ref{prop:ext-size} gives the corresponding $H^2$ estimates
for the extension, since $k_s(A)$ has bounded derivatives through order
two on the separated center set. Mixed amplitude derivatives are bounded
after the amplitude normalization.
Elliptic inversion now bounds the second differential of the extended
sum in these frozen coordinates on a small fixed ball.

In the normalized increment $x$, Taylor's theorem now gives
$\|\widetilde U(x)-\widetilde U(0)\|\ge c|x|$ on a fixed small ball.
Half the Hessian of $\|u-\widetilde U(x)\|^2$ is the tangent Gram
matrix plus $O(\|u-\widetilde U(x)\|)$. Thus the squared distance is
strictly convex for $u$ sufficiently close to $\widetilde U(0)$,
and its minimum is interior by the boundary lower bound. The critical
equations are the required orthogonalities. Comparison with $x=0$
bounds both the increment and the normal error.
\end{proof}

\section{Entry of exact solutions and coefficients of the dilation test}\label{sec:rate-low}
\subsection{Exact blow-up, matching, and entry}
For $u\in\mathcal H$ write $F=2\kappa\mu_u$. For an exact sequence
in the mass gauge \eqref{eq:mass-gauge} this is
$F_k=2\kappa_k K_De^{4v_k}dV_g$.

\begin{lemma}[Localized virial identity]\label{lem:virial}
For a smooth vector field $X$ vanishing at the marks, set
\[
 A_X=\operatorname{div}_gX+X\log K_D,\qquad
 B_X(x,y)=(X_x+X_y)G(x,y).
\]
Let $Z_X$ be the mean-zero solution of $P_gZ_X=-4\operatorname{div}(XF)$. Then
\begin{equation}\label{eq:offshellvir}
 \langle\mathcal R_\kappa(u),Z_X\rangle=-V_X(F),\qquad
 V_X(F)=\int A_X\,dF+\frac1{8\pi^2}\iint B_X\,dF\,dF .
\end{equation}
In particular every exact solution satisfies
\begin{equation}\label{eq:symvirial}
                            V_X(F_k)=0.
\end{equation}
For a concentrating sequence of exact solutions, allow the smooth
positive coefficients to converge in $C^\infty$ to $h>0$.
At a blow-up point $z_s$ the limiting atom has mass $M_s=16\pi^2\mu_s$.
The limiting regular centers satisfy $\nabla_A\cF_{p,I}=0$ for the
interaction of the limiting coefficient.
\end{lemma}
\begin{proof}
Since $X$ vanishes at the marks, $A_X$ is bounded, and $B_X$ is bounded
by first variation. Near the diagonal the simultaneous
derivative of $\log d_g(x,y)$ is the arclength average of
$\tfrac12(\mathcal L_X g)(\dot c,\dot c)$ on the unit-speed minimizing
geodesic; the first derivatives of the regular part are bounded.

First take smooth $u$. The flux of $K_De^{4u}X$ through the sphere
of radius $r$ about $q_j$ is $O_u(r^{4\mu_j})$. Integrating its
divergence and applying \eqref{eq:green} give
\[
 \langle P_gu,Z_X\rangle=4\int Xu\,dF=-\int A_X\,dF,\qquad
 Z_X(y)=\frac1{4\pi^2}\int X_xG(x,y)\,dF(x).
\]
For two comparable dyadic shells of radius $r$ at a mark, the
rescaled densities are $O_u(r^{4\mu_j})$ and the differentiated
kernel is bounded by $C/|x'-y'|$ in rescaled variables. Summing
$O_u(r^{8\mu_j})$ and using $|X(x)|/d_g(x,y)=O(1)$ on noncomparable
shells gives absolute integrability. Symmetrization and the zero
mean of $Z_X$ give \eqref{eq:offshellvir}.

For smooth $u_n\to u$ in $\mathcal H$, the densities of $F_n=2\kappa\mu_{u_n}$
converge in $L^s$, $1<s<s_*$, by Section~\ref{sec:green}.
Near $q_j$, $|X|K_D=O(r^{4\mu_j-3})$, hence $XK_D\in L^{4/3+\epsilon}$
for some $\epsilon>0$. H\"older and exponential convergence give
$XF_n\to XF$ in $L^{4/3}$. Since $H^2\hookrightarrow W^{1,4}$,
the corresponding test functions converge in $H^2$ by mean-zero inversion.
The residuals converge in $H^{-2}$ and the bounded kernels pass to
the limit, proving the identity on $\mathcal H$.

For the normalization in \cite[(1.9)--(1.11)]{AWZ} put
\[
 U_k=v_k-8\pi^2\sum_j\gamma_jG_0(\cdot,q_j),\qquad
 c_k=\int_M e^{4U_k}\dd .
\]
Since $\int_M h_ke^{4f}e^{4U_k}\dd=1$, $c_k$ is bounded above and
away from zero. The function $U_k-\frac14\log c_k$ has unit
exponential mass; the coefficient and background of its equation
are $c_k\kappa_k h_ke^{4f}$ and $\kappa_k/V$. After subselection they
converge in $C^1$ with positive limiting coefficient. Theorem~5.1
of \cite{AWZ} now gives $F_k\rightharpoonup\sum_s\sigma_s\delta_{z_s}$
in the blow-up alternative.

The coefficient term $X\log K_{D,k}$ differs from $X\log K_D$ by
$o(1)$ uniformly, by convergence of the positive smooth factors.
Take a cutoff dilation $X=x$ in normal coordinates about one atom,
supported away from every other atom and every other mark. At its center $A_X$
tends to $4\mu_s$; on its self-pair $B_X$ tends uniformly to $-2$
by first variation, also when the two arguments coalesce. All other
pair limits are zero. Concentration of $F_k$ and $F_k\otimes F_k$
therefore gives
\[
 0=4\mu_s\sigma_s-\frac{\sigma_s^2}{4\pi^2},\qquad
 \sigma_s>0.
\]
In particular an ordinary atom is simple. For an ordinary center
take instead $X=e$ near $z_s$, constant in its normal coordinates.
Now $\nabla X(z_s)=0$, so the singular part of the self-pair
has limit zero. Write $\partial_2$ for differentiation in the second argument.
Dividing by $\sigma_s$ gives
\[
 0=\partial_e\log\widehat K_s(z_s)
   +\frac{\sigma_s}{4\pi^2}\partial_{2,e}H(z_s,z_s)
   +\sum_{t\ne s}\frac{\sigma_t}{4\pi^2}\partial_{2,e}G(z_t,z_s)
   =\partial_e\log\Theta_s(z_s).
\]
Equation \eqref{eq:gradF} proves the center condition.
\end{proof}

\begin{theorem}[Strong entry and matching for exact solutions]\label{thm:strongentry}
Let $\kappa_k\to\kappa_*>0$ and let $u_k$ be exact solutions of
\eqref{eq:normalized} with $\|u_k\|_{H^2}\to\infty$. Allow the
smooth positive coefficient $h_k$ in this equation to converge in
$C^\infty(M)$ to a positive $h$. In sequence terms the weight,
matching factor and projected family use $h_k$. In limits they use
$h$. Use \eqref{eq:mass-gauge} for $v_k$. After subselection, $u_k$ has an
orthogonal representation in the extended family with natural amplitudes
$\zeta_{s,k}=(\alpha_{s,k}-1)\sqrt{\log\lambda_{s,k}}$ satisfying
\[
 \sum_s\zeta_{s,k}^2+|\beta_k|^2+\|w_k\|_{P^+}^2\longrightarrow0.
\]
Moreover $(\alpha_{s,k}-1)\log\lambda_{s,k}\to0$ for every component.
For both the peak scales and the resulting chart scales,
\begin{equation}\label{eq:matching-identity}
  \kappa_k\tau_k\,\Theta_s^{A_k,I}(z_{s,k})
         e^{-4c_s}\lambda_{s,k}^{4\mu_s}\longrightarrow 3.
\end{equation}
For the peak scales, on each fixed small core
$B_s=B(z_{s,k},r_0)$, with $r_s=d_g(x,z_{s,k})$,
\begin{equation}\label{eq:exacttail}
 \frac{dF_k}{dV_g}\le C\tau_k r_s^{-4-4\mu_s},\qquad
 F_k(B_s)=M_s+O(e^{-c|\bar v_k|})
\end{equation}
for some $c>0$. On fixed complements of the limiting poles
the densities of $F_k/\tau_k$ converge in $L^1$ to
$2\kappa_*K_De^{4\Phi_{A,I}}$. For
$L_k=\tfrac14\log(\tau_k^{-1}/\kappa_*)$, one has $L_k\to\infty$
and
\begin{equation}\label{eq:peak-energy}
 |\kappa_k-\kappa_*|=O(e^{-c L_k}),\qquad
 \langle P_gu_k,u_k\rangle=4\kappa_k L_k+O(1).
\end{equation}
\end{theorem}
\begin{proof}
Use the mass gauge and normalization in Lemma~\ref{lem:virial}.
The bounded alternative of
\cite[Theorem~5.1]{AWZ} contradicts
$\|u_k\|_{H^2}\to\infty$ by elliptic inversion, and uniform decay
contradicts unit weighted mass. Lemma~\ref{lem:virial} therefore
identifies the limiting atoms. The regular factors in
\eqref{eq:Theta} are uniformly bounded in $C^2$ on the separated
core charts and converge with their first two derivatives.
Fix $K\Subset U$ off these atoms and
let $m_k=F_k(U)\to0$, choosing $\overline U$ disjoint from the atoms.
In the Green representation of $u_k$, the exterior source has bounded
potential on $K$. Jensen's inequality for $F_k|_U/m_k$ bounds
$\int_K e^{4a u_k}\dd$ for each fixed $a>1$, since the kernel
$d_g(x,y)^{-am_k/(2\pi^2)}$ is uniformly integrable in $x$.
As $K_D\in L^b$ for some $b>1$, this bounds $K_De^{4u_k}$ locally in
$L^q$, $q>1$, also at inactive marks. Jensen bounds $\tau_k$ above,
so interpolation with the vanishing local mass gives $F_k\to0$
locally in $L^{q_1}$ for $1<q_1<q$. The nearby logarithmic potential
then tends uniformly to zero by H\"older, and the separated potential
converges uniformly by weak convergence and a smooth cutoff. Thus
$u_k\to\Phi_{A,I}$ uniformly on $K$. On a fixed ball off the atoms
the limiting weighted integral is positive and its $F_k$-mass
vanishes, whence $\tau_k\to0$ and
$F_k/\tau_k\to2\kappa_*K_De^{4\Phi_{A,I}}dV_g$ in $L^1(K)$.
Choose disjoint fixed balls about the limiting atoms. On each let $m$
be the maximum of $v_k$, attained in its interior for large $k$.
Put $z=q_j$ at a singular atom and take $z$ to be a maximum point at
a regular atom. In this proof $z,\lambda,\gamma,\mu$ refer to that
component, and
\[
 \lambda^{4\mu}=\kappa_k\widehat K(z)e^{4m}/3,\qquad
 V_k(y)=v_k(\exp_z(y/\lambda))-m .
\]
At a singular pole, the spherical Harnack estimate
\cite[Theorem~5.5]{AWZ} bounds $m+\mu\log d_g(x_{\max},z)$,
where $x_{\max}$ is the maximizer. The regular parts in the
two gauges differ by a uniformly bounded constant. Hence
$\lambda d_g(x_{\max},z)$ is bounded.
At an ordinary pole the maximum point is $y=0$.
On every fixed rescaled ball the source density is bounded by
$C_R |y|^{4\gamma}$ since $V_k\le0$. Subtract the Green formula at
the maximum. Its contribution from rescaled distances exceeding $R$
is bounded by $C_L/R$ for observation points in $B_L$, $R>2L$
enlarged to contain the maximum. This follows from
$|\nabla_xG|\le C/d_g(x,y)$ and the total mass bound. On the inner
ball, convolution with the logarithm is locally equicontinuous for
densities bounded in $L^{1+\epsilon}$, by translation continuity in
the dual finite $L^q$ norm. The regular parts converge in differences.
After extraction $V_k$ converges locally uniformly to a profile
$V$ with maximum zero. The preceding tail estimate permits passage
to the limit in the subtracted logarithmic representation, so
$\Delta^2 V=6|y|^{4\gamma}e^{4V}$ and $V$ is normal. The same
representation gives $V(y)\ge-C(1+\log(2+|y|))$ by
$|y-z'|\le(1+|y|)(1+|z'|)$ and local logarithmic integrability
at the subtraction point. Together with $V\le0$ this gives
$|V(y)|=o(|y|^2)$. The source has finite mass by Fatou's lemma.
Thus \cite[Theorem~1.3]{AWZ} at a marked point and
\cite[Theorem~1.1]{LinClass} at an ordinary point give the normalized
profile $U_\gamma$. Its center is the mark in the former case and
the chosen maximum in the latter. The singular
profile has its unique maximum there by Lemma~\ref{lem:concavity}.
The limiting mass is $M_s$ and hence exhausts this atom by
Lemma~\ref{lem:virial}. Rescaled densities converge in local
$L^{1+\epsilon}$ for $\epsilon>0$ sufficiently small.

We next obtain an integrable majorant on the entire core.
The crude estimate $r_s^4 dF_k/dV_g\le C$ follows at a singular
point from \cite[Thm.~5.5]{AWZ}. At a regular point, if
$r_s e^{v_k(x)}$ were unbounded, maximize it on a smaller fixed
ball. On the boundary and on fixed annuli away from the pole $v_k$
is bounded above. At a maximizing point $z'$ put $d=d_g(z,z')$ and
$\lambda'=e^{v_k(z')}$. Then $d\lambda'\to\infty$, $d\to0$ and
$v_k\le v_k(z')+\log2$ on $B(z',d/2)$. In coordinates rescaled by
$\lambda'$ the density on $B_2$ is uniformly bounded. Subtract the
Green formula at $z'$. On $B_2$ the logarithmic kernel differences
have bounded $L^1$ norms uniformly for observation in $B_1$.
Outside $B_2$ the first-derivative estimate and the total mass bound
give an $O(1)$ contribution. Thus $v_k-v_k(z')\ge-C$ on this
rescaled $B_1$. Since the regular coefficient is bounded below there,
its source mass is bounded below by a fixed positive number.
Also $\lambda\ge c\lambda'$ by the peak definition. The ball
$B(z',1/\lambda')$ therefore misses $B(z,R/\lambda)$ for every fixed
$R$ eventually. Both balls converge to the same atom, whose entire
mass is already exhausted as $R\to\infty$, a contradiction.

Subtract the Green formula at $z$. On $|z'|\le R$ in rescaled
coordinates the logarithmic ratio is
$\log(|z'|/|y-z'|)+O(1)$. For $|y|>2R$ its integral, with the
factor $1/(8\pi^2)$, is bounded above by
$-(2\mu-\eta)\log |y|+C_R$ by exhaustion, for $R$ sufficiently
large. On $R<|z'|<|y|/2$ and $|z'|>2|y|$ the positive part of
this ratio is bounded. On $|z'|\asymp |y|$ it has bounded
integral against a density bounded by $C|y|^{-4}$, by the crude
estimate. Contributions off the coordinate ball are bounded on a
smaller core. Thus, for any small $\eta>0$,
\[
 V_k(y)\le C_\eta-(2\mu-\eta)\log(1+|y|),\qquad
 dF_k(y)\le C |y|^{4\gamma}(1+|y|)^{-8\mu+4\eta}\,dy .
\]
Here $v_k(z)-m=o(1)$ by the centered limit. Choose $\eta<\mu$ on
every core, so the majorant satisfies the hypotheses of \eqref{eq:core-log-potential}.

Let $\sigma_{s,k}=F_k(B_s)$. Evaluation at $z$ now gives
$-\bar v_k=(\sigma_{s,k}/(8\pi^2)-\mu)\log\lambda+O(1)$,
since the rescaled logarithmic moments are uniformly integrable.
In particular all log-scales are comparable to $|\bar v_k|$. The
preceding majorant bounds the integral of $r_s^b$ on $B_s$ by
$C\lambda^{-b}$ for $0<b<\min(1,4\mu-4\eta)$. First distance
moments are consequently $O(e^{-c|\bar v_k|})$. Apply
\eqref{eq:symvirial} with the dilation centered at the current
$z_{s,k}$. In $B_s$ the errors in $A_X=4\mu+O(r_s)$ and
$B_X(x,x')=-2+O(r_s(x)+r_s(x'))$ follow from first variation and
normal coordinates. On cross-pairs with another core the kernel is
$O(r_s(x))$, and pairs meeting the exterior give $O(\tau_k)$.
More explicitly, the individual identity reads
\[
 4\mu_s\sigma_{s,k}-\sigma_{s,k}^2/(4\pi^2)
 =O\!\left(\int_{B_s}r_s\,dF_k+\tau_k\right).
\]
The derivative of the left side at $M_s=16\pi^2\mu_s$ is
$-4\mu_s\ne0$. Since $\sigma_{s,k}\to M_s$, the first-moment
bound gives $|\sigma_{s,k}-M_s|=O(e^{-c|\bar v_k|})$ as asserted.

Retain the constant in the Green formula at $z$. The log moment tends
to $\int\log |y|\,d\nu_\gamma(y)$, and $-2\mu$ times this moment
is $U_\gamma(0)-c_\gamma=-c_\gamma$ by the entire representation.
The mass error times $\log\lambda$ tends to zero. Thus
\[
 v_k(z)-\bar v_k=2\mu\log\lambda-c_\gamma+
 \mu H(z,z)+\sum_{t\ne s}\mu_t G(z_t,z)+o(1).
\]
The peak definition proves \eqref{eq:matching-identity}. At any core point,
\[
 (v_k-\bar v_k)(\exp_z(y/\lambda))
 =\frac{\sigma_{s,k}}{8\pi^2}\log\frac{\lambda}{1+|y|}+O(1).
\]
This follows from \eqref{eq:core-log-potential} for the rescaled
source; metric distance comparison on a normal ball changes its
logarithm by $O(1)$. Use a slightly larger core for this estimate,
so that the remaining source has bounded potential; changing the
fixed core radius changes its mass by $O(\tau_k)$.
For $1\le|y|\le r_0\lambda_s$, the mass error
multiplied by $\log\lambda_s$ tends to zero. Exponentiating the
last display and using $K_D=r_s^{4\mu_s-4}\widehat K_s$ gives
$dF_k/dV_g\le C\tau_k r_s^{-4-4\mu_s}$ there. For $|y|<1$,
the peak normalization gives a rescaled density bounded by
$C|y|^{4\mu_s-4}$, while matching bounds
$\tau_k\lambda_s^{4\mu_s}$ above and below. The same inequality
follows because $|y|^{4\mu_s-4}\le |y|^{-4-4\mu_s}$.
This proves the density estimate in \eqref{eq:exacttail} on the
whole core. Since $L_k=-\bar v_k-\tfrac14\log\kappa_*$, the summed
core mass estimates and the $O(\tau_k)$ exterior mass give the
first bound in \eqref{eq:peak-energy}. Matching gives
$\mu_s\log\lambda_s=L_k+O(1)$, so the core formula yields
$u_k=2L_k+O(1+\log(1+|y|))$ there. The error has bounded integral
against $F_k$ by the core majorant. Since $u_k$ is bounded on the
exterior and $F_k(M)=2\kappa_k$, testing against $u_k$ gives
$\langle P_gu_k,u_k\rangle=\int u_k\,dF_k=4\kappa_k L_k+O(1)$.
The modulation argument that follows is uniform for this coefficient
family. Indeed the corrected source depends on $h_k$ through
the first jet of its logarithm. At fixed marked poles those jets
are bounded, and differentiation in the separated ordinary centers
acts only on the Green interactions.

Take $\alpha_s^0=\sigma_{s,k}/M_s$ and keep the projection cutoffs
of \eqref{eq:cutoff-measure}. Pull $F_k|_{B_s}$ and the full transported
probability $\nu_s^M$ back to rescaled coordinates, extending both
densities by zero to $\R^4$. Denote them by $\widetilde F_{s,k}$
and $\nu^{\mathrm{res}}_{s,k}$. Their difference
$\widetilde F_{s,k}-\sigma_{s,k}\nu^{\mathrm{res}}_{s,k}$ has zero
mass even when the cutoff extends beyond $B_s$. With $A_\ell=\{2^\ell<|y|<2^{\ell+1}\}$, the estimates
above and \eqref{eq:matching-identity} give, for some $\delta_s>0$,
\begin{equation}
 \sup_{A_\ell}\frac{d\widetilde F_{s,k}}{dy}
 +\sup_{A_\ell}\frac{d(\sigma_{s,k}\nu^{\mathrm{res}}_{s,k})}{dy}
 \le C\,2^{-4\ell}2^{-\delta_s|\ell|}.
\end{equation}
For $\ell\ge0$ use \eqref{eq:exacttail} and matching, and for $\ell<0$
use $v_k\le m$ in the rescaled density. The model satisfies the same
bound on its entire support. More explicitly, one may take
$\delta_s=2\mu_s$: on negative shells the density is
$O(|y|^{-4+4\mu_s})$, and on positive shells it is
$O(|y|^{-4-4\mu_s})$. Subtract the first-shell mean of a unit
$H^2$ test from its pairing with the zero-mass difference. On
$|\ell|\le L$, the two rescaled densities converge to the same
model in $L^{1+\epsilon}$. By H\"older and \eqref{eq:entry-shell-test},
their pairing with $f_\lambda-m_0$ is $o_L(1)$ uniformly over
$\|f\|_{H^2}\le1$.
The remaining pairing is bounded by
$C\sum_{|\ell|>L}(1+|\ell|)^{3/4}2^{-\delta_s|\ell|}$, including
shells outside $B_s$ where only the model remains. Hence
$F_k|_{B_s}-\sigma_{s,k}\nu_s^M\to0$ in $H^{-2}$. On the exterior
the local $L^q$ bound on $F_k/\tau_k$ gives an $O(\tau_k)$ dual
error, also after subtracting its mean.
Inverting $P_g$ on the mean-zero slice gives strong $H^2$ approximation
by the projected sum with amplitudes $\alpha_s^0$.
Proposition~\ref{prop:ext-size} bounds the difference between this sum
and the extended sum by $o(1)$ in $H^2$.
Lemma~\ref{lem:modulate} changes log-scales by $o(1)$, the amplitudes
by $o((\log\lambda_s)^{-1/2})$, and regular centers by
$o(\lambda_s^{-1})$, preserving \eqref{eq:matching-identity}.
Since $\alpha_s^0-1=O(e^{-c|\bar v_k|})$ and
$\log\lambda_s\asymp|\bar v_k|$, it follows that $\zeta_s\to0$.
Modulation from the preliminary sum also gives $\beta=o(1)$ and
$\|w\|=o(1)$. At the chart scales, \eqref{eq:unnormalized-core}
therefore gives
\[
 F_k(B_s)=2\kappa_k\tau_k L'_s\bigl(M_s/6+o(1)\bigr).
\]
Changing from the peak ball to the chart ball costs $O(\tau_k)$
on a fixed annulus. Since the left side tends to $M_s>0$ and
\eqref{eq:matching-identity} gives $\kappa_k\tau_k L_s\to3$,
$L'_s/L_s\to1$ for each $s$. Taking logarithms in
\eqref{eq:individual-height} yields
$(\alpha_s-1)l_s=o(1)$; all $L_s$ are asymptotically equal.
\end{proof}

For every component set $h_s=\Theta_s(z_s)e^{-4c_s}/3$ and
$\Lambda_s^4=h_s\lambda_s^{4\mu_s}$. Let
$T_0=n^{-1}\sum_s\log\Lambda_s$ for the peak parameters, $n=p+|I|$. Thus
$\kappa_k\tau_k e^{4T_0}\to1$ along exact blow-up sequences.

For a type with largest weight $m=\max_s\mu_s$ write
\begin{equation}\label{eq:Eorder}
 s(T)=s_m(e^T),
\end{equation}
where $s_m$ is the finite, logarithmic or low-weight order given
in the introduction.
\begin{proposition}[Sharp componentwise Green estimates]\label{prop:green-sizes}
For the separated configurations of Section~\ref{sec:chart}, let
$K\subset M$ be compact. If $d_g(K,z_s)\ge\delta>0$, then for
sufficiently large $\lambda_s$ the individual projected bubble obeys
\[
 \|\varphi_s-\mu_sG(z_s,\cdot)\|_{C^3(K)}
 +\|\partial_{l_s}\varphi_s\|_{C^3(K)}\le C_\delta\,\mathfrak e_{\mu_s}(\lambda_s),
\]
with $\mathfrak e_\mu$ given by \eqref{eq:e-mu}, including $\mu_s=1$.
If the distance bound holds for every component and
$L^{-1}\Lambda\le\lambda_s^{\mu_s}\le L\Lambda$ with $L\ge1$ fixed, then
\begin{equation}\label{eq:green-size-total}
 \left\|\sum_s\varphi_s-\Phi_{A,I}\right\|_{C^3(K)}
       \le C_{\delta,L}s_m(\Lambda).
\end{equation}
\end{proposition}

\begin{proof}
Green representation gives
$\varphi_s(x)-\mu_sG(z_s,x)=\mu_s\int(G(y,x)-G(z_s,x))\,d\nu_s^M(y)$.
For each individual pole $z_s$ insert a radial cutoff supported in
$B_{\delta/2}(z_s)$ and equal to one in a smaller ball. On its support
$\nabla_x^mG(y,x)$, $0\le m\le3$, has uniformly bounded source
derivatives for $x\in K$. Apply Theorem~\ref{thm:moment}, or the finite
second-moment calculation at a regular pole, to this cutoff test.
Outside the smaller ball the density is $O(\lambda_s^{-4\mu_s})$;
the observation derivatives of orders zero through three are locally
integrable uniformly in $x$ by
Section~\ref{sec:green}. The cutoff normalization has the same error.
For the log-scale derivative the signed radial source has zero mass,
vanishing linear moments and the same bounds, as in
\eqref{eq:core-dilation}. The same exterior Taylor estimate applies.
Under the two-sided scale bound, summing the individual estimates and
substituting $\lambda_s\asymp\Lambda^{1/\mu_s}$ gives \eqref{eq:green-size-total}.
\end{proof}

\begin{lemma}[The self-pair kernel at second order]\label{lem:metric-pair}
In a fixed normal chart at $z$, let
$X(x)=ax+(b\cdot x)x-|x|^2b/2$ near zero. Put
$d=|x|+|y|$. Uniformly for $x\ne y$ in a sufficiently small ball,
\begin{equation}\label{eq:kernel-second-order}
 (X_x+X_y)G(x,y)
 =P_2(x,y)+\frac{2a}{3}
       \frac{\operatorname{Rm}_z(x,y,y,x)}{|x-y|^2}
       +o(d^2),
\end{equation}
where $P_2$ is a polynomial of degree at most two. Its constant term
is $-2a$. The small order is uniform up to the diagonal and uniform
for $z$ and $(a,b)$ in the compact sets used here.
\end{lemma}
\begin{proof}
Let $q=|x-y|^2$. The normal-coordinate expansion of squared
distance \cite[(9)]{Pennec} has degree-four Taylor polynomial
$q-\operatorname{Rm}_z(x,y,y,x)/3$ in our curvature convention.
For uniformity as $x\to y$, put $v=(x+y)/2$, $h=x-y$. The smooth
remainder $R$ vanishes with its first $h$-derivative on $h=0$.
Taylor's formula in $h$ gives
\[
 R=h_i h_j E_{ij}(v,h),\qquad
 |E_{ij}|\le Cd^3,\qquad |\nabla E_{ij}|\le Cd^2 ,
\]
since all derivatives of $R$ through order four vanish at $(0,0)$.
In particular, with $\mathcal D=x\cdot\partial_x+y\cdot\partial_y$,
\begin{equation}\label{eq:distance-relative-jet}
 d_g(x,y)^2=q-\tfrac13\operatorname{Rm}_z(x,y,y,x)+R,\qquad
 |R|+|\mathcal D R|\le Cd^3q .
\end{equation}

Since $|\operatorname{Rm}_z(x,y,y,x)|\le Cd^2q$, radial
differentiation of the logarithm in \eqref{eq:distance-relative-jet}
gives
\[
 \mathcal D\log d_g(x,y)
 =1-\frac{\operatorname{Rm}_z(x,y,y,x)}{3q}+O(d^3).
\]
For the quadratic part $S_b$ of $X$, its Euclidean conformal identity
is
$(S_b(x)\cdot\partial_x+S_b(y)\cdot\partial_y)\log|x-y|
=b\cdot(x+y)/2$.
Along the simultaneous flow of $S_b$, $\dot v=O(d^2)$ and
$\dot h=O(d|h|)$, so $|\frac{d}{dt}(h/|h|)|=O(d)$.
The curvature quotient is
$\operatorname{Rm}_z(v,h/|h|,h/|h|,v)$, whose derivative is
$O(d^3)$. The formula for $E_{ij}$ bounds the derivative of
$R/|h|^2$ by $O(d^4)$. These bounds are uniform for $h\ne0$. Multiply by the coefficient $-2$ of $\log d_g$.
Finally $(X_x+X_y)H(x,y)$ has a quadratic Taylor polynomial with
$o(d^2)$ remainder, by the joint $C^2$ regularity established in
Section~\ref{sec:green} and $X(0)=0$. Combining these identities proves
\eqref{eq:kernel-second-order} and its stated constant.
\end{proof}

For the common dilation let $F_*=\sum_sM_s\delta_{z_s}$ at the
current centers, $F=2\kappa\mu_u$ and $\omega=F-F_*$. Put $W=K_D e^{4\sum_s\mu_sG(z_s,\cdot)}$.
Use a smooth field supported in disjoint core balls, equal to
$X=a_sx+(b_s\cdot x)x-|x|^2b_s/2$ near $z_s$, with
$a_s=(4\mu_s)^{-1}$. The functional $V_X$ in \eqref{eq:offshellvir}
vanishes at $F_*$: its single-core terms are
$4\mu_sa_sM_s-2a_sM_s^2/(8\pi^2)=0$, and the cross terms
vanish. Its linear term at $F_*$ is pairing with
$h_X=\operatorname{div}_g X+X\log W$ off the poles, with
$h_X(z_s)=-4\mu_sa_s=-1$. The latter is also the value obtained
from the continuous self-pair limit $-2a_s$. Since
$\omega(M)=2(\kappa-\kappa_*)$, expansion of the quadratic
functional $V_X$ gives, for every such $u$,
\begin{equation}\label{eq:secondvir}
 2(\kappa-\kappa_*)+V_X(F)=\int(1+h_X)\,dF+
       (8\pi^2)^{-1}\iint B_X\,d\omega\,d\omega .
\end{equation}
By \eqref{eq:offshellvir}, $V_X(F)=\langle\mathcal R_\kappa(u),-Z_X\rangle$
and hence vanishes at every exact solution.

\subsection{Coefficients of the common dilation}
Write $k_s=\nabla\log\Theta_s(z_s)$. In an orthonormal normal frame let
\[
 b_s= -a_s k_s/[2(1-\mu_s)]\quad(\mu_s<1),\qquad b_s=0\quad(\mu_s=1),
 \qquad S_b(x)=(b\cdot x)x-|x|^2b/2 .
\]
Define the quadratic polynomial and the bounded-on-unit-pairs kernel
\begin{align*}
 Q_s(x)&=a_s\big(\nabla^2\log\Theta_s[x,x]-\tfrac13
                      \operatorname{Ric}_s[x,x]\big)+k_s\cdot S_{b_s}(x),\\
 D_s(x,y)&=\frac{\operatorname{Rm}_s(x,y,y,x)}{|x-y|^2}.
\end{align*}
The curvature convention makes $\operatorname{Rm}(x,y,y,x)$ positive
on a positively curved two-plane. Values on $x=y$ in integrals against
densities are immaterial; put $D_s(x,0)=D_s(0,y)=0$. For $4\mu_s>2$
write
\[
 J_s=M_s\int Q_s\,d\nu_s+
       \frac{a_s M_s^2}{12\pi^2}
                \iint D_s(x,y)\,d\nu_s(x)d\nu_s(y).
\]
Here $d\nu_s$ is the entire probability in the normalization of
Section~\ref{sec:entire}, also at a regular point. These expressions are independent
of the orthonormal frame. Use $m$ and $s$ from \eqref{eq:Eorder}.
Define $C(A)$ for all separated configurations by
\begin{equation}\label{eq:allC}
 C(A)=
 \begin{cases}
 \displaystyle\frac1{2\kappa_*}
    \sum_{\mu_s=m} h_s^{1/(2m)}J_s,&m>1/2,\\[2mm]
 \displaystyle\sum_{\mu_s=1/2}
      \frac{\Theta_s(z_s)}{\kappa_*\mu_s}
                   \int_{\Sph^3}Q_s(\theta)\,d\theta,&m=1/2,\\[2mm]
 D_I^{(4)}/\kappa_*,&m<1/2.
 \end{cases}
\end{equation}

For the finite- and logarithmic-moment coefficients there is a scalar
form useful when perturbing $h$. Set
\[
 \mathcal L_s=\Delta_g\log\Theta_s(z_s)
   -\frac{2(2-\mu_s)}9 R_g(z_s)-k_s\cdot b_s/a_s .
\]
Then, with $m_{2,s}=\int |x|^2d\nu_s$ in the first case,
\begin{equation}\label{eq:trace-coefficient}
 J_s=\tfrac14 a_sM_s m_{2,s}\mathcal L_s\quad(\mu_s>1/2),\qquad
 \int_{\Sph^3}Q_s=\tfrac14 a_s|\Sph^3|\mathcal L_s\quad(\mu_s=1/2).
\end{equation}
Indeed spherical averaging of $Q_s$ gives the indicated trace with
curvature term $-R_g/3$, already the asserted value at $\mu_s=1/2$.
For the curvature pair at a finite moment use independent $x,y$ of
law $\nu_s$ and $D_0=(|x|^2|y|^2-(x\cdot y)^2)/|x-y|^2$.
The flat dilation identity with $X(x)=|x|^2x$ gives
\[
 (2+4\mu_s)m_{2,s}
 =\frac{M_s}{4\pi^2}\,\mathbb E\bigl(|x|^2+|y|^2+x\cdot y-2D_0\bigr).
\]
This follows by testing the normal logarithmic representation as in
\eqref{eq:offshellvir}; the quotient
$(X(x)-X(y))\cdot(x-y)/|x-y|^2$ equals the expression in parentheses.
Cutoffs at radius $R$ have $|DX_R(z)|\le C|z|^2$, so the finite
second moment permits their removal. Hence
$\mathbb E D_0=(2\mu_s-1)m_{2,s}/(4\mu_s)$.
By rotation invariance the curvature pair is $(R_g(z_s)/12)\mathbb E D_0$,
since the mean sectional curvature in four dimensions is $R_g/12$.
Substitution gives the finite-moment identity in \eqref{eq:trace-coefficient}.

The field $X$ used in \eqref{eq:secondvir} now has a specified quadratic
part. In normal coordinates at a core, $X=a_s x+S_{b_s}(x)$ and
$\log\sqrt{|g|}=-\operatorname{Ric}_s[x,x]/6+O(|x|^3)$.
Since $\operatorname{div}_{\R^4}S_b=4b\cdot x$ and
$S_b\log|x|=b\cdot x/2$, direct expansion gives
\begin{equation}
 1+h_X=(2(1-\mu_s)b_s+a_s k_s)\cdot x+Q_s(x)+o(|x|^2).
\end{equation}
In particular the linear term vanishes at each active mark. For pure
low-weight types $W(1+h_X)$ is integrable. Integration of
$h_XW=\operatorname{div}_g(WX)$ outside small balls gives
\begin{equation}\label{eq:finite-flux}
 \int_M(1+h_X)W\,dV_g=D_I^{(4)} .
\end{equation}
Indeed the inward boundary flux at $q_j$ is
$-|\Sph^3|\Theta_j(q_j)a_j d_j^{-4\mu_j}
+O(d_j^{2-4\mu_j})$: the linear terms average to zero on each
sphere. These are exactly the subtractions in \eqref{eq:D4} for independent
deletion radii $d_j$.

\section{Degree conventions and the first chamber}\label{sec:degree}
Using \eqref{eq:Pplus}, write $\mathcal N_\kappa(u)=2\kappa(\mu_u-dV_g/V)$.
For $1<s<s_*$ use $\mathcal X_s=W^{4,s}(M)\cap\{\int_Mu\dd=0\}$. Define
\[
 \cT_\kappa=(P_g^+)^{-1}(P_g-\mathcal N_\kappa),\qquad
 \widehat{\cT}_\kappa=I-P_g^{-1}\mathcal N_\kappa .
\]
By Section~\ref{sec:green},
$K=P_g^{-1}\mathcal N_\kappa:\mathcal H\to \mathcal X_s$ is compact on
bounded sets. For the inclusion $j:\mathcal X_s\hookrightarrow
\mathcal H$, the fixed-point maps on the two spaces are $I-jK$ and
$I-Kj$. Here $K$ and $j$ are compact on bounded sets; compactness of $j$
follows from $s>1$ and Rellich embedding. Commutativity of the fixed point index \cite[Section~8]{Granas}
identifies their
degrees on isolating neighborhoods of the same compact zero set.
On $\mathcal H$, $\cT_\kappa=\tfrac12\nabla_{P^+}J_\kappa$ by
\eqref{eq:energy}. In either space it is obtained from the fixed-point map by
postcomposition with $(P_g^+)^{-1}P_g$, the identity on positive modes
and minus the identity on the $\bar k$ negative modes. Write $d$ for
the degree of $\cT_\kappa$; the degree of $\widehat{\cT}_\kappa$ is
$(-1)^{\bar k}d$. Nonresonant compactness in
Lemma~\ref{lem:compact} makes both degrees constant in positive chambers.

\begin{lemma}[Base degree]\label{lem:base-degree}
For positive parameters below the first wall $d=(-1)^{\bar k}$.
\end{lemma}
\begin{proof}
It suffices to include zero in the parameter continuation. For a
mean-zero solution $u=\rho\int G(\cdot,y)d\mu_u(y)$,
$\rho=\kappa/(8\pi^2)$. For fixed $a>1$, Jensen and the uniform
logarithmic bound on $G$ give $\int e^{4au}dV_g\le C_a$ when
$8a\rho<4$ with a fixed margin. The normalization integral is bounded
below by $(\inf K_D)V$ by Jensen. Choose $1<s<s_*$, then $b\in(s,s_*)$
and $a$ large enough that H\"older bounds $K_De^{4u}$ in $L^s$.
For sufficiently small $\kappa$ this gives a common $\mathcal X_s$ bound.
Combine with positive nonresonant compactness and continue to zero,
where the linear map is the identity with $\bar k$ signs reversed.
\end{proof}
\section{Gradient expansions and the degree across a wall}\label{sec:faces}
We follow the reduction and degree calculation of
Ahmedou--Hu--Wang \cite[Propositions~3.1--3.7 and Lemmas~4.1--4.2]{AHW}.
The gradient is tested in the normal, individual dilation, amplitude
and ordinary translation directions of the corrected sum. The
logarithmic heights use the balanced coordinates of
Ahmedou--Ndiaye \cite[Section~5]{ANboundary}. The relative and common
dilations have different orders. The latter requires a second
expansion, on the whole block of the fixed-parameter homotopy.

Fix a type, put $n=p+|I|$, $\theta_s=M_s/(2\kappa_*)$, and fix a
small closed ball of separated configurations with one local ordering.
The estimates on this ball do not require criticality. For the local
degree it will isolate one nondegenerate critical orbit. For $g\in\R^n$ let
$P_0g=g-\mathbf1\sum g_s/n$. Set $B=\sum\beta_jv_j$ and use the
natural amplitudes and actual logarithmic core heights
\begin{equation}\label{eq:newbalance}
 \begin{split}
 l_s&=\log\lambda_s,\qquad \zeta_s=(\alpha_s-1)\sqrt{l_s},\\
 H_s&=\log\Lambda_s+2\mu_s\zeta_s\sqrt{l_s}+B(z_s),\\
 T&=n^{-1}\sum_sH_s,\qquad \widehat r=P_0H .
 \end{split}
\end{equation}
Thus $H_s=T+\widehat r_s$. The variables
$\zeta,\beta,\widehat r$ and $t=(\kappa-\kappa_*)/\kappa_*$ range
in fixed small balls, independent of $T$. Since $\log\Lambda_s=\mu_sl_s+\frac14\log h_s$, the box permits
$(\alpha_s-1)l_s=O(\sqrt T)$. The coordinates are
well-defined for large $T$: explicitly,
\begin{equation}\label{eq:natural-inverse}
 \sqrt{l_s}=-\zeta_s+
 \sqrt{\zeta_s^2+
 (T+\widehat r_s-\tfrac14\log h_s-B(z_s))/\mu_s}.
\end{equation}
Uniformly on the box,
$l_s=T/\mu_s+O(\sqrt T)$ and $\alpha_s=1+O(T^{-1/2})$.
The change of coordinates preserves orientation. First replace
$\alpha_s$ by $\zeta_s$ at fixed $l_s$, a positive diagonal change.
At fixed $\zeta$, the derivative of $H_s$ in $l_s$ is
$\mu_s(1+\zeta_s/\sqrt{l_s})>0$; the dependence on $\beta$ and
centers adds only off-diagonal blocks. Finally split $H$ into its
mean-zero and common parts, with the same orientation in source and
target.

The estimates below are first uniform on this fixed box. We sharpen
the common dilation on a smaller box which still contains every
escaping exact solution near the chosen configuration.

\subsection{The normal equation and the transverse gradient}
\begin{lemma}[Quantitative normal graph on the natural amplitude box]\label{lem:fixedgraph}
For a sufficiently small fixed transverse box and all sufficiently
large $T$, the normal equation at $\kappa=\kappa_*(1+t)$ has a
unique solution $\bar w$ in a fixed small normal ball. Write
$a=\max_s|\alpha_s-1|$. For some $c>0$ one has
$\|\bar w\|\le C(a+e^{-cT})$ uniformly on the box. For
$u=\widetilde U+w$ put $\tau=(\int_M K_De^{4u}\dd)^{-1}$,
$\sigma_s=2\kappa\mu_u(B_s)$ on the fixed core balls.
With $\eta=a+\|w\|+e^{-cT}$, both on the graph and for any
$w\in E$ with $\|w\|\le T^{-1/4}$, one has
\begin{equation}\label{eq:boxmass}
 \frac{\tau^{-1}}{\kappa_* e^{4T}}=\sum_j\theta_j e^{4\widehat r_j}+O(\eta),\qquad
 \sigma_s=\frac{2\kappa\theta_s e^{4\widehat r_s}}
                {\sum_j\theta_j e^{4\widehat r_j}}+O(\eta).
\end{equation}
The conditional rescaled densities
converge locally in $L^{1+\epsilon}$ to $d\nu_s/dy$ for some
$\epsilon>0$; their first physical distance moments and the normalized
exterior mass are $O(e^{-cT})$.
\end{lemma}
\begin{proof}
By \eqref{eq:individual-height} and \eqref{eq:newbalance},
$L'_s=3e^{4(T+\widehat r_s)}$. Test the shell comparison in the proof
of \eqref{eq:source-defect} against $1$ before normalization. Since
$l_s\asymp T$, the model mass $M_s/6$ has error $O(a+e^{-cT})$:
the amplitude costs $O(a)$ and the remaining shell errors decay
exponentially in $T$. Multiplication by $e^{4w}$ adds $O(\|w\|)$
by \eqref{eq:weighted-exp}. Thus the core integral is
$e^{4(T+\widehat r_s)}(M_s/2+O(\eta))$. The exterior integral is
$O(1)$ by H\"older and Adams. Summing with $M_s=2\kappa_*\theta_s$
and dividing each core by the total proves \eqref{eq:boxmass}.
The same shell comparison gives an $O(\eta)$ error on conditional
averages of tests whose shell $L^q$ means grow polynomially for all
finite $q$, since the geometric tails absorb these factors.

At $w=0$, the core probabilities differ from $\theta_s$ by
$O(|\widehat r|)+o(1)$. Their conditional laws satisfy
\eqref{eq:weighted-shell} by \eqref{eq:individual-shape}.
Lemma~\ref{lem:weighted-testing} compares their first and second
moments on unit normal tests with the model moments up to $o(1)$.
The normalized exterior costs $O(e^{-4T})$ by H\"older.
Inserting these estimates into the variance
gives, in operator norm on the current normal space,
\[
 \big|8\kappa\operatorname{Var}_{\mu_{\widetilde U}}(f)
      -8\kappa_*\operatorname{Var}_{\eta^0}(f)\big|
  \le [C(|t|+|\widehat r|)+o(1)]\|f\|^2 .
\]
By \eqref{eq:reference-gap}, first choosing the box and then its
lower height bound gives $L=D_w\mathcal N_E(\boldsymbol\xi,0)\ge cI$,
where
$\mathcal N_E(\boldsymbol\xi,w)=\Pi_E(P_g^+)^{-1}
\mathcal R_\kappa(\widetilde U+w)$ and $\boldsymbol\xi$ denotes the
parameters of the sum.

Differentiation of the normalized exponential gives the first nonlinear
variation $-8\kappa\operatorname{Cov}_{\mu_u}(h,f)$ and the second
$-32\kappa\int(h-\bar h)(k-\bar k)(f-\bar f)\,d\mu_u$,
where bars denote $\mu_u$-means. Let $\nu=\mu_{\widetilde U}$.
The shell bounds and Lemma~\ref{lem:weighted-testing} give, uniformly
on $E$ and for every fixed $b,r<\infty$,
$\int |f|^r e^{b|f|/\|f\|}\,d\nu\le C_{b,r}\|f\|^r$.
The fixed exterior is controlled by global Adams and its normalized
mass is exponentially small. For $\|w\|$ in a fixed small ball,
Jensen gives $\int e^{4w}\,d\nu\ge
\exp(4\int w\,d\nu)\ge c>0$. H\"older therefore bounds
$\int |hkf|\,d\mu_u$ by $C\|h\|\|k\|\|f\|$, uniformly.
The same bound applies to all terms produced by centering. Hence the
normal derivative is uniformly Lipschitz.
Write $\eta_0=a+e^{-c_1T}$. The residual at zero is $O(\eta_0)$
by \eqref{eq:source-defect} since $l_s\asymp T$. Uniform Lipschitz
continuity and $L\ge cI$ make $w\mapsto w-L^{-1}\mathcal N_E(\boldsymbol\xi,w)$
a contraction on the ball of radius
$2\|L^{-1}\|\|\mathcal N_E(\boldsymbol\xi,0)\|$.
Indeed, if the Lipschitz constant of the derivative is $C$ and
$\|L^{-1}\|\le c^{-1}$, the derivative of the iteration is bounded
by $C c^{-1}\|w\|$. The radius is $O(\eta_0)$ uniformly, so it
is eventually smaller than $c/(2C)$, and the iteration maps that ball
into itself. Its fixed point satisfies
\begin{equation}\label{eq:quant-normal}
 \|\bar w\|\le C\|\mathcal N_E(\boldsymbol\xi,0)\|\le C\eta_0 .
\end{equation}
Positivity on the fixed small normal ball gives uniqueness there.
The implicit function theorem in the modulation frame gives $C^1$ dependence.
Integration of the positive normal derivative yields on the same ball
\begin{equation}\label{eq:normal-during}
 \big\langle (1-h)\mathcal N_E(\boldsymbol\xi,w)+hw,w\big\rangle
 \ge c\|w\|^2-C(1-h)\eta_0\|w\|,\qquad 0\le h\le1 .
\end{equation}

Finally, on each shell multiplication by $e^{4w}$ preserves
\eqref{eq:weighted-shell} with a smaller $p>1$ and decay exponent,
by H\"older and the annular exponential bound in the proof of
Lemma~\ref{lem:weighted-testing}. On fixed annuli $e^{4w}\to1$ in
every finite averaged $L^q$. At the pole choose $1<p'<p$ with
$4(1-1/p')$ smaller than the shell decay exponent to sum the
ordinary $L^{p'}$ norms. This proves the local density convergence.
For a sufficiently small $b>0$, summation of the physical
$b$-moment costs $O(\lambda_s^{-b})$ and bounds the first distance
moment on the fixed ball. The normalized exterior has exponentially
small mass by H\"older and global Adams.
\end{proof}

Write $R=\mathcal R_\kappa(\widetilde U+w)$ and retain the normal row
$\mathcal N_E(\boldsymbol\xi,w)$. The finite rows below are pairings of $R$; the
dilation and translation tests are those of Lemma~\ref{lem:virial}.

\begin{proposition}[Transverse equations and their orientation]\label{prop:transrows}
For $u=\widetilde U+w$ with $w\in E$, $\|w\|\le T^{-1/4}$, a
continuous orientation-preserving choice of tests, together with the
normal row, gives the following components at the prescribed parameter,
uniformly as $T\to\infty$. Put
$\eta=\max_s|\alpha_s-1|+\|w\|+e^{-cT}$:
\begin{align}
 N_j&=-\vartheta_j\beta_j+O(\eta),&
 B_{\alpha,s}&=2\mu_sM_s\zeta_s+O(\eta/\sqrt{l_s}),\label{eq:transfirst}\\
 B_A&=-\nabla_A\cF_{p,I}(A)+O(e^{-cT}),&
 B_0&=-t+O(e^{-cT}),\notag\\
 D_{\rm rel}&=-(1+t)\,P_0\!\left(
       (e^{4\widehat r_s})_s/\textstyle\sum_j\theta_j e^{4\widehat r_j}
                  \right)+O(\eta).\label{eq:translast}
\end{align}
Here $B_0$ is the $\theta$-mean and $D_{\rm rel}$ the arithmetic
mean-zero part of the normalized individual scale rows. The common
test appropriate at the smaller order $s(T)$ is specified below.
\end{proposition}
\begin{proof}
For each ordinary core minimize
$I_i(z)=\int_{B_i}d_g(z,x)^2\,dF(x)$ on a fixed small closed convex
normal ball about the chart center $a_i$. Choose the radii so that
$\nabla_z^2d_g(z,x)^2\ge c g$ throughout for $x\in B_i$.
Then $\nabla^2 I_i\ge c\sigma_i g$. The conditional rescaled law
has a radial limit and its first moment on a shell of radius $R\ge1$
is bounded by $CR^{-3+o_T(1)}\exp(C(1+\log R)^{3/4})$, by
\eqref{eq:individual-shape} and Lemma~\ref{lem:fixedgraph}.
The summable tail and the formula
$\nabla I_i(a_i)=-2\int_{B_i}\exp_{a_i}^{-1}(x)\,dF(x)$ give
$|\nabla I_i(a_i)|=o(\sigma_i/\lambda_i)$. On the boundary the
positive quadratic term in the Taylor inequality dominates this
linear term. The unique minimizer $z_i^\circ$ is interior and
$d_g(a_i,z_i^\circ)\le|\nabla I_i(a_i)|/(c\sigma_i)=o(\lambda_i^{-1})$.
In particular $\int_{B_i}\exp_{z_i^\circ}^{-1}(x)\,dF(x)=0$.
The density varies continuously in $L^1$, and the boundaries of the
moving balls have volume zero. Hence the integrals defining $I_i$
vary uniformly on the closed balls in local coordinates. Uniqueness
gives continuous dependence of $z_i^\circ$.

Set $z_s^\circ=z_s$ at the marks. In the normal charts at $z_s^\circ$,
multiply $x+S_{b_s}(x)/a_s$ by a fixed radial cutoff equal to one
near the origin and supported inside the core to obtain $Y_s$.
Use the coefficients of $Q_s$ from
Section~\ref{sec:rate-low} evaluated at the recentered configuration.
In particular $b_s=0$ at ordinary centers. The individual dilation rows satisfy
\begin{equation}
 \mathfrak d_s=(4\pi^2/\sigma_s)V_{Y_s}(F)
              =M_s-\sigma_s+O(e^{-cT}).
\end{equation}
Indeed the self-pair gives $-2+O(|x|+|y|)$ and the one-point term
$4\mu_s+O(|x|)$, other atom pairs have zero limiting constant, and
the first-moment and exterior bounds of Lemma~\ref{lem:fixedgraph}
control every error. If a core radius is changed within a fixed annulus,
its mass changes by $O(e^{-cT})$.

For the translation field $X_{ij}$ take the constant coordinate vector
$e_j$ near the chart center $a_i$, with a cutoff inside the core.
Here $\nabla X_{ij}(a_i)=0$ and $\mathcal L_{X_{ij}}g=O(r)$, so the
logarithmic self-pair error is bounded by the sum of the two distances
to $a_i$. Expanding at the chart centers $z_s$ gives
\[
 \frac{V_{X_{ij}}(F)}{\sigma_i}
 =\partial_j\log\widehat K_i(a_i)
   +\frac{\sigma_i}{4\pi^2}\partial_{2,j}H(a_i,a_i)
   +\sum_{s\ne i}\frac{\sigma_s}{4\pi^2}\partial_{2,j}G(z_s,a_i)
   +O(e^{-cT}).
\]
The error follows from the first distance moments and the exterior
mass bound of Lemma~\ref{lem:fixedgraph}. Take its negative and
subtract linear combinations of $\mathfrak d_s$ to replace $\sigma_s$
by $M_s$. Then \eqref{eq:gradF} gives $B_A$.

With the normal row retained, the new tests have the orientation
of the parameter derivatives. In $H^2$
$-Z_{Y_s}$ differs by $o(1)$ from
$4(\sigma_s/M_s)\partial_{\log\lambda_s}\widetilde\varphi_s$,
and $Z_{X_{ij}}/\lambda_i$ differs by $o(1)$ from
$4(\sigma_i/M_i)\partial_{a_{ij}}\widetilde U/(\alpha_i\lambda_i)$.
For the radial test this comparison is differentiation of the
transported measure against the velocity $-x$ at fixed center; for
the center test the velocity is $e_j+O(r)$ in its normal coordinates.
By $H^2\hookrightarrow W^{1,4}$ it suffices to estimate
$X(F-\sigma_s\nu_s^M)$ in $L^{4/3}$, with $X/\lambda_i$ for a
translation. A change of dilation center of $o(1/\lambda_i)$ costs
$o(1)$ times the normalized translation bound; the quadratic field
adds a factor $O(r)$ on a physical shell of radius $r$. Cutoff
changes on fixed annuli cost $O(e^{-cT})$.
On fixed rescaled annuli the source difference tends to zero in
$L^{4/3}$ by the local regularity of the model density and the
averaged exponential bound for $w$. On a shell of rescaled radius $R$, writing $m_R$
for the mass majorant $C\min(R^{4\mu_s},R^{-4\mu_s+o(1)})$,
the $L^{4/3}$ norm for a dilation is at most
$C m_R\exp(C(1+|\log R|)^{3/4})$: density scale $r^{-4}$,
vector scale $r$, and the volume factor $r^3$ cancel. For a
normalized regular translation there is one extra $R^{-1}$, summable
also at zero since $\mu_i=1$. On the exterior the cutoff fields are
supported away from inactive marked points and the density is
exponentially small in $L^{4/3}$. Euler derivatives of the cutoff are supported on a fixed annulus
where the model density is $O(\lambda_s^{-4\mu_s})$. Differentiating
its normalizing constant adds $O(\lambda_s^{-4\mu_s})$ times the model
probability. A term of this
mass times the first-shell mean is still exponentially small in $T$:
\eqref{eq:entry-shell-test} on a fixed outer annulus gives
$|m_0|\le C(1+\log\lambda_s)^{3/4}$ for a unit $H^2$ test.
Proposition~\ref{prop:ext-size} bounds the additional derivatives of
the extension, including its dependence on the regular centers.
Inversion of $P_g$ on mean-zero functions gives the two test comparisons.
To check orientation, project each replacement onto the tangent span
and use the uniform Gram bound of Lemma~\ref{lem:modulate}. After the indicated positive rescalings the coefficient matrix
is $I+o(1)$ in the bounded normalized Gram frame. Straight interpolation
therefore remains invertible. The normal parts of the new tests add
bounded linear combinations of the retained normal component, which
are triangular row operations.

The negative-mode test gives
\[
 \langle R,v_j\rangle-\sum_s \mathfrak d_s v_j(z_s)
                           =-\vartheta_j\beta_j+O(\eta)=:N_j .
\]
Indeed a smooth test differs from its center value in each average
by $O(e^{-cT})$, by the distance moments. The factors $\alpha_s$
in the kinetic averages cost $O(a)$ after subtracting the mass rows;
$w$ is orthogonal to the tangent and negative modes.
For the amplitude test $\widetilde\varphi_s$ put $c_{sk}=\mu_sG(z_s,z_k)$ for
$k\ne s$ and $c_{ss}=\int(U_s-c_s+\mu_s H(z_s,z_s))\,d\nu_s$.
By \eqref{eq:physical-projection} and Proposition~\ref{prop:green-sizes},
the mean of $\widetilde\varphi_s$ on core $k$ is
$2\mu_s l_s\mathbf1_{s=k}+c_{sk}+O(\eta)$ under the conditional
nonlinear law and the corrected source normalized on that ball. Its missing mass
outside the ball is exponentially small, also after multiplication
by $l_s$. For the self
term, $U_s$ grows at most linearly in the shell index, so the
weighted comparison of Lemma~\ref{lem:fixedgraph} applies.
Variation of the smooth factors costs a first distance moment,
and the signed extension has zero mass on every concentric ball
and exponentially small logarithmic moments by
Proposition~\ref{prop:ext-size}. The smooth cross terms allow
recentring with error $O(e^{-cT})$. Since $w$ is normal and
$\int \widetilde\varphi_s\,dV_g=0$, the kinetic and density pairings therefore give
\[
 \langle R,\widetilde\varphi_s\rangle
 =2\mu_s l_s(\alpha_s M_s-\sigma_s)
 +\sum_k c_{sk}(\alpha_k M_k-\sigma_k)
 -\sum_j\vartheta_j\beta_j\!\int_M v_j\widetilde\varphi_s\,dV_g+O(\eta).
\]
The integrals against $v_j$ are bounded by the Green formula.
Substitute the negative and individual scale rows just computed;
$l_se^{-cT}$ is absorbed by decreasing $c$. This yields
\begin{equation}\label{eq:amplitude-before-division}
 \langle R,\widetilde\varphi_s\rangle-2\mu_s l_s\mathfrak d_s
 -\sum_k c_{sk}\mathfrak d_k-\sum_j N_j\int_M v_j\widetilde\varphi_s\,dV_g
 =2\mu_sM_s(\alpha_s-1)l_s+O(\eta).
\end{equation}
Division by $\sqrt{l_s}$ gives the amplitude row in
\eqref{eq:transfirst}. All subtractions are shears.

Finally divide each $\mathfrak d_s$ by $M_s$ and take the pair formed by its
$P_0$ projection and its $\theta$-weighted mean. Equations
\eqref{eq:boxmass} and $\sum_s\sigma_s=2\kappa+O(e^{-cT})$
give \eqref{eq:translast} using the $O(\eta)$ conditional
mass approximation; for the mean one uses the sum of the unnormalized
masses directly and the exponentially small exterior error. Relative to arithmetic mean in the domain
the change to weighted mean in the target adds a linear combination of
relative rows to the last row. Every other subtraction above is a shear.
This proves the orientation assertion.
\end{proof}

\begin{definition}[Nondegenerate wall data]\label{def:controlled}
At a wall $\kappa_*$ require that for each representation the function
$\cF_{p,I}$ be Morse and that $C(A)\ne0$ at its critical
configurations, with $C$ given by \eqref{eq:allC}. For $p=0$ the
configuration is a point of index zero.
\end{definition}

\subsection{The common-scale gradient at fixed mass}
Let $D_{\rm rel}=P_0(\mathfrak d_s/M_s)_s$ be the relative-scale
row of Proposition~\ref{prop:transrows}. In the virial calculation
write $z_s$ for $z_s^\circ$ of that proposition and use these centers
in $\Theta_s$. For each core write
$Y_s$ for the cutoff field used in $\mathfrak d_s$ and put
\[
 X=\sum_s a_sY_s,\quad a_s=(4\mu_s)^{-1},\qquad
 G_{\rm sc}=\frac{V_X(F)}{2\kappa_*}
            =\sum_s\frac{\sigma_s}{2\kappa_*}\frac{\mathfrak d_s}{M_s}.
\]
Thus $G_{\rm sc}$ is a tested component of the gradient at the prescribed
$\kappa$, by \eqref{eq:offshellvir}. Replacing the $\theta$-mean by this
row preserves orientation: relative components are fixed and the new
row is positive on the constant vector, since $\sum_s\sigma_s>0$.

\begin{proposition}[Refined common-scale gradient]\label{prop:sharp-fixed}
Fix a compact separated center ball and $0<K<\infty$. Put $S=s(T)$.
Let $T\to\infty$, $t\to0$ with arbitrary relative rates, and take
$u=\widetilde U+w$, $w\in E$, with
$|\zeta|+|\beta|+|\widehat r|+\|w\|\le K T^{-3}$.
Writing $\operatorname{Var}_\theta d=\sum_s\theta_s
(d_s-\sum_j\theta_j d_j)^2$, one has uniformly
\begin{equation}\label{eq:fixed-common}
\begin{split}
 G_{\rm sc}={}&-t-t^2-\operatorname{Var}_\theta D_{\rm rel}+C(A)S\\
 &\quad+o\!\left(S+(|t|+|D_{\rm rel}|)\sqrt S\right).
\end{split}
\end{equation}
Here $C(A)$ is given by \eqref{eq:allC} at the chart centers,
also away from a critical configuration. In particular the remainder
in $G_{\rm sc}+t-C(A)S$ is $o(S)$ on
$D_{\rm rel}=0$ if $|t|=O(S)$. This estimate holds throughout the
normal ball specified above, without solving a gradient equation in $E$.
\end{proposition}

\begin{proof}
Write $\tau=(\int K_D e^{4u}\,dV_g)^{-1}$ and
$F_*=\sum_s M_s\delta_{z_s}$ at the recentered points,
$\omega=F-F_*$. The conditions give $a=O(T^{-7/2})$,
$\max_s|\alpha_s-1|l_s=o(1)$ and $\|w\|\le KT^{-3}$.
Equations \eqref{eq:boxmass} and \eqref{eq:newbalance} now give
$\kappa_*\tau e^{4T}\to1$, $\log\Lambda_s=T+o(1)$ and $\sigma_s/M_s\to1$.
Set $W=K_De^{4\sum_s\mu_sG(z_s,\cdot)}$ and
$h_X=\operatorname{div}_gX+X\log W$. Identity \eqref{eq:secondvir} reads
\begin{equation}\label{eq:offcommonidentity}
 2\kappa_*(G_{\rm sc}+t)=\int(1+h_X)\,dF
                    +(8\pi^2)^{-1}\iint B_X\,d\omega\,d\omega .
\end{equation}

On a rescaled shell $A_\ell$ in a core, Lemmas~\ref{lem:tight}
and~\ref{lem:weighted-testing} bound the mean of a unit normal
test by $C(1+|\ell|)^{3/4}$ and its exponential oscillations in
every fixed averaged power. Hence on $0\le\ell\le C_1T$,
\begin{equation}\label{eq:shell-w-small}
 \left(\operatorname{avg}_{A_\ell}|e^{4w}-1|^q\right)^{1/q}
 \le C_q\|w\|(1+\ell)^{3/4}
            \exp\{C_q\|w\|(1+\ell)^{3/4}\}=o(1).
\end{equation}
For negative shells the subexponential factor is absorbed by their
geometrically decreasing model mass. The rescaled core laws thus
converge locally in $L^1$ to $M_s\nu_s$. On shells of physical
radius $r\ge1/\lambda_s$ their mass is bounded by
$C\tau r^{-4\mu_s}$: use \eqref{eq:raw-core-potential},
$|\alpha_s-1|l_s=o(1)$, the bounded regular fields and
\eqref{eq:shell-w-small}. This bound also holds for fixed-radius
annuli after recentering; for a regular point the center displacement
is $o(1/\lambda_s)$. Off the cores $dF/(\tau\,dV_g)\to 2\kappa_* W$
in $L^1$ along convergent centers, taking $W$ at their limits.
This follows from the exterior Green estimates, the uniform smallness
of the extensions, and global exponential convergence for $w\to0$;
H\"older applies at the inactive singularities.

On $1/\lambda_s\ll r\ll1$, divide the density before multiplication
by $e^{4w}$ by $2\kappa\tau r^{-4-4\mu_s}$. Its limit is
$\Theta_s(z_s)$ by \eqref{eq:physical-projection} and \eqref{eq:far-field}:
the self-potential is $-2\mu_s\log r+\mu_sH(z_s,z_s)+o(1)$, the
amplitude error is $o(1/l_s)$, and the correction tends to zero.
The ordinary recentering $o(1/\lambda_s)$ preserves this limit.

At a borderline pole, fix $\epsilon>0$ and choose $R$ large and $d$
small so that the relative error on $R/\lambda_s\le r\le d$ is
at most $\epsilon$ for all sufficiently large $T$. This uniform
choice follows by contradiction from the neck limit above.
The inner and outer regions each contribute $O_{R,d}(\tau)$ to
the quadratic moment. Each intervening dyadic shell contributes
$O(\tau)$, so \eqref{eq:shell-w-small} changes the moment by
$o(\tau\log\lambda_s)$, also with a bounded angular factor.
After division by $\tau\log\lambda_s$, the remaining error has
limsup $O(\epsilon)$. Let $\epsilon\downarrow0$.

Let $e_s=\sigma_s-M_s$, $e=\max_s|e_s|=o(1)$ and $r_s=d_g(x,z_s)$.
Keep $X$ and its cutoffs fixed when varying the integration radius $d$.
Choose $d$ so that the balls $B_s(d)=B_d(z_s)$ are disjoint and lie in
the polynomial portions of $X$. Write $F_s=F|_{B_s(d)}$ and
$E_d=M\setminus\bigcup_s B_s(d)$. Summing the shell masses weighted
by $r_s^2$ and using bounded mass on $r_s\le1/\lambda_s$ gives the
second-moment bounds $C\lambda_s^{-2}$, $C\tau\log\lambda_s$, and
$C\tau$ for $4\mu_s>2$, $4\mu_s=2$, and $4\mu_s<2$, respectively.
All are $O(S)$. In the last case the geometric sum up to
$2^\ell/\lambda_s\le d$ gives the sharper bound
\[
 \int_{B_s(d)}r_s^2\,dF
 \le C\lambda_s^{-2}+C\tau d^{2-4\mu_s},\qquad
 \lambda_s^{-2}=o(\tau).
\]

For fixed $d$, the signed first moment
$p_s(d)=\int\exp_{z_s}^{-1}(x)\,dF_s(x)$ is $o_d(\sqrt S)$,
where $o_d$ takes $T\to\infty$ at fixed $d$.
For a dominant finite second moment this follows by radiality and
uniform integrability after rescaling. At the logarithmic threshold
the absolute first moment is $O(\lambda_s^{-1})$. A core with second
moment $o(S)$ is covered by Cauchy--Schwarz. In the pure low-weight
case apply Cauchy--Schwarz on $B_s(\delta)$ with the sharper bound
above. The remaining annulus has mass $O_\delta(\tau)$.
Divide by $\sqrt S$, pass to the limit in $T$, then let
$\delta\downarrow0$. At an ordinary center the signed first moment
vanishes on the ball defining the center of mass, so changing its
radius costs only $O_d(\tau)$.

In each individual scale test expand the one-point term, self-pair
and separated pairs to first order at the current centers.
Lemma~\ref{lem:metric-pair} bounds the self-pair remainder by the
sum of squared distances, also for the quadratic part of $Y_s$.
The constant terms give precisely $M_s-\sigma_s$ after normalization,
and the integrated linear terms are $o(\sqrt S)$ by the preceding
bounds. The exterior contributes $O_d(\tau)$. Consequently
\[
  \mathfrak d_s=M_s-\sigma_s+o(\sqrt S),\qquad
  e\le C\bigl(|D_{\rm rel}|+|t|\bigr)+o(\sqrt S).
\]
For the second inequality sum the masses:
$\sum_s e_s=2\kappa_*t+O_d(\tau)$. Indeed, subtracting the
$\theta$-mean from $\mathfrak d_s/M_s$ gives
\[
 e_s/M_s=t-(D_{\rm rel})_s+\sum_j\theta_j(D_{\rm rel})_j
               +o(\sqrt S).
\]
Here $\tau=O(S)=o(\sqrt S)$.

To evaluate the double integral in \eqref{eq:offcommonidentity}, set
$\omega_s=F_s-M_s\delta_{z_s}$ and
$e_s(d)=\int d\omega_s=e_s+O_d(\tau)$. For $y\in E_d$, subtract
$B_X(z_s,y)$ in the integral against $\omega_s$. On $B_s(d/2)$ the
kernel has bounded first derivative in $x$ by separation. On the
remaining annulus $B_X$ is bounded and the $F$-mass is $O_d(\tau)$.
Thus, uniformly in $y$,
\[
 \left|\int B_X(x,y)\,d\omega_s(x)-e_s(d)B_X(z_s,y)\right|
 \le C_d\int_{B_s(d/2)}r_s\,dF+C_d\tau
 \le C_d(\sqrt S+\tau).
\]
Since $F(E_d)=O_d(\tau)$, a mixed pair costs at most
$C_d\tau(e+\sqrt S+\tau)=o_d(S)$ even when $\tau\asymp S$.
Two exterior factors cost $O_d(\tau^2)$.

On two cores use the quadratic Taylor polynomial in their normal
coordinates. For distinct cores the kernel is smooth. On a self-pair
Lemma~\ref{lem:metric-pair} adds $(2a_s/3)D_s$ to the polynomial.
In both cases the remainder is bounded by
$\eta(d)(r_s^2+r_j^2)$, with $\eta(d)\to0$ uniformly on the center
set. Its integral against $|\omega_s|\otimes|\omega_j|$ is at most
$C\eta(d)S$ by the second moments and bounded total variations.
The constant integrates to $-2a_s e_s(d)^2$ on a self-pair and
vanishes on distinct cores since $X(z_s)=0$. A linear monomial gives a
component of $p_s(d)e_j(d)$ or $e_s(d)p_j(d)$. A pure quadratic costs
$O(S(|e_s(d)|+|e_j(d)|))$, and a mixed quadratic is a product of first
moments, hence $o_d(S)$. Replacing $e_s(d)^2$ by $e_s^2$ costs
$O_d(e\tau+\tau^2)$. Since $D_s$ vanishes when either argument is zero,
\begin{align*}
 \iint B_X\,d\omega\,d\omega
   &=-2\sum_s a_s e_s^2
      +\frac23\sum_s a_s\iint D_s\,dF_s\,dF_s+R_{T,d},\\
 \lim_{d\downarrow0}\limsup_{T\to\infty}
       \frac{|R_{T,d}|}{S+e\sqrt S}&=0 .
\end{align*}
These bounds apply along every sequence in the stated box. The
constant term contributes $-\sum_s e_s^2/(2\kappa_*M_s)$ to $G_{\rm sc}$.
In the single integral the linear term vanishes at a mark by the
choice of $b_s$. At an ordinary center its integral is $O_d(\tau)$
by recentering, and then $m=1$ and $\tau=o(S)$. The quadratic term
is $Q_s$ of Section~\ref{sec:rate-low}, with integrated remainder
bounded by $o_{d\downarrow0}(1)S$.

For $m>1/2$ only the cores of weight $m$ survive. Rescale by
$\lambda_s$ on each. The single integral gives $M_s\int Q_s\,d\nu_s$
in units of $\lambda_s^{-2}$, by the uniform second-moment tail.
For the pair kernel use $|D_s(x,y)|\le C\min(|x|^2,|y|^2)$ and local
$L^1$ convergence; pairings with its central atom vanish.
Including the factor in \eqref{eq:offcommonidentity} gives exactly
$J_s$ in these units. Since $\lambda_s^{-2}/S\to h_s^{1/(2m)}$
along convergent centers and the exterior term is $O_d(\tau)=o_d(S)$,
division by $2\kappa_*$ gives the first line of \eqref{eq:allC}.

For $m=1/2$ the exterior mass and the second moments at the lighter
cores are $O_d(\tau)=o_d(S)$. On a dominant core the rescaled tail
probability is $O((1+R)^{-2})$ uniformly. If $Y,Y'$ are independent with this
rescaled law, layer cake and independence give
\[
 \mathbb E\min(|Y|^2,|Y'|^2)
 =\int_0^\infty 2R\,\mathbb P(|Y|>R)^2\,dR\le C.
\]
Since $|D_s(x,y)|\le C\min(|x|^2,|y|^2)$, the double curvature
term is $O(\lambda_s^{-2})=o(S)$. The quadratic single integral has coefficient
$2\kappa_*\Theta_s(z_s)\int_{\Sph^3}Q_s$ in units of
$\tau\log\lambda_s$, by the neck estimate. Here
$\log\lambda_s/T\to1/\mu_s$, proving the second line.

For $m<1/2$ all centers are marked and $1+h_X=O(r_s^2)$ there.
The truncated second-moment bound above gives
$C\tau\sum_s d^{2-4\mu_s}+o_d(\tau)$ for the absolute contributions
of the inner single integral and the curvature pairs. The same powers
of $d$ bound the integrals of $|(1+h_X)W|$ on the deleted balls.
On $E_d$ the exterior $L^1$ convergence of $F/\tau$ gives the limit
of the single integral there. Thus the full single integral divided
by $2\kappa_*\tau$ tends to $\int(1+h_X)W\,dV_g=D_I^{(4)}$ by
\eqref{eq:finite-flux}, and $\tau/S\to1/\kappa_*$.
In all three cases first take $T\to\infty$ and then $d\downarrow0$.
Continuity gives the same $C(A)S$ up to $o(S)$ at the chart centers.
The remaining error is $o(S+e\sqrt S)$. To evaluate the constant
self-pair, put $d_s=(D_{\rm rel})_s$ and $\bar d=\sum_s\theta_sd_s$.
The individual mass formula reads $e_s/M_s=t-(d_s-\bar d)+o(\sqrt S)$.
Square and sum with weights $\theta_s=M_s/(2\kappa_*)$.
Since $\sum_s\theta_s(d_s-\bar d)=0$, this gives
\[
 \sum_s\frac{e_s^2}{2\kappa_*M_s}
 =t^2+\operatorname{Var}_\theta D_{\rm rel}
    +o\big(S+(|t|+|D_{\rm rel}|)\sqrt S\big).
\]
Together with the bound on $e$, this proves \eqref{eq:fixed-common}
uniformly on the stated box.
\end{proof}

\subsection{Localization and sharp rates for exact zeros}
\begin{lemma}[Transverse localization of exact zeros]
\label{lem:uniform-exact-rate}
In the fixed transverse box and fixed small normal ball of
Lemma~\ref{lem:fixedgraph}, every sequence of exact zeros with
$T\to\infty$ has $|t|+|\zeta|+|\beta|+|\widehat r|+\|w\|=O(e^{-cT})$
for some $c>0$. Its peak and chart heights satisfy $T_0=T+o(1)$.
\end{lemma}
\begin{proof}
On the normal graph \eqref{eq:quant-normal} gives
$\|w\|=O(a+e^{-cT})=O(T^{-1/2})$, so
Proposition~\ref{prop:transrows} applies with $\eta=O(a+e^{-cT})$.
At a zero each tested row vanishes. Equation
\eqref{eq:amplitude-before-division} and $l_s\asymp T$ give
$Ta\le C(a+e^{-cT})$. Hence $a$, $|\zeta|$ and $\|w\|$ are
exponentially small, after decreasing $c$. The negative and mean
dilation rows give the same bound for $\beta$ and $t$.
The relative row has leading term $-4(1+t)\widehat r+
O(|\widehat r|^2)$ on $\mathbf1^\perp$; absorbing the quadratic term
gives $|\widehat r|=O(e^{-cT})$.

Equation \eqref{eq:boxmass} gives $\kappa_*\tau e^{4T}\to1$.
The core laws identify the type and the amplitude Gram estimate gives
divergence in $H^2$. For its peak scales, \eqref{eq:matching-identity}
gives $\kappa\tau e^{4T_0}\to1$. Their ratio yields $T_0-T=o(1)$.
\end{proof}

\begin{theorem}[The mass defect of an exact solution]\label{thm:allrates}
Along every exact escaping sequence with limiting configuration $A$,
$t=C(A)s(T_0)+o(s(T_0))$, with the coefficient in \eqref{eq:allC}.
\end{theorem}
\begin{proof}
By Theorem~\ref{thm:strongentry}, height matching and
$(\alpha_s-1)l_s=o(1)$ place the sequence in the fixed box.
Lemma~\ref{lem:uniform-exact-rate} places it in the shrinking box and
gives $T_0-T=o(1)$. At an exact solution the common and relative
rows vanish, so \eqref{eq:fixed-common} gives
$t+t^2=C(A)s(T)+o(s(T)+|t|\sqrt{s(T)})$.
Absorption first gives $t=O(s(T))$, then $t=C(A)s(T)+o(s(T))$.
In each regime $s(T_0)/s(T)\to1$. Uniformity in a closed ball
isolating a critical orbit follows by subsequence extraction and
the center equation.
\end{proof}

\begin{theorem}[Exact rate at a pure low-weight wall]\label{thm:low-rate}
For every escaping exact sequence of pure singular type $I$ with
$0<\mu_j<1/2$ at all its poles,
\begin{equation}\label{eq:lowactual}
 \kappa_k-\kappa_*=\kappa_*D_I^{(4)}\tau_k+o(\tau_k).
\end{equation}
Equivalently, the coefficient of $e^{-4T_0}$ in
$t=(\kappa_k-\kappa_*)/\kappa_*$ is $D_I^{(4)}/\kappa_*$.
If $D_I^{(4)}\ne0$, then
$\operatorname{sgn}(\kappa_k-\kappa_*)=\operatorname{sgn}D_I^{(4)}$
eventually. Under this condition there is no escaping sequence of
type $I$ at the fixed parameter $\kappa_*$.
\end{theorem}
\begin{proof}
Substitute $C(A)=D_I^{(4)}/\kappa_*$ and
$s(T_0)=e^{-4T_0}\sim\kappa_*\tau_k$ in Theorem~\ref{thm:allrates}.
When $D_I^{(4)}\ne0$, the leading term in \eqref{eq:lowactual} has
its sign, excluding escape at $\kappa=\kappa_*$.
\end{proof}

\subsection{The homotopy and the jump}
\begin{theorem}[Degree of an escaping configuration]\label{thm:faces}
Suppose Definition~\ref{def:controlled} holds, and let $A_*$ be an
unordered critical configuration of index $m_*$. For all sufficiently
small $t=(\kappa-\kappa_*)/\kappa_*\ne0$, the local degree of
$\cT_\kappa$ on the block of its escaping solutions on the side $t C(A_*)>0$ is
\[
       (-1)^{\bar k+n+m_*}\operatorname{sgn}t .
\]
On the opposite side this configuration admits no escaping sequence.
\end{theorem}
\begin{proof}
Put $C=C(A_*)$. On the indicated side use $T_-<T_+$ given by
$s(T_-)=2|t/C|$, $s(T_+)=|t/C|/2$. Both tend to infinity, and their
difference stays bounded by \eqref{eq:Eorder}. Use transverse balls
of radius $\delta=T_-^{-3}$ for $\zeta,\beta,\widehat r$ and a
fixed sufficiently small coordinate ball for $A$ containing just
$A_*$. For fixed $t$ the truncated tube is bounded. Throughout it
$|t|\le K s(T)$ and $\delta\le K T^{-3}$ for a fixed $K$.

Let $\mathcal U_t$ be the image of the strict base inequalities
and $\|w\|<\delta$ under $(\xi,w)\mapsto \widetilde U(\xi)+w$,
$w\in E_\xi$. This map is injective on the closed block for small
$|t|$. Otherwise two representations of the same $u$ along $t\to0$
have reference sums within $2\delta$ in $H^2$. On the closed block
$\mu_s l_s=T+O(1)$ by \eqref{eq:newbalance}, so corresponding
scales have bounded ratios. Test the source difference against smooth
bumps on the first set of rescaled cores, with bounded $H^2$ norm.
An unbounded rescaled center displacement would give zero mass
from the second representation on every fixed ball of that core,
contradicting the positive model mass there. For bounded displacements,
the two limiting model laws agree. Uniqueness of their centers and
scales makes the normalized center and log-scale increments $o(1)$.
The normalized amplitude increments and the negative-mode increments
are $o(1)$ by the transverse bounds. Lemma~\ref{lem:modulate} then
identifies the two representations.

Extend the base coordinates across the faces. The equations
$\langle u-\widetilde U(\xi),\partial_{\xi_j}\widetilde U(\xi)
\rangle_{P^+}=0$ have derivative in $\xi$ equal to the negative
Gram matrix plus $O(\|w\|)$ in the frozen normalized coordinates
of Lemma~\ref{lem:modulate}. The implicit function theorem yields
a local inverse on the closed block and makes $\mathcal U_t$ open.
If images $u_k$ of the closed block converge, extract $\xi_k\to\xi$
in the compact base. Then $w_k=u_k-\widetilde U(\xi_k)$ converges
to a point of $E_\xi$ of norm at most $\delta$, by continuity of the
orthogonal projections. Thus the image of the closed block is closed,
and injectivity gives a continuous inverse there. In particular
$\partial\mathcal U_t$ lies in the image of the coordinate faces.

On the closed block put $\mathcal B=(N,B_\alpha,B_A)$ and
\[
 \mathcal B_*(\beta,\zeta,A)=
 \big((-\vartheta_j\beta_j)_j,\ (2\mu_sM_s\zeta_s)_s,\ -\nabla_A\cF_{p,I}(A)\big).
\]
The actual tested rows are evaluated at $u=\widetilde U+w$, at the
fixed mass $\kappa$. For $0\le h\le1$ prescribe all components by
\begin{equation}\label{eq:fixed-wall-homotopy}
 \left(\begin{array}{c}
 (1-h)\mathcal N_E(\boldsymbol\xi,w)+hw\\
 (1-h)\mathcal B+h\mathcal B_*\\
 D_{\rm rel}\\
 (1-h)G_{\rm sc}+h\{Cs(T)-t\}
 \end{array}\right).
\end{equation}
In a second interval keep the first, second and last rows at their
endpoints and replace the relative row by
$(1-h)D_{\rm rel}-4h\widehat r$. The initial map represents
the gradient by Proposition~\ref{prop:transrows}. The terminal map
is a product on the parameter box and its normal ball.
Lift the finite rows to the Hilbert space by
$L_\xi b=J_\xi(J_\xi^*J_\xi)^{-1}b$, where $J_\xi$ is the derivative
of $\widetilde U$ in $(\beta,\zeta,A,\widehat r,T)$ and the adjoint uses
the $P^+$ metric. Add the normal row. Its kinetic term is $w$ by
orthogonality to the amplitude and negative modes. The lifted homotopy
minus $u$ is the sum of $-\widetilde U(\xi)$, the lifted finite rows
and the projected nonlinear term. The finite rows are bounded by
Proposition~\ref{prop:transrows}, \eqref{eq:offshellvir} and the
positive lower bound for $\sigma_s$. The first two terms range over
relatively compact sets by compactness of the base and continuity of
$\widetilde U$ and $L_\xi$. The last has relatively compact range by
Section~\ref{sec:green} and continuity of the normal projections on
the compact base. This gives a compact perturbation of the identity.
Proposition~\ref{prop:transrows}
deforms the initial tests by invertible interpolations, shears and positive
rescalings. The height change \eqref{eq:natural-inverse} has positive
determinant, so on the contractible base its matrix is homotopic to the
identity through invertible matrices. The initial degree is that of $\cT_\kappa$.

At a zero of the first homotopy \eqref{eq:normal-during} gives
$\|w\|\le C(a+e^{-cT})=o(\delta)$ since $a=O(\delta/\sqrt{T_-})$.
In the second homotopy $w=0$. Thus $\eta=o(\delta)$ at any possible
zero. By \eqref{eq:transfirst}--\eqref{eq:translast} the negative-mode
and amplitude rows differ from their linear models by $o(\delta)$,
and the relative row is $-4\widehat r+o(\delta)$. This excludes
their respective spheres. The center sphere is excluded by the positive
minimum of $|\nabla_A\cF_{p,I}|$ there.

In the first homotopy $D_{\rm rel}=0$ at a zero, so
Proposition~\ref{prop:sharp-fixed} gives
$G_{\rm sc}=C(A)s(T)-t+o(s(T))$. Shrink the center ball so that
$|C(A)-C|<|C|/16$. At $T_-$ and $T_+$ the model row is $t$ and
$-t/2$, respectively. The error $|C|s(T)/16+o(s(T))$ preserves
both signs under interpolation. In the second homotopy the common
row already equals its model. This excludes every boundary zero.

The endpoint has a single zero, at $w=0$, $\beta=\zeta=\widehat r=0$,
$A=A_*$ and $Cs(T)=t$. Its derivative is the identity on $E$.
On the tangent span, using $J_\xi$ to identify both copies with
parameter space, it is $(J_\xi^*J_\xi)^{-1}$ times the derivative of
the finite model. Derivatives of the moving lift multiply the zero
finite row and hence drop out. The Gram determinant is positive.
The negative modes give $(-1)^{\bar k}$, the amplitudes $+1$, the
center block $(-1)^{m_*}$ since its dimension is $4p$, and the
relative heights $(-1)^{n-1}$. The remaining derivative $Cs'(T)$
has sign $-\operatorname{sgn}C$. The classical local
Leray--Schauder index of this invertible derivative equals the sign
of its finite determinant. As $tC>0$, the product is
$(-1)^{\bar k+n+m_*}\operatorname{sgn}t$. On the other side
Theorem~\ref{thm:allrates} excludes escape at this configuration.
\end{proof}

\begin{theorem}[Jump between adjacent chambers]\label{thm:jump}
Assume Definition~\ref{def:controlled} at $\kappa_*$. Let $d_-$ and
$d_+$ be the degrees of the $|P_g|$-preconditioned operator $\cT$
in the chambers immediately below and above $\kappa_*$. Then
\[
 d_+-d_- = \sum_{(p,I),[A]} J_A,\qquad
 J_A=(-1)^{\bar k+p+|I|+\morse(\cF_{p,I},A)}.
\]
The sum includes every representation $(p,I)$ of $\kappa_*$ and
every unordered critical configuration $[A]$ of that type. For the
fixed-point operator $\widehat{\cT}$ the same jump formula holds
with $J_A=(-1)^{p+|I|+\morse(\cF_{p,I},A)}$.
\end{theorem}
\begin{proof}
Write $\kappa=\kappa_*(1+t)$. At $t=0$ the nonvanishing coefficients
and Theorem~\ref{thm:allrates} bound all solutions. Enclose them in
an open $H^2$ ball $B_R$. Compactness of the nonlinear map on bounded
sets excludes boundary zeros for small $|t|$, hence the degrees on
$B_R$ agree at the two endpoints.

There are finitely many critical orbits by Lemma~\ref{lem:proper}
and the Morse assumption. Choose the center balls with disjoint
closures in the unordered configuration space for each type. The
corresponding sets of limiting atomic measures are disjoint, also for
distinct types, since ordinary and marked centers stay separated.
The blocks are therefore disjoint for all sufficiently small $|t|$:
if two met along $t\to0$, testing $P_gu$ on smooth functions and
passing to a subsequence would give the same limiting atomic measure
for both. The amplitude Gram estimate places all blocks outside
$B_R$. Any missed sequence of exterior zeros with $t\to0$ must have
an unbounded subsequence by compactness. Theorem~\ref{thm:strongentry} places it
near a critical orbit, Lemma~\ref{lem:uniform-exact-rate} gives
transverse size $O(e^{-cT})$, and Theorem~\ref{thm:allrates} gives
$s(T)/|t/C(A)|\to1$ on the prescribed side. These are precisely the
strict block bounds of Theorem~\ref{thm:faces}.

By excision an orbit with $C(A)>0$ adds $J_A$ to $d_+$, while one with
$C(A)<0$ adds $-J_A$ to $d_-$. Each therefore contributes $J_A$ to
the jump. Section~\ref{sec:degree} supplies the factor $(-1)^{\bar k}$
for the fixed-point convention.
\end{proof}
\subsection{Antipodal concentration on the round sphere}
On the unit round $\Sph^4$ take antipodal marks $p,q$ of equal
strength $\gamma=\mu-1$, with $0<\mu<1$, and constant $h>0$.
Use stereographic coordinates $p=0$, $q=\infty$ and $s=|x|$.
The chordal Green function in \eqref{eq:green} is
\[
 G(x,y)=-2\log|x-y|+\log(1+|x|^2)+\log(1+|y|^2)+c .
\]
Indeed conformal covariance reduces the equation off the pole to
$\Delta^2\log(1+|x|^2)=-96(1+|x|^2)^{-4}$, the flat singularity
has mass $16\pi^2$, and the mean fixes $c$. Thus
$K_D=B s^{4\gamma}(1+s^2)^{-4\gamma}$ for a constant $B>0$,
and $dV_g=16(1+s^2)^{-4}dx$. With $U_{\gamma,\lambda}$ from
\eqref{eq:scale} put
\[
 F_\lambda=6s^{4\gamma}e^{4U_{\gamma,\lambda}}dx,\qquad
 u_\lambda=(16\pi^2)^{-1}\int G(\cdot,y)\,dF_\lambda(y).
\]
For fixed $\lambda$, $dF_\lambda/dV_g=O_\lambda(d_g(\cdot,z)^{4\gamma})$
near either mark $z$, by the origin bound and \eqref{eq:far-field}. Thus
$u_\lambda\in W^{4,p}$ for $1<p<-1/\gamma$.
The entire normal representation gives $u_\lambda=
U_{\gamma,\lambda}+\mu\log(1+s^2)+c_\lambda$.
Since $4\mu-4\gamma-4=0$, normalization gives
$\mu_{u_\lambda}=F_\lambda/(16\pi^2\mu)$.
Hence $u_\lambda$ solves \eqref{eq:normalized} at the fixed wall
$\kappa_*=8\pi^2\mu$. Dilation gives
$F_\lambda\rightharpoonup16\pi^2\mu\,\delta_p$, with $q$ inactive.

For $\mu<1/2$, the same formula for $G$ gives
$K_De^{4\mu G(p,\cdot)}dV_g=B' s^{-4-4\mu}dx$ with $B'>0$.
The exterior integral in \eqref{eq:D4} is
$B'|\Sph^3|\tan(r/2)^{-4\mu}/(4\mu)$ for geodesic radius $r$.
Its leading coefficient is $2^{4\mu}B'|\Sph^3|/(4\mu)$ and the
remainder is $O(r^{2-4\mu})$. Thus $D_{\{p\}}^{(4)}=0$, as required
by \eqref{eq:lowactual} along this fixed-parameter family.

This example also realizes both signs of the mass defect. Fix
$0<\mu<1/2$ and take $h_\epsilon=h e^{\epsilon v}$ with
$0\le v\in C^\infty(M)$, $v\not\equiv0$, supported away from $p,q$.
For $z=p,q$ the finite-part counterterm at $z$ is unchanged. Subtracting
\eqref{eq:D4} at the constant $h$ gives the exact identity
\begin{equation}\label{eq:antipodal-variation}
 D_{\{z\}}^{(4)}(h_\epsilon)
   =\int_M(e^{\epsilon v}-1)K_D e^{4\mu G(z,\cdot)}\,dV_g ,
\end{equation}
where $K_D$ is the unperturbed weight. Both coefficients have the sign of
$\epsilon$ for $\epsilon\ne0$. The only representations of this wall
are the two singletons, each with center index zero. For fixed
$\epsilon\ne0$, Theorems~\ref{thm:low-rate} and~\ref{thm:faces} therefore
give a solution in each singleton block for every sufficiently close
$\kappa$ with $\epsilon(\kappa-\kappa_*)>0$. As $\kappa\to\kappa_*$
these solutions concentrate at their respective poles. The opposite side admits no bubbling. Each pole contributes $-1$ to
the fixed-point jump by Theorem~\ref{thm:jump}.

For the unperturbed data the local coefficient can also be computed
directly. In geodesic distance $r$ from the active pole,
$K_D$ is proportional to $(\sin r)^{4\gamma}$ and
$H(p,x)=H(p,p)+r^2/12+O(r^4)$. Hence
$\log\Theta_p(x)=\log\Theta_p(p)+(2-\mu)r^2/3+O(r^4)$.
Thus $k_p=0$ and $\Delta_g\log\Theta_p(p)=8(2-\mu)/3$.
Since $R_g=12$, the curvature term in \eqref{eq:trace-coefficient}
cancels this trace. Both the finite-moment and the borderline
coefficients vanish, in agreement with the fixed-wall family.

\section{Poincar\'e--Hopf evaluation and the singular generating series}

Put $X=M\setminus\Sigma$ and write $\Conf_p(X)$ for the ordered
$p$-tuples of distinct points of $X$, with $\Conf_0(X)$ a point.

\begin{lemma}[Ends of the center configuration space]\label{lem:proper}
$\cF_{p,I}$ tends to $+\infty$ on divergent configurations and has no
critical sequence escaping to an end.
\end{lemma}
\begin{proof}
Each vanishing distance $d$ contributes $c\log(1/d)$.
For an unordered regular pair, $c=8$. For a regular point approaching
$q_j$, $c=4(1-\mu_j)$ if $j\notin I$ and $c=4(1+\mu_j)$ if $j\in I$.
These coefficients follow from \eqref{eq:Hamiltonian} and are all
positive. Sum all these terms if several distances vanish.
For a putative escaping critical sequence pass to limits of all centers.
Choose one limiting cluster with at least two regular points or one
regular point at a puncture. In a normal chart dilate all its regular
points about the fixed limiting point, keeping the other clusters fixed.
For a pair within the cluster the derivative of $\log d$ for this
variation is $1+o(1)$, even for nested scales. Indeed the Lie derivative
of the metric under the radial vector field is $2g+o(1)$ uniformly along
each minimizing segment in the shrinking chart. The same holds for a
puncture in the cluster. Interactions with other clusters contribute
$o(1)$. The regular parts have bounded first derivatives in these
coordinates. The derivative of $\cF$ tends to minus the positive sum of
the internal coefficients, contradicting criticality.
\end{proof}
\begin{proposition}[Configuration-space index sum]\label{prop:PH}
If $\cF_{p,I}$ is Morse, then
\begin{equation}\label{eq:PH}
 \sum_{A\in\Crit(\cF_{p,I})}
 (-1)^{\morse(\cF_{p,I},A)}
 =\chi_c(\Conf_p(X))
 =\prod_{l=0}^{p-1}\bigl(\chi(M)-N-l\bigr),
\end{equation}
where the empty product is $1$.
\end{proposition}

\begin{proof}
The case $p=0$ is immediate. Otherwise Lemma~\ref{lem:proper} gives
a compact sublevel at a regular value $L$ above all critical values.
Every larger finite level band is compact and free of critical points.
By \cite[Theorem~3.1]{Milnor}, inclusion into every larger sublevel
induces a homology isomorphism. Every singular cycle and bounding
chain has compact support, so the same holds on $\Conf_p(X)$. Poincar\'e--Hopf for its outward gradient gives the
ordinary Euler characteristic as the sum of the Hessian signs.

The configuration spaces have finite CW type by compatible
triangulation of the diagonals and marked strata in a finite
polyhedron \cite[Theorem~9.2.1]{bcr}. Forgetting the last of
$p\ge2$ points is the Fadell--Neuwirth fibration
\cite[Theorem~3]{FN}, with fiber $X$ minus $p-1$ points.
The fiber has Euler characteristic $\chi(M)-N-p+1$: compactly
supported Euler characteristic is additive, and it equals the
ordinary one on even-dimensional manifolds of finite type by
Poincar\'e duality over $\mathbb F_2$. Euler characteristic is
multiplicative for this fibration (equivalently, take the Euler
characteristic of the Serre spectral sequence with its local
coefficients). Starting with $\Conf_1(X)=X$ gives \eqref{eq:PH},
using the same duality in dimension $4p$.
\end{proof}

\begin{proof}[Proof of Theorem~\ref{thm:main-degree}]
First assume Definition~\ref{def:controlled} at each wall below
$\kappa$. By the binomial theorem the coefficient in
\eqref{eq:generating} is
\begin{equation}
 b_\sigma=\sum_{p+\sum_{j\in I}\mu_j=\sigma}
 \frac{(-1)^{p+|I|}}{p!}\prod_{l=0}^{p-1}(\chi(M)-N-l).
\end{equation}
The sum ranges over $p\in\mathbb N_0$ and subsets $I$ of the marked
indices. Only finitely many pairs occur below any fixed exponent.
Each unordered configuration has $p!$ ordered lifts. Thus
\eqref{eq:PH} and Theorem~\ref{thm:jump} give $b_\sigma$ for the
fixed-point degree jump across $8\pi^2\sigma$. Summing from the
first chamber, where the fixed-point degree is one by
Lemma~\ref{lem:base-degree}, gives \eqref{eq:generating}. For $N=0$
the partial binomial sum is $\prod_{l=1}^m(l-\chi(M))/m!$ in
$m<\rho<m+1$.

To remove the wall hypotheses, perturb $h$ to $he^v$ as in
\cite[Section~7]{ANboundary}. For the finitely many types below
$\kappa$, the cluster dilation in Lemma~\ref{lem:proper} confines
all critical centers to compact separated sets uniformly for $v$
small in $C^1$. The variation of the gradient at a configuration
is the collection of jets $\nabla v(a_i)$. At each configuration,
cutoff affine functions supported in disjoint balls span these jets.
The spanning property persists in a neighborhood. Compactness and
finiteness of the types give finitely many such functions whose span
$V$ works simultaneously on every relevant center set.

For each type with $p>0$, the map
$(A,v)\mapsto\nabla_A\mathcal F_{p,I}^{he^v}(A)$ on
$\operatorname{Conf}_p(M\setminus\Sigma)\times V$ is therefore
transverse to the zero section over the compact set. Its zero set is
a finite-dimensional $C^1$ manifold of dimension $\dim V$. Sard's
theorem applied to its projection onto $V$ gives arbitrarily small
$v$ for which every zero is nondegenerate. Intersect the finitely
many resulting sets of regular values. When $p=0$ there is no center
condition.

There are now finitely many critical configurations. Restrict the next
variation to functions vanishing to first order at their centers and at
all marks. Nonsingularity and compactness keep precisely these critical
configurations for small $C^2$ variations. For a finite or logarithmic
moment choose a bump at one dominant pole with prescribed Hessian trace,
supported away from the other poles. By \eqref{eq:trace-coefficient}
and \eqref{eq:allC} the coefficient variation is a positive multiple of
that trace. For a low-weight type, a nonzero nonnegative bump away
from all poles has variation $\int vK_De^{4\Phi_I}\dd>0$ in the
finite part. Thus each zero coefficient has a proper kernel in the
variation space. Choose a direction outside their finite union and
perturb slightly. This makes all coefficients nonzero and preserves
those already nonzero. Lemma~\ref{lem:compact} bounds the solutions
along the path back to $h$ at the fixed nonresonant mass, so the degree
formula follows by homotopy invariance.
\end{proof}

\begin{proof}[Proof of Corollary~\ref{cor:nonzero-chambers}]
A type $(p,I)$ has $p+|I|$ points and cost $p+\sum_{j\in I}\mu_j$.
In \textup{(i)} every type with at most $k$ points has cost at most $k$.
Among types with $k+1$ points the least cost is attained by taking the
$\min(k+1,N)$ smallest marked weights and filling the remaining places
with ordinary points. This gives the upper endpoint in \textup{(i)}.
Deleting a point strictly decreases the cost. Thus $\rho$ is nonresonant and
the sum in \eqref{eq:generating} includes exactly the types with at
most $k$ points. Recording the number of points by $y$, their signed
coefficients in each degree are those of
$(1-y)^{\chi(M)-N}(1-y)^N=(1-y)^{\chi(M)}$.
Summing through degree $k$ gives
$\sum_{m=0}^k(-1)^m\binom{\chi(M)}m=(-1)^k\binom{\chi(M)-1}k$.
By Section~\ref{sec:green}, $2\rho\le\chi(M)$ for both choices of
$\kappa$, since the Weyl term is nonnegative and $\sum_j\gamma_j\le0$.
Hence $\chi(M)>2k$ and the binomial coefficient is positive.

In \textup{(ii)} the only possible wall exponents below $\rho$
are zero and the individual weights $\mu_j<\rho$. Each such mark
contributes $-1$, including when several weights coincide. The strict
inequalities exclude every wall at $\rho$, giving $1-r$ by
\eqref{eq:generating}. In both cases the existence assertion follows
from Theorem~\ref{thm:main-degree} and the lifts in Section~\ref{sec:green}.
\end{proof}

\section{Existence from a one-sided barycenter map}
\subsection{Compactness and the negative eigenspace}\label{sec:red}
On $\mathcal H$ write $Q(u,v)=\langle Pu,v\rangle$ and $Q(u)=Q(u,u)$.

\begin{lemma}[Compactness for perturbed equations]\label{lem:compact}
Fix the metric, marked points and weights. Let $\kappa_n\to\kappa>0$
with $\kappa/(8\pi^2)\notin\mathcal W$, let $h_n\to h>0$ in $C^1$
with $h_n>0$, and let the mean-zero backgrounds $b_n$ converge in
$C^1$. Use $\mu_{n,u}$ from \eqref{eq:normalized} with $h_n$ in
the weight. Mean-zero solutions of
\[
 P_g u_n=2\kappa_n(\mu_{n,u_n}-V^{-1}dV_g)+b_n\,dV_g
\]
are precompact in $W^{4,s}$ for $1<s<s_*$. The bound is uniform
on compact sets of such nonresonant parameters and data.
\end{lemma}
\begin{proof}
Use the mass gauge with $K_{D,n}=(h_n/h)K_D$. The normalization
in Lemma~\ref{lem:virial} gives coefficient $c_n\kappa_n h_n e^{4f}$
and background $\kappa_n/V-b_n/2$. A subsequence of $c_n$ has a
positive finite limit since $h_ne^{4f}$ is uniformly bounded above and
away from zero. Thus the data converge in $C^1$ and the limiting
coefficient is positive. In Theorem~5.1 of \cite{AWZ}, concentration
would give $\kappa/(8\pi^2)\in\mathcal W$, and uniform decay would
contradict unit weighted mass. The regular parts, and hence $u_n$,
are bounded. Elliptic estimates and compact embedding first give
subsequential uniform convergence; $K_{D,n}\to K_D$ in $L^s$ then
gives convergence of the right sides in $L^s$ and strong $W^{4,s}$
convergence by mean-zero inversion. Compact families follow by extraction.
\end{proof}

Write $\mathcal H=E^+\oplus E^-$ for the $L^2$-orthogonal spectral splitting of $P$ and set $\|w\|_+^2=Q(w)$, $\|z\|_-^2=-Q(z)$ on the respective spaces. Their sum of squares defines an equivalent Hilbert norm on $\mathcal H$. Write $\Pi_\pm$ for the spectral projections, omitting constants. By \eqref{eq:energy} and differentiation of the exponential,
$D\log Z(u)[v]=4\int v\,d\mu_u$ and
$D^2\log Z(u)[v,w]=16\operatorname{Cov}_{\mu_u}(v,w)$.
The Hilbert gradient of $\log Z$ is compact on bounded sets by
the $L^s$ compactness of the densities in Section~\ref{sec:green}.

\begin{lemma}[Reduction]\label{lem:red}
For $\kappa>0$, $z\mapsto J_\kappa(w+z)$ on $E^-$ has a unique maximum $z(w)$ for every $w\in E^+$. The map $z$ is $C^1$ and uniformly bounded. The functional
\begin{equation}
 F_\kappa(w)=J_\kappa(w+z(w))=\|w\|_+^2-K(w)
\end{equation}
is $C^2$ with $\nabla K$ compact. Write $F_\kappa^a=\{w:F_\kappa(w)\le a\}$. Its critical points lift bijectively to solutions of \eqref{eq:normalized}, and
\begin{equation}\label{eq:bounded}
 J_\kappa(w)\le F_\kappa(w)\le J_\kappa(w)+C_\kappa .
\end{equation}
\end{lemma}
\begin{proof}
On the negative space the Hessian is at most $-2\|\cdot\|_-^2$ by the covariance formula. The bound $|\log Z(w+z)-\log Z(w)|\le4\|z\|_\infty$ gives coercivity from above. At the unique maximum,
\[
 Q(z(w),e)=2\kappa\int e\,d\mu_{w+z(w)}\qquad(e\in E^-).
\]
Since $E^-$ is finite dimensional and smooth, this yields $\|z(w)\|_-\le C\kappa$. The implicit function theorem gives $C^1$ dependence and \eqref{eq:bounded} follows. Differentiation in $w$ gives $DK(w)[v]=4\kappa\int v\,d\mu_{w+z(w)}$ for $v\in E^+$. The compactness just proved applies since $z(w)$ is bounded. The projected Euler equations in the two spectral spaces together give \eqref{eq:normalized}.
\end{proof}
\begin{lemma}[Sublevel deformation]\label{lem:lucia}
At a positive nonresonant parameter, if $[a,b]$ contains no critical value of $F_\kappa$, then $F_\kappa^a$ is a deformation retract of $F_\kappa^b$.
\end{lemma}
\begin{proof}
Apply \cite[Proposition~1.1]{lucia} to
$I(t,w)=(\|w\|_+^2-tK(w))/2$ at $t=1$, keeping the reduction at
$\kappa$ fixed. Lemma~\ref{lem:red} gives $K$ a compact, locally
Lipschitz gradient. At a critical point of $I(t,\cdot)$, the two
projected equations for $u=w+z(w)$ give
\[
 Pu=2t\kappa(\mu_u-V^{-1}dV_g)+2(1-t)\kappa\Pi_-\mu_u .
\]
Since $\Pi_-\mu_u=\sum_r v_r\int v_r\,d\mu_u$ is uniformly bounded
in $C^1$, the added background tends to zero for $t_n\le1$,
$t_n\to1$. Lemma~\ref{lem:compact} and $W^{4,s}\hookrightarrow H^2$
for $s>1$ give a convergent subsequence of the lifted critical points.
Both $I(t_n,w_n)$ and $I(1,w_n)$ tend to $F_\kappa(w)/2$ at the
limiting critical point $w$. This lies outside $[a/2,b/2]$, verifying the
exclusion in Lucia's proposition.
\end{proof}

\subsection{Low sublevels and weighted barycenters}
Fix a positive nonresonant $\kappa=8\pi^2\rho$ and write $F=F_\kappa$. We use the reference probability $m_D$ of \eqref{eq:weight} on $M$.

\begin{lemma}[Localized inequality and concentration]\label{lem:low}
Let $U_i\Subset V_i$ be finitely many regions with the $\overline{V_i}$ pairwise disjoint and each containing at most one mark. Set $c_i=1$ if $\overline{V_i}\cap\Sigma=\varnothing$ and $c_i=\mu_j$ if $\overline{V_i}\cap\Sigma=\{q_j\}$. For every $\epsilon>0$ and $\eta>0$ there is $C$ depending on the regions and these data such that
\begin{equation}\label{eq:improved}
 8\pi^2\big(\textstyle\sum_i c_i-\epsilon\big)\log Z(w)
 \le \|w\|_+^2+C
 \quad\text{if }\mu_w(U_i)\ge\eta\text{ for all }i.
\end{equation}
For every weak neighborhood $\mathcal U$ of $B_\rho$ in the probability measures, $F(w)\le-L$ implies $\mu_w\in\mathcal U$ once $L$ is sufficiently large.
\end{lemma}
\begin{proof}
Choose smooth cutoffs $0\le\chi_i\le1$ supported in $V_i$, equal to one on $U_i$. Split $w=l+v$ into finitely many low positive modes and their spectral complement. For any $\delta>0$ the cutoff in the spectrum can be chosen so that
\[
 \|v\|_{H^1}\le\delta\|w\|_+,\qquad
 \|\Delta_g v\|_2^2\le(1+C\delta^2)\|w\|_+^2 .
\]
The first estimate follows from compact embedding of the unit ball in $H^1$ and orthogonality to increasingly many modes; the second follows from the lower-order term in $Q$ and $Q(v)\le Q(w)$. The low part obeys $\|l\|_\infty\le C_\delta\|w\|_+$. The mass bound gives
\[
 \log Z(w)\le 4\|l\|_\infty+C_\eta+
                \log\int_{U_i} e^{4\chi_iv}\,dm_D .
\]
At a regular region use boundedness of the density and set $a_i=1$. At $q_j$ use H\"older with conjugate exponents $a_i$ and $a_i/(a_i-1)<-1/\gamma_j$, choosing $1/a_i$ arbitrarily close to $\mu_j$ from below. By \eqref{eq:adams},
\[
 \log Z(w)\le C+ C_\delta\|w\|_++\frac{a_i}{8\pi^2}\|\Delta_g(\chi_i v)\|_2^2 .
\]
The average of $\chi_iv$ is absorbed in the linear term. Multiply by $8\pi^2/a_i$ and sum. By the product rule and disjointness the sum of the squared Laplacians is at most $(1+\theta)\|\Delta_gv\|_2^2+C_\theta\|v\|_{H^1}^2$ for any $\theta>0$. First choose the exponents, then $\theta$, then $\delta$, and absorb the linear term by Young's inequality. This gives \eqref{eq:improved}, using the uniform lower bound on $\log Z$.

If $F(w_n)\to-\infty$, take a weak limit $\nu$ of $\mu_{w_n}$. If $\nu\notin B_\rho$, its support contains distinct points whose total weight exceeds $\rho$. For infinite support use $\omega_0>0$; for finite support use nonresonance. Take disjoint balls about those points as above, with positive limiting mass in each smaller open ball. The mass criterion holds eventually. Choose $\epsilon$ smaller than the weight excess. Substitution of \eqref{eq:improved} into $F(w_n)\ge J_\kappa(w_n)$ bounds these energies below. This contradiction proves the assertion.
\end{proof}

If $\rho<\omega_0$, the same H\"older and spectral-splitting argument on all of $M$ with $1/a$ close to $\omega_0$ gives $8\pi^2(\omega_0-\epsilon)\log Z(w)\le\|w\|_+^2+C_\epsilon$. Thus $F$ is coercive. It is weakly lower semicontinuous: for each fixed $z$, $w\mapsto J_\kappa(w+z)$ has this property by the exponential bounds of Section~\ref{sec:green}, and $F$ is the supremum over $z$. A minimum gives the first assertion of Theorem~\ref{th:exist} in this range. If there is no negative space, $F=J_\kappa$.

Suppose now $B_\rho\ne\varnothing$ and set $k=\lfloor\rho/\omega_0\rfloor$.
Every admissible barycenter has at most $k$ support points. The possible costs below $\rho$ have a maximum $b<\rho$. The space $B_\rho$ is compact: on passing to limits, points can merge, amplitudes can disappear, and the weight at a marked point can drop; none increases the cost. More explicitly, take subsequences in the compact parameter space $M^k\times\Delta^{k-1}$ and fix the marked incidences along the sequence.

Triangulate the marked manifold as a finite polyhedron in $\R^d$.
Coordinate moments through degree $2k-1$ give a continuous map
$m_*$ on probabilities that separates measures with at most $k$
atoms: on the union of two supports one can isolate each point by
a product of at most $2k-1$ affine factors. Thus $m_*$ embeds the
compact $B_\rho$. Its image is semialgebraic by polynomial
parametrization on the finitely many simplex and incidence strata.
By compatible triangulation \cite[Theorem~9.2.1]{bcr} it is a
neighborhood retract of a surrounding box (after barycentric
subdivision, normalize the coordinates on the full subcomplex).
Composing that retraction with $m_*$ and its inverse on $B_\rho$
gives a retraction $r$ from a weak neighborhood $\mathcal U$.

\begin{proposition}[Test functions with a homotopy left inverse]\label{prop:test}
For every sufficiently large $L$ there are continuous maps
\[
 \Phi:B_\rho\longrightarrow F^{-L},\qquad
 \Psi:F^{-L}\longrightarrow B_\rho,\qquad \Psi\Phi\simeq\mathrm{id}_{B_\rho}.
\]
\end{proposition}
\begin{proof}
Set $\Psi(w)=r(\mu_w)$ by Lemma~\ref{lem:low}. Choose
$m_D(B(x,r_T(x)))=e^{-4T}$ and let $\nu_{T,x}$ be $m_D$
conditioned on that ball. Ball mass is continuous in center and
radius by dominated convergence, strictly increasing below the
injectivity radius, and uniformly positive at any fixed small radius.
Thus for $T$ large the radius is unique and continuous.
Sphere boundaries have zero volume, so the resulting densities are
continuous in $x$ in $L^s$, $1<s<s_*$, by dominated convergence.
Volume comparison and \eqref{eq:weight} give
\begin{equation}\label{eq:radii}
 r_T(x)\le C e^{-T}\quad(x\in M),\qquad
 r_T(q_j)\le C e^{-T/\mu_j}.
\end{equation}
For $\sigma=\sum_i t_i\delta_{x_i}$ put $\nu=\sum_i t_i\nu_{T,x_i}$ and
\[
 \Phi_T(\sigma)=2\kappa G_+\nu,\qquad
 G_+(x,y)=G_0(x,y)+\sum_{r=1}^{\bar k}\vartheta_r^{-1}v_r(x)v_r(y),
\]
Here $G_+$ inverts $P_g$ on $E^+$ and vanishes on constants and
$E^-$ by \eqref{eq:negative}. For fixed $T$, extract converging
locations and weights in a representation with $k$ slots. The
mixtures converge in $L^s$ to the mixture of the limiting barycenter,
also at collisions or vanishing weights. Elliptic inversion gives
continuity of $\Phi_T$ in $H^2$.

Write $\mathcal G(\nu)=\iint G_+(x,y)\,d\nu(x)d\nu(y)$. The kernel is bounded below everywhere. On a ball of radius $r_T(x_i)$ its diagonal interaction is at least $(8\pi^2)^{-1}\log(1/(2r_T(x_i)))-C$. It follows that
\begin{equation}
 \mathcal G(\nu)\ge \frac{T}{8\pi^2}\sum_i\frac{t_i^2}{\omega(x_i)}-C
 \ge\frac{T}{8\pi^2 b}-C ,
\end{equation}
by Cauchy--Schwarz. This bound uses no separation between support points. Also $d\nu/dm_D\le e^{4T}$, so $\Ent_{m_D}(\nu)=\int\log(d\nu/dm_D)\,d\nu\le4T$. Jensen's inequality under $\nu$ gives $\log Z(v)\ge 4\int v\,d\nu-\Ent_{m_D}(\nu)$. At $v=2\kappa G_+\nu$ we have $Q(v)=4\kappa^2\mathcal G(\nu)$. Consequently
\begin{equation}\label{eq:depth}
 F(\Phi_T(\sigma))\le J_\kappa(v)+C
 \le 4\kappa T-4\kappa^2\mathcal G(\nu)+C
 \le -4\kappa(\rho/b-1)T+C .
\end{equation}

The power weights give $m_D(B(a,s))\le Cs^{4\omega_0}$ uniformly for
$0<s<1$: a mark within $2s$ is treated on $B(q_j,3s)$, and otherwise
its distance is at least $s$ on the ball. Hence
$\nu_{T,x}(B(a,s))\le\min(1,Ce^{4T}s^{4\omega_0})$. Layer-cake
integration of $\log^+(1/d_g)$ bounds $\|G_+\nu_{T,x}\|_\infty$
by $C(1+T)$. Away from $x$ the potential is uniformly bounded by
\eqref{eq:radii}.

Take $T_n\to\infty$ and $\sigma_n\to\sigma$ in $B_\rho$.
Off any neighborhood of $\supp\sigma$, atoms of positive limiting
weight give bounded potentials and those of vanishing weight give
$o(T_n)$ after weighting. But \eqref{eq:depth} and $Q(\Phi_{T_n})\ge0$
give $\log Z(\Phi_{T_n})\ge cT_n-C$. All limiting normalized
densities are therefore supported on $\supp\sigma$. Every limit of
the segments from $\sigma_n$ to $\mu_{\Phi_{T_n}(\sigma_n)}$ belongs
to $B_\rho$, since taking a subset of the support cannot increase its
cost. Compactness puts all such segments in $\mathcal U$ uniformly
for large $T$. Composing with $r$ gives the homotopy, and
\eqref{eq:depth} puts $\Phi_T(B_\rho)$ in $F^{-L}$ for the same $T$
chosen sufficiently large.
\end{proof}

\begin{proof}[Completion of the proof of Theorem~\ref{th:exist}]
Let a finite cycle represent a nonzero reduced homology class of $B_\rho$. Its image in $F^{-L}$ is nonzero by Proposition~\ref{prop:test}. Since $E^+$ is contractible, the image cycle bounds a finite chain in $E^+$, hence in some $F^b$ by compactness of its support. If $F$ had no critical point, Lemma~\ref{lem:lucia} would retract $F^b$ onto $F^{-L}$, a contradiction. Lemma~\ref{lem:red} and \eqref{eq:lift} give the conclusions.
\end{proof}

\medskip
\noindent\begin{minipage}{\linewidth}
\begingroup
\small\raggedright
\noindent\textbf{Mohameden Ahmedou}\\
Department of Mathematics, Justus Liebig University Giessen, Arndtstra\ss e 2, 35392 Gie\ss en, Germany.\\
\href{mailto:Mohameden.Ahmedou@math.uni-giessen.de}{Mohameden.Ahmedou@math.uni-giessen.de}\par\smallskip

\noindent\textbf{Lina Wu}\\
Department of Mathematics, Beijing Jiaotong University, Beijing 100044, China.\\
\href{mailto:lnwu@bjtu.edu.cn}{lnwu@bjtu.edu.cn}\par\smallskip

\noindent\textbf{Lei Zhang}\\
Department of Mathematics, University of Florida, 1400 Stadium Rd, Gainesville, FL 32611, USA.\\
\href{mailto:leizhang@ufl.edu}{leizhang@ufl.edu}\par
\endgroup
\end{minipage}

\end{document}